\documentclass[12p,oneside,a4paper]{scrartcl}
\pdfoutput=1
\def\bf{\normalfont\bfseries}
\setkomafont{sectioning}{\bfseries}
\usepackage{mathptmx}
\usepackage[T1]{fontenc}\usepackage[usenames,dvipsnames]{xcolor}
\usepackage[utf8]{inputenc}
\RequirePackage[OT1]{fontenc}
\usepackage[sans]{dsfont}
\usepackage{amsbsy,bbm,fontawesome,tikzsymbols}
\usepackage{multirow,yhmath, hyperref}
\usepackage{amsmath,amssymb,amsthm}
\allowdisplaybreaks
\usepackage{enumerate,mathtools}
\usepackage{graphics}
\usepackage{natbib}
\makeatletter
\@addtoreset{equation}{section}

\def\N{{\mathbbm {N}}}
\def\R{{\mathbbm {R}}}
\def\Z{{\mathbbm {Z}}}
\def\I{{\mathbbm {1}}}
\makeatother

\newcommand{\ind}{\mathds{1}_}

\theoremstyle{plain}
\newtheorem{theorem}{Theorem}[section]
\newtheorem{rem}[theorem]{Remark}
\newtheorem{lemma}[theorem]{Lemma} 
\newtheorem{assumption}{Assumption}
 
\newtheorem{example}[theorem]{Example} 
\newtheorem{alg}{Algorithm}
 
\theoremstyle{definition}
\newtheorem{defi}[theorem]{Definition} 

\newcommand{\var}{\operatorname{Var}}
\newcommand{\cov}{\operatorname{Cov}}
\newcommand\un{\underline}

\newcommand{\X}{\underline{X}}
\newcommand{\s}{\underline{s}}
\newcommand{\f}{\bar{f}}

\newcommand\be{\begin{equation}}
	\newcommand\ee{\end{equation}}
\newcommand\bee{\begin{equation*}}
	\newcommand\eee{\end{equation*}}
\newcommand\ba{\begin{align}}
	\newcommand\ea{\end{align}}
\newcommand\bae{\begin{align*}}
	\newcommand\eae{\end{align*}}
\newcommand\bi{\begin{itemize}}
	\newcommand\ei{\end{itemize}}
\newcommand\ben{\begin{enumerate}}
	\newcommand\een{\end{enumerate}}
\newcommand\mbf{\mathbf}

\newcommand{\bq}{\begin{eqnarray*}}
\newcommand{\eq}{\end{eqnarray*}}

\newcommand{\Tinfty}{\mathop{\longrightarrow}\limits_{T\to\infty}}
\newcommand{\rull}{\mathop{\longrightarrow}\limits_{r\to 0}}
\newcommand{\dto}{\overset{d}{\longrightarrow}}
\newcommand{\pto}{\overset{P}{\longrightarrow}}

\newcommand{\go}{\mathcal{O}}

\renewcommand{\hat}{\widehat}

\renewcommand{\t}{\mathbf{t}}

\newcommand{\neu}{\color{blue}}

\newcommand{\x}{\mathbf{x}}

\usepackage{mathrsfs,nicefrac}

\usepackage{geometry}
\definecolor{persianblue}{rgb}{0.11, 0.22, 0.73}\definecolor{raspberry}{rgb}{0.89, 0.04, 0.36}

\begin{document}
	
\title{A bootstrap functional central limit theorem for time-varying linear processes}
\author{Carina Beering\thanks{Research group of Mathematics and Statistics, Helmut-Schmidt-Universität/Universität der Bundeswehr Hamburg, Holstenhofweg 85, 22043 Hamburg, Germany} \and Anne Leucht\thanks{Research group of Statistics and Mathematics, Otto-Friedrich-Universität Bamberg, Feldkirchenstr.~21, 96052 Bamberg, Germany}}

\maketitle

\begin{abstract}
We provide a functional central limit theorem for a broad class of smooth functions for possibly noncausal multivariate linear processes with time-varying coefficients. Since the limiting processes depend on unknown quantities, we propose a local block bootstrap procedure to circumvent this inconvenience in practical applications. In particular, we prove bootstrap validity for a very broad class of processes. Our results are illustrated by some numerical examples.
\end{abstract}

 \bigskip\noindent
 \textit{2020 Mathematical Subject Class}: 60F17, 62G09, 62G20


  \noindent  
 \textit{Short title:} A bootstrap functional central limit theorem 
 \medskip

\section{Introduction}\label{sec_1}

Unifying asymptotic theory is a powerful tool to develop statistical test procedures or to quantify the uncertainty of parameter estimators. Still, in many applications the limiting random objects of interest depend on unknown quantities that rely on the data generating process, e.g.~its variance or the underlying dependence structure. Therefore, they cannot be used directly to construct (asymptotically) valid confidence sets or critical values of hypothesis tests. The bootstrap offers a convenient way to overcome these difficulties and is therefore the key to enable practical use of asymptotic results. From a statistical perspective, (functional) central limit theorems and their bootstrap counterparts are particularly appealing as they can be applied to approximate confidence sets for parameters or critical values of $L_2$-test statistics. 
While there exists a large body of literature on bootstrap validity for empirical processes based stationary processes, see e.g.~\citet{Kue89}, \citet{NNR94}, \citet{Bue95}, \citet{DLLN15}, and  \citet{W16}, there are no comparable results for locally stationary processes. The goal of the present paper is to fill this gap as the assumption of a gradually changing probabilistic structure over time is much more realistic than a stationary setting in many applications. Typical examples of use are medical or economical data, see \citet{DPV11}, \citet{AS19}, and \citet{JLMB20}, for instance. The idea of approximating non-stationary time series on segments by stationary ones can be found in \citet{P65}, whereas the concept of local stationarity, on which this paper is based, goes back to \citet{D97}. There, the definition of local stationarity was introduced, and this conception opened the way for momentous theory via the use of an asymptotic (in-fill) framework. An overview of the state of the art is to be found in \citet{D12}. More recently, \citet{DRW19} and \citet{PR22} developed a broad asymptotic theory including laws of large numbers and (functional) central limit theorems for nonlinear, causal locally stationary processes. However, on the bootstrap side methods and validity checks for locally stationary data are mostly tailor-made for specific applications so far: \citet{SP09} used a semiparametric bootstrap in a testing framework, while wild bootstrap methods have been considered for instance by \citet{V15}, \citet{B22} and \citet{KRW22}. \citet{SP08} established a frequency domain bootstrap for statistics of the local periodogram, and \citet{KP15} developed a time-frequency domain bootstrap for a broad class of periodogram based statistics.  In the present paper, we verify that a local version of the block bootstrap can be applied successfully to mimic the distributional behaviour of a large class of empirical processes indexed by smooth functions. This local block bootstrap was initially proposed by   \citet{PP02} for the mean and by  \citet{DPP13} for heteroscedastic time series with trend. Moreover, it was used by \citet{FK22} for indirect estimation in locally stationary structured models. In our case, we investigate possibly two-sided time-varying linear processes. The reasons are two-fold: On the one hand, to the best of our knowledge, there is no functional central limit theorem (FCLT) for noncausal locally stationary linear processes. However, noncausal models play an important role in the context of economic data. Examples for the application of noncausal AR models to stock volume data can be found in \citet{BDT01} and \citet{ACD07}. Moving on to the multivariate case, \citet{LS14} expounded validity verification of prevalent test procedures using the causal VAR model as a base in testing economic hypotheses as an important application of the noncausal VAR model.
On the other hand, a generalization to nonlinear processes would rely on high-level assumptions. Since the technical notation in the proofs would blow up, we restrict ourselves to linear processes here and leave an adaptation of the proof to causal non-linear processes for future research. It is important to note that the method of proof used by \citet{PR22} to verify a FCLT for locally stationary Bernoulli shifts cannot be adapted in a straight-forward manner to prove validity of block bootstrap methods as it relies on martingale difference approaches for Bernoulli shifts. However, the bootstrap process does not inherit this structure.

The rest of the paper is organized as follows: 
In Section~2, we describe the setting under consideration and provide a FCLT for smooth functions of locally stationary linear processes. Following this, Section~3 is devoted to a local block bootstrap procedure for empirical processes. Afterwards, we illustrate the finite-sample performance of our method in Section~4. The proofs of the main results are carried out in Section~5, while several auxiliary results including their proofs are deferred to the Appendix.

\paragraph{Notation}
	Let $|\cdot|_1$ denote the max column sum matrix norm, i.e.~$|M|_1=\max_{1 \leq j \leq r} \sum_{i=1}^{d} |m^{(i,j)}|$ for some $(d\times r)$ matrix $M=(m^{(i,j)})$.
	Note that this matrix norm is submultiplicative. For a $d$-dimensional vector $\un v=(v_1,\ldots,v_d)^\prime$, its $\ell_p$-norm is denoted by $|\un v|_p$,~$p\in[1,\infty]$, whereas the $\mathcal{L}^p$-norm with $ p\in[1,\infty)$ for $d$-dimensional random vectors $\X$ is signified by $\left\|\X\right\|_p:=\left(E\left|\X\right|_p\right)^{1/p}$. Additionally, said notation is transferred to the bootstrap world by defining $\left\|\X^\star\right\|_{p,\star}:=\left(E^\star\left|\X^\star\right|_p\right)^{1/p}$ as the (conditional) bootstrap $\mathcal{L}^p$-norm of $\X^\star$.
	Besides, the Lipschitz seminorm $|\cdot|_{\mathrm{Lip}}$ of a function $\mathsf{f}$ as above signifies
	\begin{equation*}
		|\mathsf{f}|_{\mathrm{Lip}}:=\sup_{\substack{\mathsf{\un x},\mathsf{\un y}\in\mathsf{T}\\\mathsf{\un x}\neq\mathsf{\un y}}}\frac{\left|f\left(\mathsf{\un x}\right)-f\left(\mathsf{\un y}\right)\right|}{\left|\mathsf{\un x}-\mathsf{\un y}\right|_1}.
	\end{equation*}
	For $x\in\R$, let $\lfloor x\rfloor$ be the largest integer smaller as or equal to $x$ and $\lceil x\rceil$ be the smallest integer greater than or equal to $x$.

\section{A functional central limit theorem for locally stationary linear processes}

Let $\left(\un\varepsilon_t\right)_{t\in\Z}$ be a sequence of i.i.d.~centred $\R^d$-valued random vectors and $\left(\X_{t,T}\right)_{t=1}^T$ a $d$-variate (possibly) two-sided linear process
\begin{equation}\label{eq: MAR}
	\X_{t,T}=\un\mu\left(\frac{t}{T}\right)+\sum_{j\in\Z}A_{t,T}(j)\,\un\varepsilon_{t-j},
\end{equation}
where $\un\mu=\left(\mu_1,\dots,\mu_d\right)^\prime$ is a $d$-variate time-varying mean function and $\left(A_{t,T}(j)\right)_{j\in\Z}$ are coefficient matrices of dimension $(d\times d)$. Note that the sequence $\left(\un\varepsilon_t\right)_{t\in\Z}$ is not necessarily a white noise because we do not assume $E\un\varepsilon_t^2<\infty$. Hence, $\left(\X_{t,T}\right)_{t=1}^T$ is only having a MA($\infty$)-representation but is not automatically a MA($\infty$)-process. To ensure that the afore-defined process exists, we need the series in (\ref{eq: MAR}) to converge almost surely, which means the coefficient matrices $\left(A_{t,T}(j)\right)_{j\in\Z}$ have to decay in a sufficiently fast manner as $|j|$ tends to $\infty$. Additionally, we do not allow for rapid changes in the coefficients over time to assure a meaningful statistical methodology.
This leads to the following assumptions concerning the process $\left(\X_{t,T}\right)_{t=1}^T$:
 
\begin{assumption}[Locally stationary linear processes]\label{as: Pro}\ \\
	The process $\left(\X_{t,T}\right)_{t=1}^T$ is of form (\ref{eq: MAR}) with the following specifications:
	\begin{enumerate}[(i)]
		\item The innovations $\left(\un\varepsilon_t\right)_{t\in\Z}$ are i.i.d., centred and $E|\varepsilon_1|_1<\infty$.
		\item For some $\vartheta\in(0,1)$ and a constant $B<\infty$ 
		\begin{equation}\label{eq: Pro1}
			\sup_{t,T} \left| A_{t,T}(j) \right|_1 \,\leq\, B\,\vartheta^j.	
		\end{equation}
		Further, for each $j\in\Z$ there exists an entry-wise continuously differentiable function $A(\cdot,j)\colon[0,1]\to\R^{d\times d}$ such that for all $p,q=1,\dots,d$ and $ s\leq k \in\{0,1\} $, it holds
		\begin{equation}\label{eq: Pro2}
			\sup_{u}\left|\frac{\partial^s a^{(p,q)}(u,j)}{\partial u^s}\right|\,\leq\, B\,\vartheta^j
		\quad\text{and}\quad
			\sup_{t,T} T\left|A_{t,T}(j)-A\left(\frac{t}{T}, j\right)\right|_1\,\leq\, B\,\vartheta^j
		\end{equation}
		with $A(u,j)=\left(a^{(p,q)}(u,j)\right)_{p,q=1,\dots,d}$.
		\item Each component of the mean function $\underline{\mu}$ is continuously differentiable.
	\end{enumerate}
\end{assumption}

 \begin{rem}
		\begin{enumerate}[(i)]
			\item This kind of assumptions represents a classical framework for statistical inference having  estimators pertaining to locally stationary processes as a base, see \citet{D12} and \citet{JLMB20}. Still, many papers, e.g.~\citet{CN10} and \citet{DS06}, only require a polynomial decay instead of a geometric one as in \eqref{eq: Pro1} and \eqref{eq: Pro2}. In fact, a polynomial decay is sufficient in the present context as well. However, the degree of decay depends on the presumed order of absolute moments belonging to the function $f$, which will be introduced later on in a complicated manner, see \citet{B21} for details in a comparable context. For sake of notational simplicity, we stick to the exponential decay here.
			\item Other and more general definitions of local stationarity are invoked for example in \citet{V12} and \citet{DRW19}. They do not require a linear representation of the process $\left(\X_{t,T}\right)_{t=1}^T$ to the price of presupposed causality. 
	\end{enumerate}
\end{rem}

Having introduced nonstationarity as in Assumption~\ref{as: Pro}, the process $\left(\X_{t,T}\right)_{t=1}^T$ can be approximated locally by a (strictly) stationary linear process, its so-called \emph{companion process} $\left(\widetilde{\X}_t(u)\right)_{t\in\Z}$
\begin{equation}\label{eq: DefCP}
 	\widetilde{\X}_t(u)=\un\mu(u)+\sum_{j\in\Z}A(u,j)\,\un\varepsilon_{t-j},
\end{equation}
 as long as $u$ is close to the rescaled time $t/T$.
 
Hence, compared to the original process the mean function $\un\mu$ stays the same for $u=t/T$, whereas the function $A_{t,T}(j)$ is replaced by $A(u,j)$. 
From Assumption~\ref{as: Pro}(ii), we can conclude
\begin{equation}\label{eq.au}
	\sup_{u\in[0,1]}\left|A(u,j)\right|_1\leq\widetilde{B}\,\vartheta^j
\end{equation}
for some finite constant $\widetilde{B}$. This inequality connotes that $(\ref{eq: DefCP})$ possesses a strictly stationary solution for each fixed $u$ while Assumption~\ref{as: Pro} is satisfied. 

\begin{rem}\label{le: AP}
	\begin{enumerate}[(i)]
		\item Closeness of the locally stationary process and its companion process as well as closeness of companion processes for nearby rescaled time points can be specified. More precisely, from \citet[Lemma~2.1]{JLMB20} we obtain 
		\begin{equation*}
			\sup_{1\leq t\leq T}\left\|\X_{t,T}-\widetilde{\X}_t\left(\frac{t}{T}\right)\right\|_1=O(T^{-1})
		\end{equation*}
	if Assumption~\ref{as: Pro} holds true for $k=0$ and 
\begin{equation*}
\left\|\widetilde{\X}_0\left(u_1\right)-\widetilde{\X}_0\left(u_2\right)\right\|_1\leq C\left|u_1-u_2\right|_1\quad\forall u_1,u_2\in[0,1]
\end{equation*}
 for some $C<\infty$ if Assumption~\ref{as: Pro} holds true for $k=1$.
		\item Although the construction with $A_{t,T} (j)$ and $A(t/T, j)$ appears to be unnecessarily complicated, it is required to include time-varying ARMA-processes, see \citet{D12} for details.
	\end{enumerate}
\end{rem}

%
Statistical methods for locally stationary processes can either rely on \emph{local} or \emph{global} characteristics of the process. Local quantities of interest are, for instance, the local variance $\var(\widetilde X_0(u))$ for any fixed $u\in(0,1)$ or the local characteristic function, introduced in \citet{JLMB20} as
$ 
\varphi(u,\underline{s}):=E\left(e^{i\langle\underline{s}\, , \,\widetilde{\underline{X}}_1(u)\rangle}\right), ~\s\in\R^d$.
These quantities can be estimated using kernel estimates based on the observations $X_{1,T},\dots, X_{T,T}$; e.g.~$\varphi(u,\un s)$ can be estimated by the local empirical characteristic function (ECF)
\begin{equation}\label{phi_hat}
\widehat\varphi(u,\underline{s})\,=\,\frac{1}{b_T\,T} \,\sum_{t=1}^T K\left(\frac{t/T-u}{b_T}\right) \, e^{i\langle\underline{s} \, , \, \underline{X}_{t,T}\rangle},	 
\quad \underline{s}\in\R^d,
\end{equation}
using a suitable kernel function $K$ and an appropriate bandwidth $b_T$.
 A simple and prominent example for global quantity of interest is the integrated volatility in high-frequency finance, see \citet[Section~2]{F15} and references therein. Consider a discrete-time model for the intraday log returns $X_{t,T}=\sqrt{T}^{-1}\sigma(t/T)\,\varepsilon_t$ for a smooth deterministic spot volatility function $\sigma$ and a  centred stationary process $(\varepsilon_t)_t$. Then, a natural estimator for the integrated volatility is given by the realized volatility
\begin{equation}\label{eq: rv}
 RV\,=\,\sum_{t=1}^T\, X_{t,T}^2.
\end{equation}
Our goal is to derive a FCLT that is flexible enough to deduce the asymptotic distribution of both quantities \eqref{phi_hat} and \eqref{eq: rv}. Therefore, we consider function classes changing with~$T$
\begin{equation*}
	\mathcal F_T:=\biggl\{	\sum_{t=1}^{T}w_{t,T}\,f\left(\s,\X_{t,T}\right)\colon s\in\mathcal S\,\bigg|\, (w_{t,T})_t \text{ satisfying Ass.~\ref{as: We}},~f \text{ satisfying Ass.~{\neu~\ref{as: Fu2}}} \biggr\}
\end{equation*}
with $\mathcal S\subseteq\overline\R^d$ and specify the underlying Assumptions~\ref{as: We} and~\ref{as: Fu2} below. 
Obviously, the local ECF is included considering its real and imaginary part separately by setting  $w_{t,T}:=\left(b_TT\right)^{-1/2}K\left(\frac{t/T-u}{b_T}\right)$ and $f\left(\s,\X_{t,T}\right)=\cos(\s'\X_{t,T})$ or $f\left(\s,\X_{t,T}\right)=\sin(\s'\X_{t,T})$, respectively. Secondly, if we choose $w_{t,T}:=T^{-1/2}$ and $f(\s,\X_{t,T})=X_{t,T}^2$, we end up with $RV$. The first example illustrates perfectly that several weights may be zero. To obtain asymptotic normality, we control for the number of  zero weights and the magnitude of the non-zero weights as follows:

\begin{assumption}[Weights]\label{as: We}\ \\
	The sequence of non-negative weights $\left(w_{t,T}\right)_{t=1}^T$ fulfils  $\sup_{1\leq t\leq T} w_{t,T}\leq C_w\,d_T^{-1/2}$ for some finite constant $C_w$,
	where $d_T\Tinfty \infty$ denotes the number of non-zero weights in  $\left(w_{t,T}\right)_{t=1}^T$.
\end{assumption}

Before we specify the class of functions for our Donsker-type result, we need to introduce some auxiliary quantities to properly handle the dependence structure of the observed and the companion processes within the proofs relying on truncation arguments.
	Considering a truncation parameter $M\in\N$, we set
	\begin{align}\label{eq: Trun1}
		\X_{t,T}^{(M)}&:=\un\mu\left(\frac{t}{T}\right)+\sum_{|j|< M}A_{t,T}(j)\,\un\varepsilon_{t-j}\quad\text{  and }\quad
		\widetilde{\X}_t^{(M)}(u) :=\un\mu(u)+\sum_{|j|<M}A(u,j)\,\un\varepsilon_{t-j}.
	\end{align}
%
%
\begin{assumption}[Function I]\label{as: Fu}\ \\
 Let $(\mathcal S,\rho)$ be a compact semimetric space with $\mathcal S\subseteq \overline\R^d$. The function $f\colon \mathcal S\times\R^d\to \R$ satisfies 
		\begin{equation}\label{eq: Lip}
		 \sup_{\s\in\mathcal S}	\left|f\left(\s,\underline x\right)-f\left(\s,\underline x^\circ\right)\right|\leq C_{Lip}\left|\underline x-\underline x^\circ\right|_1,\quad\underline x,\,\underline x^\circ\in \R^d,
		\end{equation}
		for  some $C_{Lip}<\infty$.
		Additionally, for some $\delta\in(0,1/2)$ it holds
		\begin{equation*}
			\sup_{t\leq  T,\,\s\in\mathcal S,\; M\in\N}\,E[|f(\s,X_{t,T})|^{2+\delta}+
			|f(\s,X_{t,T}^{(M)})|^{2+\delta}]<\infty
		\end{equation*}
	and
	\begin{equation*}
		\sup_{u\in[0,1],\,\s\in\mathcal S,\, M\in\N}\,E[|f(\s,\widetilde X_0(u))|^{2+\delta}+
		|f(\s,\widetilde X_0^{(M)}(u))|^{2+\delta}]<\infty.
	\end{equation*}
\end{assumption}

\begin{rem}\label{le: AF}
	\begin{enumerate}[(i)]
		\item Obviously, assuming validity of Assumption~\ref{as: Fu} additionally to~Assumption~\ref{as: Pro} allows for a generalization of the results stated in Remark~\ref{le: AP}(i) towards
		\begin{equation*}
			\sup_{1\leq t\leq T}\left\|f\left(\s,\X_{t,T}\right)-f\left(\s,\widetilde{\X}_t\left(\frac{t}{T}\right)\right)\right\|_1=O(T^{-1}).
		\end{equation*}
	\item Our assumptions are slightly different compared to those in \citet{PR22}, where they allow for H\"older continuity of $f$ with respect to (w.r.t.)~$x$. We expect that it is possible to relax our assumption in a similar way, which, however, would have an effect on the choice of tuning parameters of the bootstrap procedure in Section~\ref{sec: boot}. For sake of notational simplicity, we stick to Lipschitz continuity here. Moreover, note that we work under weaker moment constraints regarding the data generating process in their case $s=1$.
	\end{enumerate}
\end{rem}

	We abbreviate the centred version of $f$ by  
\begin{equation}\label{eq: Cent}
	\f\left(\s,\cdot\right):=f\left(\s,\cdot\right)-Ef\left(\s,\cdot\right)
\end{equation}
and state a CLT for the finite-dimensional distributions first. Note that this result is sufficient to deduce asymptotic normality of the realized volatility defined in~\eqref{eq: rv}.

\begin{theorem}[Central Limit Theorem]\label{th: CLT}\ \\
	Suppose Assumption~\ref{as: Pro} holds true for $k=1$ and Assumptions~\ref{as: We} as well as~\ref{as: Fu} are valid. Then, for any $J\in\N$ and $\s_j\in\mathcal S$, $j=1,\dots,J$, we have
	\begin{equation*}
		\biggl(\sum_{t=1}^{T}w_{t,T}\,\f\left(\s_j,\X_{t,T}\right),j=1,\dots,J\biggr)\dto\mathcal{N}\left(\un 0,\mbf{V}\right)
	\end{equation*}
	as $T\to\infty$, where $\mathbf{V}:=\left(\mathbf{V}\left(\s_{j_1},\s_{j_2}\right)\right)_{j_1,j_2=1,\dots,J}$ is a $(J\times J)$ covariance matrix with
	\begin{equation*}
		\mathbf{V}\left(\s_{j_1},\s_{j_2}\right)=\sum_{h\in\Z}\lim_{T\to\infty} \sum_{t=1}^{T}w_{t,T}\,w_{t+h,T}\, \cov\left(f\left(\s_{j_1},\widetilde{\X}_0\left(\frac{t}{T}\right)\right),f\left(\s_{j_2},\widetilde{\X}_h\left(\frac{t}{T}\right)\right)\right).
	\end{equation*}
\end{theorem}
\smallskip

We aim at deriving a FCLT. For that purpose, some additional conditions on $f$ are imposed to assure tightness. To this end, let
\begin{equation*}
	D\left(u,\mathcal S,\rho\right)=\max\left\{\#\mathcal T_0\,\middle|\, \mathcal T_0\subseteq \mathcal S,\;\rho\left(\s_1,\s_2\right)> u\;\forall\, \s_1 \neq \s_2\in \mathcal T_0\right\}
\end{equation*}
denote the usual packing number defined e.g.~in Definition 2.2.3 of \citet{VW00}.  
 
\begin{assumption}[Function II]\label{as: Fu2}\ 
	\begin{enumerate}[(i)]
		\item Additionally to Assumption~\ref{as: Fu}, it holds $\left\|f(\cdot,\underline x)\right\|_\infty<\infty$ for any $\underline x\in\R^d$.
 \item It holds $\left|f\left(\s,\x\right)-f\left(\s^\circ,\x\right)\right|\leq g(\x)\, \rho\left(\s,\,\s^\circ\right)$ with some function $g\colon \R^d\to\R_{\geq 0}$ satisfying one of the following conditions:
 	\begin{enumerate}[(a)]
 	\item $E[g^{2+\delta}(Y)]<K$ for some $K<\infty$ with $Y\in\left\{\widetilde X_t(u), \widetilde X_t^{(M)}(u)\right\}_{t\in\Z, M\in\N}$,
 	\item $\|f\|_\infty <\infty$  and $E[g^{1+\delta/2}(Y)]<K$ for some $K<\infty$ with $Y\in\{\widetilde X_t(u), \widetilde X_t^{(M)}(u)\}_{t\in\Z, M\in\N}$.
 \end{enumerate}
\item For any $u>0$, let $
	D\left(u,\mathcal S,\rho\right)\,\leq\, C_D\,(1+u^{-1})^d$. 
\end{enumerate}
\end{assumption}

Note that in Assumption~\ref{as: Fu2}(ii), we distinguish between bounded and unbounded functions $f$. At first glance, boundedness seems much more restrictive as unboundedness, but the different assumptions concerning the moments of the function $g$ open up the field of applications. As an example, consider the ECF case in \citet{JLMB20}, where the function $g$ is equal to the $|\cdot|_1$-norm. Combined with an $\alpha$-stable distribution with $\alpha=1.5$, which we will use in our simulation study later on, we are not able to fulfil case (a) of the second part of Assumption~\ref{as: Fu2} due to the lack of second absolute moments. Especially in finance, the absence of those moments is quite common. Thus, instead of being mostly excluding, the separate handling of bounded functions broadens the scope.

Assumption~\ref{as: Fu2}(iii) holds, for instance, for $(\mathcal S,\rho)= \left([-S,S]^d,\rho\left(\s_1, \s_2\right)=\left|\s_1-\s_2\right|_1\right)$ for any $S\in(0,\infty)$ as well as for $(\mathcal S,\rho)= \left(\overline \R^d,\rho\left(\s_1, \s_2\right)=\sum_{i=1}^d\left|\arctan\left(\s_{1,i}\right)-\arctan\left(\s_{2,i}\right)\right|_1\right)$.

\begin{theorem}[Functional Central Limit Theorem]\label{th: FCLT}\ \\
	Suppose Assumption~\ref{as: Pro} holds true for $k=1$ and Assumptions \ref{as: We} as well as \ref{as: Fu2} are valid. Then,
	\begin{equation*}
		\biggl(\sum_{t=1}^Tw_{t,T}\,\f\left(\s,\X_{t,T}\right)\biggr)_{\s \in \mathcal S}\dto\left(G\left(\s\right)\right)_{\s \in\mathcal S}
	\end{equation*}
	 as $T$ tends to $\infty$, where $\left(G\left(\s\right)\right)_{\s \in \mathcal S}$ is a centred Gaussian process with continuous sample paths with respect to $\rho$ and covariance function $V\left(\s,\s^\circ\right)$ originating in Theorem~\ref{th: CLT}.
\end{theorem}

\section{Locally blockwise bootstrapped empirical processes}\label{sec: boot}

As already pointed out in \citet{PP02}, the classical block bootstrap algorithm for stationary time series has to be modified in the case of locally stationary time series to capture not only the dependence structure but also the time-changing characteristics of the process. More precisely, given $X_{1,T},\dots, X_{T,T}$, a block of a  bootstrap analogue starting at time point $t$ should only consist of a stretch of the original time series with time index close to $t$. This is achieved by the introduction of an additional tuning parameter, the so-called window parameter, that controls for the range of observations a certain bootstrap block can be drawn from. An adaption of the local block bootstrap (LBB) proposed by \citet{DPP13} to the present setting  reads as follows:

\begin{alg}[Bootstrap Algorithm]\label{al: BS}\  
			\begin{enumerate}[(a)]
				\item Consider a blocklength $L_T$ depending on $T$.
				\item Select a window parameter $D_T\in(0,1)$ such that $TD_T\in\mathbb N$.
				\item Generate i.i.d.~integers $k_0,\dots,k_{\left\lfloor T/L_T\right\rfloor-1}$ using a discrete uniform distribution on $\left[-TD_T,TD_T\right]$.
				\item  For $i=0,\dots\lfloor T/L_T \rfloor-1$, define $\X^\star_{1,T},\dots,\X^\star_{T,T}$ by
				\begin{equation*}
					\X^\star_{j+iL_T,T}:=\X_{j+iL_T+k_i,T}\quad\text{for }j=1,\dots,L_T
				\end{equation*}
		 if the resulting set of indices is in $[1,T]$ and use $-k_i$ instead of $k_i$ otherwise.			
				\item Construct the bootstrap estimator by replacing $\X_{t,T}$ with $\X^\star_{t,T}$, that is
				\begin{equation*}
					\sum_{t=1}^{T}w_{t,T}\,f\left(\s,\X^\star_{t,T}\right).
				\end{equation*}
			\end{enumerate}
	\end{alg}
%
	
	\begin{rem}
		\begin{enumerate}[(i)]
			\item The distribution used to generate $k_0,\dots,k_{\left\lfloor T/L_T\right\rfloor-1}$ does not need to assign uniform weights to every choice of $k_i$, see \citet{PP02}. Here, we choose the discrete uniform distribution as it is easy to handle, analogously to \citet{DPP13}. Besides, it matches the choice made for the ordinary moving block bootstrap algorithm designed for stationary processes.
			\item  The case differentiation in part (d) of Algorithm \ref{al: BS} ensures that if a block is in danger of going over the edge, there is a sound way out. By adjusting the sign of $k_i$ for the whole block, the interrelated dependence structure is preserved, and moreover, no observation is used twice in the same block.
		\end{enumerate}
	\end{rem}

 To establish asymptotic validity, we have to modify our assumptions towards more restrictive moment conditions in the case of unbounded~$f$. Especially when it comes to covariance results, finite $(2+\delta)$-th absolute moments of~$f$ are not always sufficient but we require the following:

\begin{assumption}[Function III]\label{as: Ge}\ \\
	  Assumption \ref{as: Pro} for $k=1$ plus Assumptions \ref{as: We} and \ref{as: Fu} are satisfied and 
	  \begin{equation*}
	  \sup_{u\in[0,1],\,\s\in\mathcal S,\, M\in\N}\,E[|f(\s,\widetilde X_0(u))|^{4+\delta}+
	  |f(\s,\widetilde X_0^{(M)}(u))|^{4+\delta}]<\infty.
	  \end{equation*}
\end{assumption}

In contrast, in the case dealing with bounded functions~$f$ there is no need for modifications of the assumptions. 
With the blocklength $L_T$ and the window parameter $D_T$, two new parameters are involved, which need to behave good-naturedly in combination with the number of non-zero weights $d_T$:

\begin{assumption}[Bootstrap Rates]\label{as: Ra}\ \\ 
	For the blocklength $L_T\in\N$ and the window parameter $D_T$ with $TD_T\in\N$, it holds
	\begin{equation*}
		L_T\to\infty,\quad L_T=o\left(d_T^{ \frac{\delta}{2(1+\delta)}}\right),\quad d_T^{\frac{2\delta}{2+\delta}}=\go(TD_T)  
		\quad\text{and}\quad TD_T=\go\left( d_T^{\frac{1}{2+\delta}}\right).
	\end{equation*}
\end{assumption}

\begin{rem}\ 
In particular, these assumptions imply $L_T=o(TD_T)$.
\end{rem}

In the style of (\ref{eq: Cent}), we define
\begin{equation}\label{eq: CentBS}
	\f^\star\left(\s,\X_{t,T}^\star\right):=f\left(\s,\X_{t,T}^\star\right)-E^\star f\left(\s,\X_{t,T}^\star\right)
\end{equation}
for $t=1,\dots,T$ and $\s\in\R^d$ as the bootstrap version of the centred function $\f$, which enables us to state the bootstrap counterpart to Theorem~\ref{th: CLT} in a comprehensive way:

\begin{theorem}[Bootstrap Central Limit Theorem]\label{th: CLTBS}\ \\
	Under Assumptions \ref{as: Ge} and \ref{as: Ra}, it holds for any $J\in\N$ and $\s_j\in\mathcal S$, $j=1,\dots,J$, 
	\begin{equation*}
		\biggl(\sum_{t=1}^{T}w_{t,T}\,\f^\star\left(\s_j,\X_{t,T}^\star\right),\;j=1,\dots,J\biggr)\stackrel{d}{\longrightarrow}Z\sim\mathcal N(0,\mbf{V})
	\end{equation*}
	in probability, where  $\mbf{V}$ is defined in Theorem~\ref{th: CLT}.
	If additionally $\mbf{V}(\s_1,\s_1)>0$, then
	\begin{equation*}
		\sup_{v\in\R}\left|P^\star\biggl(\sum_{t=1}^{T}w_{t,T}\,\f^\star\left(\s_1,\X_{t,T}^\star\right)\leq v\biggr)-\Phi\left(\frac{v}{\mbf{V}(\s_1,\s_1)}\right)\right|\pto 0.
	\end{equation*}
\end{theorem}
\smallskip

Note that this result is in line with Theorem~3.1 in \citet{PP02}, who considered the mean of locally stationary time series $X_{t,T}=\mu+V(t/T)\varepsilon_t$ with $\mu\in\R$, a smooth function $V$, and $\alpha$-mixing stationary innovations $(\varepsilon_t)_t$ satisfying the stronger moment assumption $E|\varepsilon_t|^6<\infty$.
 
For the more general case of a functional CLT, stronger assumptions concerning the function $g$ introduced in Assumption~\ref{as: Fu2} are required, too:

\begin{assumption}[Function IV]\label{as: G-BS}\ \\
	The function $g$ originating from Assumption~\ref{as: Fu2} satisfies one of the following conditions:
		\begin{enumerate}[(a)]	\item $Eg^{4+\delta}(Y )<K$ for some $K<\infty$ with $Y \in\{\widetilde X_t(u), \widetilde X_t^{(M)}(u)\}_{t\in\Z, M\in\N}$ and for some constant $C_{g,a}<\infty$	\begin{equation*}
				\sup_{t\in\{1,\dots,T\}}\left\|g\left(\X_{t,T} \right)-g\left(\widetilde{\X}_t\left(\frac{t}{T}\right)  \right)\right\|_{4+\delta}\leq \frac{C_{g,a}}{T},
			\end{equation*}
		\item $\|f\|_\infty<\infty$, $Eg^{1+\delta/2}(Y )<K$ for some $K<\infty$ with $Y \in\{\widetilde X_t(u), \widetilde X_t^{(M)}(u)\}_{t\in\Z, M\in\N}$ and for some constant $C_{g,b}<\infty$
			\begin{equation*}
			\sup_{t\in\{1,\dots,T\}}\left\|g\left(\X_{t,T} \right)-g\left( \widetilde{\X}_t\left(\frac{t}{T}\right)\right)\right\|_\frac{2+\delta}{2}\leq \frac{C_{g,b}}{T}.
		\end{equation*}
	\end{enumerate}
\end{assumption}

Finally, we use the previously established bootstrap CLT in combination with a tightness result to prove the desired bootstrap FCLT:

\begin{theorem}[Bootstrap Functional Central Limit Theorem]\label{th: FCLTBS}\ \\
	Let Assumptions \ref{as: Fu2} to \ref{as: G-BS} be true. Then, for $\f^\star$ as in (\ref{eq: CentBS}) it holds
	\begin{equation*}
		\biggl(\sum_{t=1}^Tw_{t,T}\,\f^\star\left(\s,\X_{t,T}^\star\right)\biggr)_{\s \in\mathcal S  }\dto\left(G\left(\s\right)\right)_{\s \in \mathcal S}
	\end{equation*}
	in $P$-probability as $T$ tends to $\infty$, where $\left(G\left(\s\right)\right)_{\s \in \mathcal S}$ is defined in Theorem~\ref{th: FCLT}.
\end{theorem}
%

%
 
%

%

%

\section{Numerical results}\label{S4}

We illustrate the finite sample performance of the proposed local bootstrap procedure by two small numerical examples.
	First, coverage of bootstrap-based confidence sets are investigated for the realized volatility introduced in~\eqref{eq: rv} as a global characteristic of the data generating process. Second, we consider ECFs and study the effect different bootstrap window sizes have on coverage results. In both examples, we replicate the simulations $N = 500$ times each with $B = 500$ bootstrap resamplings. The implementations are
	carried out with the aid of the statistical software \texttt{R}; see \citet{Rcran}.

	\begin{example}
		In order to  investigate the performance of the local bootstrap for RV, we revisit a scenario similar to \citet{F15}. 
		We consider  the coverage of symmetric $90\%$ confidence intervals derived by local bootstrap for samples of size $T= 1000$ and $T=2000$ with $\sigma (u)=0.32(u-0.5)^2+0.04,~u\in[0,1]$, reflecting a volatility smile,  and  $\varepsilon_t=0.5\,\eta_{t-1}+\eta_t$.
			Here, $(\eta_t)_t$ is a sequence of i.i.d.~innovations satisfying $\eta_t\sim\mathcal N(0,\, 0.8)$.
  As it can be seen from Figure~\ref{fig: Ex1Plot}, the coverage of the confidence intervals for RV obtained by the local bootstrap is close to the desired level. In particular, the results are  robust w.r.t.~appropriate choices of blocklength and window size. Note that the optimal blocklengths of the ordinary moving block bootstrap in the sense of \citet[Corollary~7.1, $k=2$]{La03} are $L_{1000}=3$ and $L_{2000}=4$ in this example.
  \begin{figure}[h!]
  	\centering
  	\includegraphics[width=7cm]{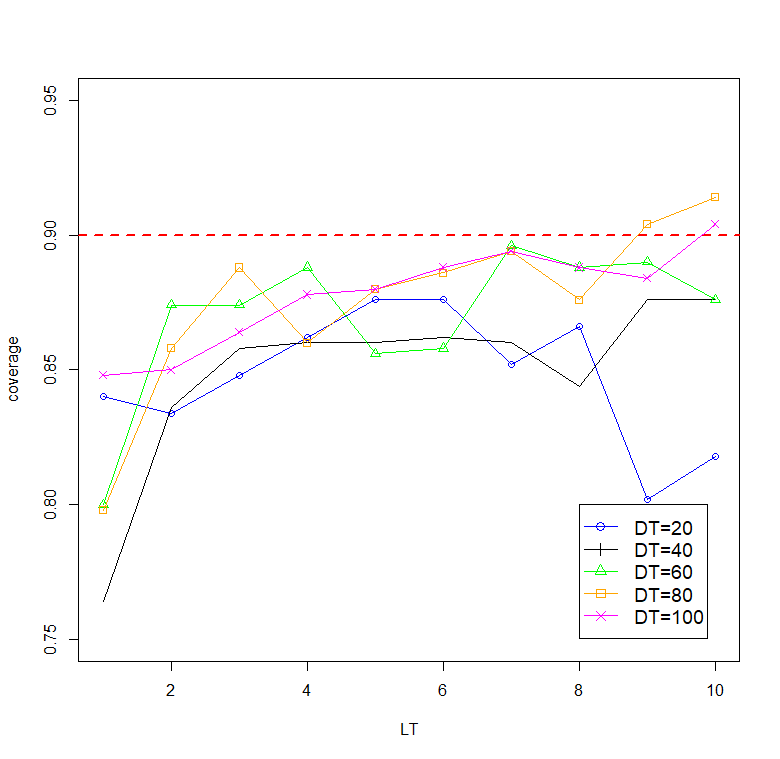}\qquad\includegraphics[width=7cm]{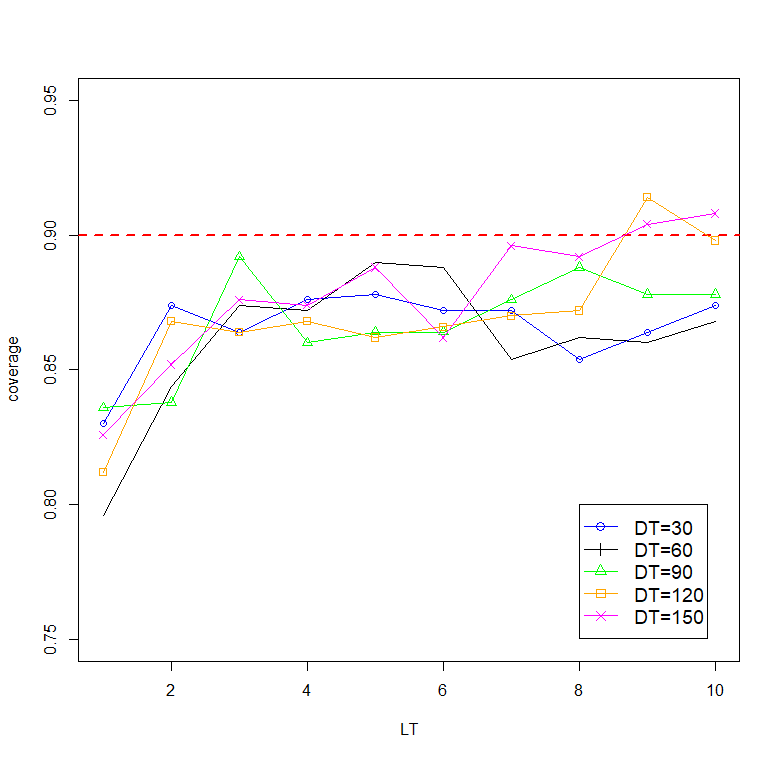}
  	\caption{Coverage results for $T=1000$ (left) and $T=2000$  (right) for different blocklengths $L_T$ and half window sizes  $DT=D_T\,T$.}
  	\label{fig: Ex1Plot}
  \end{figure}
	\end{example}

\begin{example}
	We aim to examine the impact of different sample and window size choices with regard to coverage results in the ECF setup addressed in \citet{JLMB20}. Thus, the quantity in question is $ |\left(b_TT\right)^{1/2}\left(\widehat{\varphi}_X(u;s)-\varphi_{X}(u;s)\right)|$. Similar to \citet[Section 5]{JLMB20}, we generate a locally stationary process $\left(X_{t,T}\right)_{t=1}^T$ by
	\begin{equation*}
		X_{t,T}=\begin{cases}
			0.9\sin\left(2\pi \frac{1}{T}\right) \varepsilon_0+\varepsilon_1, & t=1,\\
			0.9\sin\left(2\pi \frac{t}{T}\right) X_{t-1,T}+\varepsilon_t, & t=2,\dots,T,
		\end{cases}
	\end{equation*}
	with $\left(\varepsilon_t\right)_{t\in\Z}$ forming an i.i.d.~sequence and following an $\alpha$-stable marginal distribution with parameters $\mu=0,\,\alpha=1.5,\,\beta=0$ and $\gamma=0.5$. The innovations bequeaths the $\alpha$-stable distribution to the companion process with slightly different parameters:
	\begin{equation*}
		\widetilde{\mu}(u)=0,\quad\widetilde{\alpha}(u)=1.5,\quad\widetilde{\beta}(u)=0\quad\text{and}\quad \widetilde{\gamma}(u)=\frac{0.5}{1-\left|0.9\sin(2\pi u)\right|^{1.5}}.
	\end{equation*} 
	This leads to 
	\begin{equation*}
		\varphi(u;s)
		=\exp\left(-\frac{0.5|s|^{1.5}}{1-\left|0.9\sin(2\pi u)\right|^{1.5}}\right)
	\end{equation*}
	as belonging characteristic function for the companion process. Remembering (\ref{phi_hat}), we need to specify some other parameters, which are $s=6$, $u=0.4$ and $b_T=T^{-0.4}$. Moreover, we choose the blocklength~$L_T$ equal to $25$. Regarding the sample size, we look at both $T=2000$ and $T=5000$. Because of these choices, there is no need to consider endpoints as they are filtered out by the kernel function. Furthermore, our simulations are based on a significance level of $0.05$. Figure~\ref{fig: Ex2Plot} shows the increase of the coverage results for higher choices of the bootstrap window size towards the aimed $0.95$. While the coverage results regarding $T=5000$ grow faster, the ones belonging to $T=2000$ are closer to the target value for larger window sizes. A wider simulation study to examine parameter choice impact can be found in Chapter~5 in \citet{B21}.
\begin{figure}[h!]
	\centering
	\includegraphics{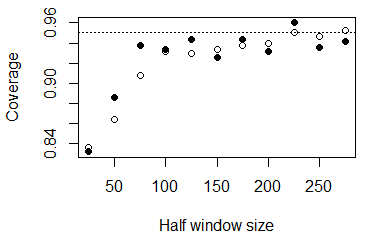}
	\caption{Coverage results for different half window sizes $TD_T$, where $\circ$ denotes the result for $T=2000$ and $\bullet$ the one for $T=5000$, respectively.}
	\label{fig: Ex2Plot}
\end{figure}
\end{example}

For a more involved example in the field of hypothesis testing based on $L_2$ statistics, we refer the reader to Chapter~7 in \citet{B21}, where a characteristic function based test for local independence is established.

\section{Proofs of the main results}\label{S5}
Throughout this section, $C$ denotes a generic constant that may change its value from line to line.
\begin{proof}[Proof of Theorem~\ref{th: CLT}]
	Throughout the proof, we use the following notation based on (\ref{eq: Cent}): 	
	\begin{equation*}
		Y_T\left(\s\right):=\sum_{t=1}^{T} w_{t,T}\,\f\left(\s,\X_{t,T}\right)\quad\text{and}\quad	Y_{t,T}:=\un c^\prime\left(w_{t,T}\,\f\left(\s_j,\X_{t,T}^{(M)}\right),j=1,\dots,J\right)
	\end{equation*}
	for $t=1,\dots,T$ and any $c\in\R^J$. In addition to that, we consider
	\begin{equation*}
		Z_T:=\un c^\prime\left(Y_T\left(\s_j\right),j=1,\dots,J\right)\quad\text{and}\quad Z_T^{(M)}:=\sum_{t=1}^T Y_{t,T}
	\end{equation*}
	and define the corresponding truncated version~${\mathbf V}_M$ of the covariance matrix~${\mathbf V}$ as
	\begin{equation*}
		\mathbf{V}_M:=\left(\mathbf{V}_M\left(\s_{j_1},\s_{j_2}\right)\right)_{j_1,j_2=1,\dots,J}
	\end{equation*}
with
	\begin{equation*}
	\mathbf{V}_M\left(\s_{j_1},\s_{j_2}\right)
	\,:=\,\sum_{h=-2(M-1)}^{2(M-1)}\lim_{T\to\infty} \sum_{t=1}^{T}w_{t,T}\,w_{t+h,T}\, \cov\left(f\left(\s_{j_1},\widetilde{\X}_0^{(M)}\left(\frac{t}{T}\right)\right),f\left(\s_{j_2},\widetilde{\X}_h^{(M)}\left(\frac{t}{T}\right)\right)\right).
\end{equation*}
Using the Cram\'er-Wold device, we are to prove 
	\begin{equation}\label{eq: altCLT}
		Z_T\overset{d}{\longrightarrow}\mathcal N\left(0,\un c^\prime\mathbf{V}\un c\right),\quad\forall c\in\R^J.
	\end{equation}
	Concerning the variance, we distinguish between two cases, namely 
	\begin{equation*}
		\text{(A)}\;\;\un c^\prime\mathbf{V}\un c=0\quad\text{and}\quad \text{(B)}\;\;\un c^\prime\mathbf{V}\un c>0.
	\end{equation*}
	Regarding case (A), we have $
			\var\left(Z_T\right)\longrightarrow \un c^\prime\mathbf{V}\un c=0
	$ by Lemma~\ref{le: Kov}, which gives 
			$Z_T\overset{d}{\longrightarrow}\mathcal N(0,0)$
			as required.
		In the sequel, let $\un c^\prime\mathbf{V}\un c>0$ (case (B)). 
		  In order to show (\ref{eq: altCLT}), Proposition 6.3.9 of \citet{BD91} imposes the verification of the following conditions:
			\begin{enumerate}[(1)]
				\item $\forall M\in\mathbb N:Z_T^{(M)}\overset{d}{\Tinfty} \mathcal N\left(0,\un c^\prime\mathbf{V}_M\un c\right)$,
				\item $\un c^\prime\mathbf{V}_M\un c\underset{M\to\infty}{\longrightarrow}\un c^\prime\mathbf{V}\un c$,
				\item $\forall \epsilon>0:\lim_{M\to\infty}\limsup_{T\to\infty}P\left(\left|Z_T-Z_T^{(M)}\right|\geq \epsilon\right)=0$.
			\end{enumerate}
	
			In terms of constraint (2), it suffices to show
				\begin{equation*}
				\mathbf{V}\left(\s_{j_1},\s_{j_2}\right)-\mathbf{V}_M\left(\s_{j_1},\s_{j_2}\right)\longrightarrow 0, \quad M\to\infty,
			\end{equation*}
			for $j_1,j_2=1,\dots,J$. The difference above can be bounded by
			\begin{align*}
				&\left|\mathbf{V}\left(\s_{j_1},\s_{j_2}\right)-\mathbf{V}_M\left(\s_{j_1},\s_{j_2}\right)\right|\\
				&\begin{multlined}[t][\linewidth]
					\leq\biggl| \sum_{|h|>2(M-1)}\lim_{T\to\infty} \sum_{t=1}^{T}w_{t,T}\,w_{t+h,T}\, \cov\left(f\left(\s_{j_1},\widetilde{\X}_0\left(\frac{t}{T}\right)\right),f\left(\s_{j_2},\widetilde{\X}_h\left(\frac{t}{T}\right)\right)\right)\biggr|\\
					+ \biggl|\sum_{h=-2(M-1)}^{2(M-1)}\lim_{T\to\infty}\sum_{t=1}^{T}w_{t,T}\,w_{t+h,T}\left(\cov\left(f\left(\s_{j_1},\widetilde{\X}_0\left(\frac{t}{T}\right)\right),f\left(\s_{j_2},\widetilde{\X}_h\left(\frac{t}{T}\right)\right)\right)\right.\\
					\left.\left.-\cov\left(f\left(\s_{j_1},\widetilde{\X}_0^{(M)}\left(\frac{t}{T}\right)\right),f\left(\s_{j_2},\widetilde{\X}_h^{(M)}\left(\frac{t}{T}\right)\right)\right)\right)\right|
				\end{multlined}\\
				&=:\text{\un I}+\text{\un{II}}.
			\end{align*}
		 Using  Lemma~\ref{le: KovAb}, we obtain asymptotic negligibility of \un I as $M\to\infty$ under our weight assumptions. Regarding \un{II}, we  obtain
			\begin{align}\label{eq: CLT2}
				&
					\left|\cov\left(f\left(\s_{j_1},\widetilde{\X}_0\left(\frac{t}{T}\right)\right),f\left(\s_{j_2},\widetilde{\X}_h\left(\frac{t}{T}\right)\right)\right)\right.
					-\left.\cov\left(f\left(\s_{j_1},\widetilde{\X}_0^{(M)}\left(\frac{t}{T}\right)\right),f\left(\s_{j_2},\widetilde{\X}_h^{(M)}\left(\frac{t}{T}\right)\right)\right)\right|\nonumber\\
				&\begin{multlined}[t][\linewidth]
					\leq\left|E\left(f\left(\s_{j_1},\widetilde{\X}_0\left(\frac{t}{T}\right)\right)f\left(\s_{j_2},\widetilde{\X}_h\left(\frac{t}{T}\right)\right)\right)
					-E\left(f\left(\s_{j_1},\widetilde{\X}_0^{(M)}\left(\frac{t}{T}\right)\right)f\left(\s_{j_2},\widetilde{\X}_h^{(M)}\left(\frac{t}{T}\right)\right)\right)\right|\\
					+\left|E\left(f\left(\s_{j_1},\widetilde{\X}_0\left(\frac{t}{T}\right)\right)\right)E\left(f\left(\s_{j_2},\widetilde{\X}_h\left(\frac{t}{T}\right)\right)\right)\right.\\
					\left.-E\left(f\left(\s_{j_1},\widetilde{\X}_0^{(M)}\left(\frac{t}{T}\right)\right)\right)E\left(f\left(\s_{j_2},\widetilde{\X}_h^{(M)}\left(\frac{t}{T}\right)\right)\right)\right|
				\end{multlined}\nonumber\\
				&=:\text{\un{II}a}+\text{\un{II}b}.
			\end{align}
			 Starting with the second summand on the right-hand side (RHS), it holds
			\begin{align*}
				\text{\un{II}b}
					%
					%
				%
				%
				%
				&\leq C\, \left\|\widetilde{\X}_0\left(\frac{t}{T}\right)-\widetilde{\X}_0^{(M)}\left(\frac{t}{T}\right)\right\|_1 \leq C\,\vartheta^M
			\end{align*}
			due to Assumption~\ref{as: Fu} and Lemma~\ref{le: AbschM}.
			The remaining summand~\un{II}a of (\ref{eq: CLT2}) can be bounded  with the use of Lemma~\ref{le: DiffProd}
			\begin{align*}
			\text{IIa}&\begin{multlined}[t][0.97\linewidth]
					\leq\left|E\left(\left(f\left(\s_{j_1},\widetilde{\X}_0\left(\frac{t}{T}\right)\right)-f\left(\s_{j_1},\widetilde{\X}_0^{(M)}\left(\frac{t}{T}\right)\right)\right)f\left(\s_{j_2},\widetilde{\X}_h\left(\frac{t}{T}\right)\right)\right.\right.\\
					\left.\left.+f\left(\s_{j_1},\widetilde{\X}_0^{(M)}\left(\frac{t}{T}\right)\right)\left(f\left(\s_{j_2},\widetilde{\X}_h\left(\frac{t}{T}\right)\right)-f\left(\s_{j_2},\widetilde{\X}_h^{(M)}\left(\frac{t}{T}\right)\right)\right)\right)\right|
				\end{multlined}\\
				&\begin{multlined}[t][0.97\linewidth]
					\leq \left|E\left(\left(f\left(\s_{j_1},\widetilde{\X}_0\left(\frac{t}{T}\right)\right)-f\left(\s_{j_1},\widetilde{\X}_0^{(M)}\left(\frac{t}{T}\right)\right)\right)f\left(\s_{j_2},\widetilde{\X}_h\left(\frac{t}{T}\right)\right)\right)\right|\\
					+\left|E\left(f\left(\s_{j_1},\widetilde{\X}_0^{(M)}\left(\frac{t}{T}\right)\right)\left(f\left(\s_{j_2},\widetilde{\X}_h\left(\frac{t}{T}\right)\right)-f\left(\s_{j_2},\widetilde{\X}_h^{(M)}\left(\frac{t}{T}\right)\right)\right)\right)\right|
				\end{multlined}\\
				&\leq C\,\vartheta^{\frac{M\delta}{1+\delta}},
			\end{align*}
		 which yields asymptotic negligibility of \un{II} and, thus, finishes the verification of condition~(2).
			
			 Below, we focus on condition (1). By~(2), it holds $\un c^\prime\mathbf{V}_M\un c>0$ for sufficiently large $M$. Hence, it is adequate to show
			\begin{equation*}
				\frac{Z_T^{(M)}}{\left(\var\left(Z_T^{(M)}\right)\right)^{1/2}}\overset{d}{\longrightarrow}\mathcal N(0,1)
			\end{equation*}
			for $T\to\infty$ because
			\begin{equation}\label{eq: CLT3a}
				\var\left(Z_T^{(M)}\right)=\un c^\prime\mathbf{V}_M\un c+o(1)
			\end{equation}
			can be demonstrated analogously to  Lemma~\ref{le: Kov}. As the number of non-zero weights equals $d_T$, $Z_T^{(M)}$ has only $d_T$ non-vanishing summands denoted by $Y_{t_1,T},\dots,Y_{t_{d_T},T}$. Consequently, we have
			\begin{equation}\label{eq: CLT3b}
				Z_T^{(M)}=\sum_{h=1}^{d_T}Y_{t_h,T},
			\end{equation}
			and $\left(Y_{t_h,T},h=1,\dots,d_T\right)$ form a triangular array of centred $\left(2(M-1)\right)$-dependent random variables such that the CLT in Theorem 2.1 in \citet{RW00} can be applied if the requirements listed therein can be fulfilled (putting their $\gamma=0$). These conditions read for some $\delta>0$ as well as finite constants $\Delta_T,K_{1,T}$ and $K_{2,T}$ depending on $T$ as follows:
			\begin{align*}
				\text{(i)}&\;\;E\left|Y_{t_h,T}\right|^{2+\delta}\leq\Delta_T\quad\forall\,h\in\left\{1,\dots,d_T\right\}, & \text{(ii)}&\;\;\frac{\var\left(\sum_{h=a}^{a+k-1}Y_{t_h,T}\right)}{k}\leq K_{1,T}\quad\forall\,a\,\forall\,k\geq 2(M-1),\\
				\text{(iii)}&\;\;\frac{\var\left(\sum_{h=1}^{d_T}Y_{t_h,T}\right)}{d_T}\geq K_{2,T}, & \text{(iv)}&\;\;\frac{K_{1,T}}{K_{2,T}}=\go(1),\\
				\text{(v)}&\;\;\frac{\Delta_T}{K_{2,T}^{1+\delta/2}}=\go(1), & \text{(vi)}&\;\;\frac{\left(2(M-1)\right)^{2+2/\delta}}{d_T}\Tinfty 0.
			\end{align*}
		
		Now we verify the validity of the conditions stated above starting with~(i):
			\begin{align*}
				E\left|Y_{t_h,T}\right|^{2+\delta}
				&\leq C\,d_T^{-\frac{2+\delta}{2}}\biggl(\sum_{j=1}^J\left|c_j\right|\biggr)^{2+\delta}=:\Delta_T.
			\end{align*}
		 Going on to requirement $\left(\text{ii}\right)$, we make use of the upper bounds for the weights as presupposed in Assumption~\ref{as: We}
			\begin{equation*}
				\var\biggl(\sum_{h=a}^{a+k-1}Y_{t_h,T}\biggr)
				=
				\sum_{h=\max\left\{-(k-1),-2M\right\}}^{\min\{k-1,2M\}}\sum_{s=\max\{1,1-h\}}^{\min\{k,k-h\}}\cov\left(Y_{t_{s+h+a-1},T},Y_{t_{s+a-1},T}\right)
				\leq (2M+1)\,C\,k\,d_T^{-1},
			\end{equation*}
		which proves the validity of (ii) with 
			$K_{1,T}=(2M+1)\, C\,d_T^{-1} $.
		For $\left(\text{iii}\right)$, we obtain from (\ref{eq: CLT3a}),  (\ref{eq: CLT3b}) and (2) that for $T$ large enough there exists some $C_L>0$ such that
			\begin{equation*}
				\var\biggl(\sum_{h=1}^{d_T}Y_{t_h,T}\biggr)=\un c^\prime\mathbf{V}_M\un c+o(1)\geq C_L
			\end{equation*}
			 for $T$ large enough. Hence, (iii) holds with $K_{2,T}= C_L\,d_T^{-1}$, which obviously satisfies (iv) and (v), too.
			Lastly, we see that requirement $(\text{vi})$ holds trivially due to the fact that $M$ is fixed.
			
			Finally, we verify the remaining constraint (3) proving
			\begin{equation*}
			\lim_{M\to\infty}\limsup_{T\to\infty}\var\left(Z_T-Z_T^{(M)}\right)%
			=\lim_{M\to\infty}\limsup_{T\to\infty}\var\biggl(\sum_{t=1}^Tc\,w_{t,T}\left(f\left(\s,\X_{t,T}\right)-f\left(\s,\X_{t,T}^{(M)}\right)\right)\biggr)
			=0.
			\end{equation*}
		We start by inserting another truncated version of the companion process
		\begin{equation*}
			\widetilde{\X}_h^{(M,v)}\left(\frac{t}{T}\right):=\un\mu\left(\frac{t}{T}\right)+\sum_{|j|<\min\{M,v\}} A\left(\frac{t}{T},j\right)\un\varepsilon_{t-j}
		\end{equation*}
	 with $v:=\left\lceil|h|/2\right\rceil$. Using Lemma~\ref{le: KovAb} and similar arguments, we obtain 
		\begin{align*}
			&
				\left|\cov\left(f\left(\s,\widetilde{\X}_0\left(\frac{t}{T}\right)\right)
				-f\left(\s,\widetilde{\X}_0^{(M)}\left(\frac{t}{T}\right)\right),f\left(\s,\widetilde{\X}_h\left(\frac{t}{T}\right)\right)-f\left(\s,\widetilde{\X}_h^{(M)}\left(\frac{t}{T}\right)\right)\right)\right|\\
			&\begin{multlined}[t][\linewidth]
				\leq\left|\cov\left(f\left(\s,\widetilde{\X}_0\left(\frac{t}{T}\right)\right)-f\left(\s,\widetilde{\X}_0^{(M)}\left(\frac{t}{T}\right)\right)\right.\right.\\
				\left.\left.	-f\left(\s,\widetilde{\X}_0^{(v)}\left(\frac{t}{T}\right)\right)-f\left(\s,\widetilde{\X}_0^{(M,v)}\left(\frac{t}{T}\right)\right)\,,\,f\left(\s,\widetilde{\X}_h\left(\frac{t}{T}\right)\right)-f\left(\s,\widetilde{\X}_h^{(M)}\left(\frac{t}{T}\right)\right)\right)\right|	\nonumber\\
				+\left|\cov\left(f\left(\s,\widetilde{\X}_0^{(v)}\left(\frac{t}{T}\right)\right)-f\left(\s,\widetilde{\X}_0^{(M,v)}\left(\frac{t}{T}\right)\right),\right.\right.\\
				\left.\left.		f\left(\s,\widetilde{\X}_h\left(\frac{t}{T}\right)\right)-f\left(\s,\widetilde{\X}_h^{(M)}\left(\frac{t}{T}\right)\right)
				-f\left(\s,\widetilde{\X}_h^{(v)}\left(\frac{t}{T}\right)\right)-f\left(\s,\widetilde{\X}_h^{(M,v)}\left(\frac{t}{T}\right)\right)\right)\right|
			\end{multlined}\\
			&\leq C\,\vartheta^{\frac{v\delta}{1+\delta}}.
		\end{align*}
			This bound is independent of $M$ and can be totalled over $h$. So, Lebesgue's theorem can be used to justify the following result based on the same argumentation as in the proof  of Lemma~\ref{le: Kov} as well as in the proof of Lemma~\ref{le: KovAb}:
			\begin{align*}
				&\lim_{M\to\infty}\limsup_{T\to\infty}\var\biggl(\sum_{t=1}^Tc\,w_{t,T}\left(f\left(\s,\X_{t,T}\right)-f\left(\s,\X_{t,T}^{(M)}\right)\right)\biggr)\\
				%
				%
				%
				%
				%
				%
				&\begin{multlined}[t][\linewidth]
					\leq c^2\sum_{h=-\infty}^\infty\lim_{M\to\infty}\limsup_{T\to\infty}\sum_{t=1}^Tw_{t,T}\,w_{t+h,T}\left|\cov\left(f\left(\s,\widetilde{\X}_0\left(\frac{t}{T}\right)\right)-f\left(\s,\widetilde{\X}_0^{(M)}\left(\frac{t}{T}\right)\right),\right.\right.\\
					\left.\left.f\left(\s,\widetilde{\X}_h\left(\frac{t}{T}\right)\right)-f\left(\s,\widetilde{\X}_h^{(M)}\left(\frac{t}{T}\right)\right)\right)\right|
				\end{multlined}\nonumber\\	&=0,
			\end{align*}
		 which concludes (c) and, hence, finishes the proof.	 
\end{proof}\medskip

\begin{proof}[Proof of Theorem~\ref{th: FCLT}]
	Following Theorem 1.5.4 of \citet{VW00}, we need to show convergence of the fidis and asymptotical tightness in order to prove process convergence. Finally, the continuity of the sample path of the limiting process can be concluded with the help of Addendum 1.5.8 of \citet{VW00}. Theorem~\ref{th: CLT} gives the required convergence of the fidis. Using Theorem~1.5.7 of \citet{VW00}, we show uniform equicontinuity. In view of 
		\begin{equation*}
			\lim_{T\to\infty}P\biggl(\sup_{\s\in\mathcal S}\biggl|\sum_{t=1}^{T} w_{t,T}\, \f\left(\s,\X_{t,T}\right)\biggr|>\lambda\biggr)
			=\lim_{T\to\infty} P\biggl(\sup_{\s\in\mathcal S}\biggl|\sum_{t=1}^{T} w_{t,T}\,\f\left(\s,\widetilde{\X}_t\left(\frac{t}{T}\right)\right)\biggr|>\frac{\lambda}{2}\biggr),
		\end{equation*}
which follows straightforwardly from
\begin{equation*}
	E\biggl(\sup_{\s\in\mathcal S}\,\biggl|\sum_{t=1}^{T}w_{t,T}\left(f\left(\s,\X_{t,T}\right)-f\left(\s,\widetilde{\X}_t\left(\frac{t}{T}\right)\right)\right)\biggr|\biggr)\leq C \,  \sum_{t=1}^{T}w_{t,T}E\left|\X_{t,T}-\widetilde{\X}_t\left(\frac{t}{T}\right)\right|_1 
	\leq C \,d_T^{-1/2},
\end{equation*}
it remains to show
\begin{equation}\label{eq: T}
	\underset{r\to 0}{\lim}\,\limsup_{T\to\infty}\, P\biggl(\underset{ \rho\left(\s_1,\s_2\right)<r}{\sup}\,\biggl|\sum_{t=1}^{T} w_{t,T} \left(\f\left(\s_1,\widetilde{\X}_t\left(\frac{t}{T}\right)\right)-\f\left(\s_2,\widetilde{\X}_t\left(\frac{t}{T}\right)\right)\right)\biggr|>\lambda\biggr)=0
\end{equation}
for any $\lambda>0$.  For this purpose, we define
\begin{equation*}\label{eq: defan}
	\kappa_{T}:= \left\lfloor d_T^{\frac{1}{2m}}\right\rfloor
	\quad\text{and}\quad
	\mu_{T}:= \left\lfloor\frac{d_T}{2\kappa_{T}}\right\rfloor
\end{equation*}
for some case-specific $m>1$ which will be particularized for the cases~(a) and~(b) in Assumption~\ref{as: Fu2} later on. Recall that $d_T$ denotes the number of positive weights, but the non-vanishing weights need not to be subsequent. In the style of \citet{AY94}, we divide our set of indices into blocks $H_t,T_t$ and $R$ in such a way that the indices of the first $\kappa_T$ non-negative weights are in $H_1$, the indices of the second $\kappa_T$ non-negative weights in $T_1$, the indices of the second $\kappa_T$ non-negative weights in $H_2$ and so on until we have eventually $\mu_T$ $H$-blocks and $\mu_T$ $T$-blocks each. The remaining indices are arranged in block~$R$. 
We establish an upper bound for the RHS of (\ref{eq: T}) considering the $H$-blocks, $T$-blocks and the $R$-block separately.
Regarding the last one, we obtain from Assumption~\ref{as: We} and the Lipschitz condition in Assumption~\ref{as: Fu2}
\begin{align*}
	&\underset{r\to 0}{\lim}\,\limsup_{T\to\infty}\, P\biggl(\underset{\substack{  \rho\left(\s_1,\s_2\right)<r}}{\sup} \sum_{i\in R} w_{i,T} \left|\f\left(\s_1,\widetilde{\X}_i\left(\frac{i}{T}\right)\right)-\f\left(\s_2,\widetilde{\X}_i\left(\frac{i}{T}\right)\right)\right|>\frac{\lambda}{3}\biggr)\\
	&\leq\lim_{r\to 0}\,\limsup_{T\to\infty}\, C\underset{\substack{  \rho\left(\s_1,\s_2\right)<r}}{\sup} \,\sum_{i\in R} w_{i,T}\, E\left[g\left( \widetilde{\X}_i\left(\frac{i}{T}\right)\right)\right]\,\rho\left(\s_1,\s_2\right)\\
	&\leq\limsup_{T\to\infty}\, C\,d_T^{-1/2}\kappa_T\\
	&=0.
\end{align*}
As the sum of the $T$-blocks can be treated analogously to the one containing the $H$-blocks, we focus on the latter. In the following, we want to make use of the block structure in such a way that the involved random variables whose indices are situated in different blocks $H_1,\dots,H_{\mu_T}$ are independent. To achieve this, we make use of the truncated variables $\left(\widetilde{\X}_i^{(M)}\left(\frac{i}{T}\right)\right)$ with $M=\lceil\kappa_{T}/2\rceil$ and divide the sum as follows:
\begin{align*}
	&\sum_{t=1}^{\mu_{T}}\sum_{i\in H_t} w_{i,T}\left(\f\left(\s_1,\widetilde{\X}_i\left(\frac{i}{T}\right)\right)-\f\left(\s_2,\widetilde{\X}_i\left(\frac{i}{T}\right)\right)\right)\nonumber\\
	&\begin{multlined}[t][\linewidth]
		=\sum_{t=1}^{\mu_{T}}\sum_{i\in H_t} w_{i,T}\left(\f\left(\s_1,\widetilde{\X}_i\left(\frac{i}{T}\right)\right)-\f\left(\s_1,\widetilde{\X}_i^{(M)}\left(\frac{i}{T}\right)\right)\right)\\
		+\sum_{t=1}^{\mu_{T}}\sum_{i\in H_t} w_{i,T}\,\left( \f\left(\s_2,\widetilde{\X}_i\left(\frac{i}{T}\right)\right)-\f\left(\s_2,\widetilde{\X}_i^{(M)}\left(\frac{i}{T}\right)\right)\right)\\
		+\sum_{t=1}^{\mu_{T}}\sum_{i\in H_t} w_{i,T}\,\left( \f\left(\s_1,\widetilde{\X}_i^{(M)}\left(\frac{i}{T}\right)\right)-\f\left(\s_2,\widetilde{\X}_i^{(M)}\left(\frac{i}{T}\right)\right)\right)
	\end{multlined}\\
	&=:\text{\underline{I}a}+\text{\underline{I}b}+\text{\underline{I}c}.
\end{align*}
First, similar arguments as used in the proof of Lemma~\ref{le: AbschM} yield for \un Ia (and similarly \un Ib)
\begin{align*}
	\underset{r\to 0}{\lim}\,\limsup_{T\to\infty}\, E\biggl(\underset{\substack{  \rho\left(\s_1,\s_2\right)<r}}{\sup}\sum_{t=1}^{\mu_{T}}\sum_{i\in H_t} w_{i,T} \left|\f\left(\s_1,\widetilde{\X}_i\left(\frac{i}{T}\right)\right)-\f\left(\s_1,\widetilde{\X}_i^{(M)}\left(\frac{i}{T}\right)\right)\right|\biggr)
	&\leq C\,\limsup_{T\to\infty}\, \sum_{t=1}^{T} w_{t,T}\,\vartheta^{\lceil\kappa_T/2\rceil}\\
	&=0.
\end{align*}
Hence, it remains to show asymptotic negligibility of
\begin{equation}\label{eq: FCLT2}
	\underset{r\to 0}{\lim}\,\limsup_{T\to\infty}\, P\biggl(\underset{\substack{  \rho\left(\s_1,\s_2\right)<r}}{\sup}\biggl|\sum_{t=1}^{\mu_{T}}\sum_{i\in H_t} w_{i,T}\left(\f\left(\s_1,\widetilde{\X}_i^{(M)}\left(\frac{i}{T}\right)\right)-\f\left(\s_2,\widetilde{\X}_i^{(M)}\left(\frac{i}{T}\right)\right)\right)\biggr|>\frac{\lambda}{9}\biggr).
\end{equation}
Since we only deal with the truncated version of process now, we obtained independence of the summands with different indices~$t$. This opens the way to the use of standard empirical process theory. Before pursuing the proof, we introduce some further notation. Consider $\s,\s_1,\s_2\in\mathcal S$ and define 
\begin{equation*}
	\nu_{T}(\s):= \sum_{t=1}^{\mu_{T}}\sum_{i\in H_t} w_{i,T}\, \f\left(\s,\widetilde{\X}_i^{(M)}\left(\frac{i}{T}\right)\right)
	\quad\text{and}\quad
	\nu_{T}\left(\s_1,\s_2\right):= \nu_{T}\left(\s_1\right) - \nu_{T}\left(\s_2\right),
\end{equation*}
respectively. We follow the main ideas of \citet{AY94} and use a classical chaining argument. For this purpose, let
\begin{equation}\label{eq: Defrk}
	r_k:= r\,2^{-k},\quad k=0,\dots,k_{T},
\end{equation}
for some $r$, $k_T$ which will be specified thereinafter. Moreover, let $\mathcal F_k\subseteq\mathcal S$   be an index set satisfying
\begin{equation*}
	\#\mathcal F_k=D(k)=D\left(r_k,\mathcal S,\rho\right)
	\quad\text{and}\quad
	\underset{\s_1\in\mathcal S}{\sup}\,\underset{\s_2\in\mathcal F_k}{\min}\, \rho\left(\s_1,\s_2\right)<r_k,\qquad k\in\left\{0,\dots,k_T\right\}.
\end{equation*}
By Assumption~\ref{as: Fu2}(iii), it holds $D(k)\leq r_k^{-d}$ for $r>0$ chosen sufficiently small. This gives us the existence of maps $\pi_k\colon \mathcal S \to \mathcal F_k$ for $k=0,\dots,k_T$ such that
\begin{equation*}
	\left|\s-\pi_k\,\s\right|_1\leq r_k\quad\forall \s\in\mathcal S.
\end{equation*}
Subsequently, we get the following two inequalities for $\s,\,\s_1,\s_2\in\mathcal S$ with $\rho\left(\s_1,\s_2\right)<r$:
\begin{equation*}
	\rho\left(\pi_0\,\s_1,\pi_0\,\s_2\right)\,\leq\,3\,r\quad\text{and}\quad	\rho\left(\pi_k\,\s,\pi_{k-1}\,\s\right)\leq 3\,r_k, \quad k\in\left\{1,\dots,k_T\right\}.
\end{equation*}
Thus, we get
$$
	\underset{\substack{\s_1,\s_2\in\mathcal S\\ \rho\left(\s_1,\s_2\right)<r}}{\sup}\left|\nu_T\left(\s_1,\s_2\right)\right| 
	%
	%
	%
	%
	%
	%
	\leq 2 \underset{\substack{\s_1,\s_2\in\mathcal S\\ \rho\left(\s_1,\s_2\right)\leq r_{k_T}}}{\sup} \left|\nu_T\left(\s_1,\s_2\right)\right|
	+\underset{\substack{\s_1,\s_2\in\mathcal F_0\\ \rho\left(\s_1,\s_2\right)\leq 3\,r}}{\sup}\left|\nu_T\left(\s_1,\s_2\right)\right|
	+ 2 \sum_{k=1}^{k_T} \underset{\substack{\s_1\in\mathcal F_k,\s_2\in\mathcal F_{k-1}\\ \rho\left(\s_1,\s_2\right)\leq 3\,r_k}}{\sup} \left|\nu_T\left(\s_1,\s_2\right)\right|.
$$
Again, we introduce some auxiliary quantities. Let
\begin{equation*}
	\lambda_k:= r_k^{\frac{1}{4(1+\delta)}} \vee \left(\frac{4}{\bar{C}}\,r_k^{\frac{1}{1+\delta}}\log D(k)\right)^{1/2}, \qquad k\in\left\{1,\dots,k_T\right\},
\end{equation*}
be defined for some finite constant $\bar{C}>0$, which will be specified later on and may take different values in cases (a) and (b). Hence, we get
%
%
\begin{equation}\label{eq: logDk}
	\log D(k)\leq\lambda_k^2\,\frac{\bar{C}}{4}\,r_k^{-\frac{1}{1+\delta}}.
\end{equation}
Additionally, let $r$ be small enough to allow for
\begin{equation}\label{eq: Sumlambdak}
	4\sum_{k\in\N}\lambda_k\leq \frac{\lambda}{27}.
\end{equation}
Since we have $D(k)=\mathcal O\left(r_k^{-d}\right)$, summability of $\left(\lambda_k\right)_{k=1}^{k_T}$ for $T\to\infty$ is assured. At this point, we come back to (\ref{eq: FCLT2}). With the preassigned notation and~(\ref{eq: Sumlambdak}), we can split up as follows:
\begin{align}\label{eq: EZD}
	&P\Biggl(\underset{\substack{\s_1,\s_2\in\mathcal S\\ \rho\left(\s_1,\s_2\right)<r}}{\sup}\sum_{t=1}^{\mu_{T}}\sum_{i\in H_t} w_{i,T} \left|\f\left(\s_1,\widetilde{\X}_i^{(M)}\left(\frac{i}{T}\right)\right)-\f\left(\s_2,\widetilde{\X}_i^{(M)}\left(\frac{i}{T}\right)\right)\right|>\frac{\lambda}{9}\Biggr)\nonumber\\
	%
	%
	%
	%
	%
	%
	&\begin{multlined}[t][\linewidth]
		\leq P\Biggl(2 \underset{\substack{\s_1,\s_2\in\mathcal S\\ \rho\left(\s_1,\s_2\right)\leq r_{k_T}}}{\sup} \left|\nu_T\left(\s_1,\s_2\right)\right|>\frac{\lambda}{27}\Biggr)\\
		+P\Biggl(2 \sum_{k=1}^{k_T} \underset{\substack{\s_1\in\mathcal F_k,\s_2\in\mathcal F_{k-1}\\ \rho\left(\s_1,\s_2\right)\leq 3\,r_k}}{\sup} \left|\nu_T\left(\s_1,\s_2\right)\right|>4\sum_{k=1}^{k_T}\lambda_k\Biggr)
		+P\Biggl(\underset{\substack{\s_1,\s_2\in\mathcal F_0\\ \rho\left(\s_1,\s_2\right)\leq 3\,r}}{\sup}\left|\nu_T\left(\s_1,\s_2\right)\right|>\frac{\lambda}{27}\Biggr)
	\end{multlined}\nonumber\\
	&=:\text{I}+\text{II}+\text{III}.
\end{align}
In the following, we treat the individual terms in two different ways. In order to show asymptotic negligibility of terms II and III, we want to make use of Bernstein's inequality for sums of independent random variables exerted on the outer sum of $\nu_T$. Term I, however, will be discussed by using a symmetrization lemma at the end of the proof.

The remaining part of the proof presumes Assumption~\ref{as: Fu2}(a) to hold. For (b) see Lemma~\ref{L: rest-FCLT} in the Appendix. 
Before starting with the examination of term II in (\ref{eq: EZD}), we specify the lower bound of $m$  as $	m>\frac{20+15\delta-4\delta^2-3\delta^3}{2\delta\left(1-\delta^2\right)}$.
 Moreover, we need $r_{k_T}$ in (\ref{eq: Defrk}) to meet the following bounding condition:
	\begin{equation}\label{eq: defkn}
		d_T^{-\frac{(2m\delta-4-3\delta)(1+\delta)}{2 m(4+3\delta)}}\leq r_{k_{T}}\leq d_T^{-\frac{2}{m(1-\delta)}}.
	\end{equation}
	Note that our choice of $m$ guaranties that the left-hand side (LHS) is strictly smaller than the RHS.		
	%
	%
	%
	%
	%
	%
	%
	%
	%
	%
	%
	%
	Now we turn our attention to the second summand in (\ref{eq: EZD}). To be able to apply Bernstein's inequality, we need to establish an upper bound for the variance of the inner sum of $\nu_T$. Consider $l:=\left|i_1-i_2\right|$. Then, we have
	\begin{equation}\label{eq: FCLT3}
		\var\left(\nu_T\left(\s_1,\s_2\right)\right)
		%
		%
		%
		%
		%
		%
		%
		%
		%
		%
		%
		%
		\begin{multlined}[t][0.85\linewidth]
			\leq\sum_{t=1}^{\mu_{T}}\sum_{i_1,i_2\in H_t}w_{i_1,T}\,w_{i_2,T}\left|\cov\left(\left(f\left(\s_1,\widetilde{\X}_{i_1}^{(M)}\left(\frac{i_1}{T}\right)\right)-f\left(\s_2,\widetilde{\X}_{i_1}^{(M)}\left(\frac{i_1}{T}\right)\right)\right)\right.\right.\\
			-\left(f\left(\s_1,\widetilde{\X}_{i_1}^{\left(M(l)\right)}\left(\frac{i_1}{T}\right)\right)-f\left(\s_2,\widetilde{\X}_{i_1}^{\left(M(l)\right)}\left(\frac{i_1}{T}\right)\right)\right),\\
			\left.\left.f\left(\s_1,\widetilde{\X}_{i_2}^{(M)}\left(\frac{i_2}{T}\right)\right)-f\left(\s_2,\widetilde{\X}_{i_2}^{(M)}\left(\frac{i_2}{T}\right)\right)\right)\right|\\
			+\sum_{t=1}^{\mu_{T}}\sum_{i_1,i_2\in H_t}w_{i_1,T}\,w_{i_2,T}\left|\cov\left(f\left(\s_1,\widetilde{\X}_{i_1}^{\left(M(l)\right)}\left(\frac{i_1}{T}\right)\right)-f\left(\s_2,\widetilde{\X}_{i_1}^{\left(M(l)\right)}\left(\frac{i_1}{T}\right)\right),\right.\right.\\
			\left(f\left(\s_1,\widetilde{\X}_{i_2}^{(M)}\left(\frac{i_2}{T}\right)\right)-f\left(\s_2,\widetilde{\X}_{i_2}^{(M)}\left(\frac{i_2}{T}\right)\right)\right)\\
			\left.\left.-\left(f\left(\s_1,\widetilde{\X}_{i_2}^{\left(M(l)\right)}\left(\frac{i_2}{T}\right)\right)-f\left(\s_2,\widetilde{\X}_{i_2}^{\left(M(l)\right)}\left(\frac{i_2}{T}\right)\right)\right)\right)\right|
		\end{multlined}
	\end{equation}
	for $M(l):=\left\lceil\min\{M,l/2\}\right\rceil$ as truncation parameter. Next, we take a closer look at only the first covariance of (\ref{eq: FCLT3}) since the second one behaves similarly. We have
	\begin{align}\label{eq: FCLT4}
		&\begin{multlined}[b][\linewidth]
			\left|\cov\left(\left(f\left(\s_1,\widetilde{\X}_{i_1}^{(M)}\left(\frac{i_1}{T}\right)\right)-f\left(\s_2,\widetilde{\X}_{i_1}^{(M)}\left(\frac{i_1}{T}\right)\right)\right)\right.\right. 
			-\left(f\left(\s_1,\widetilde{\X}_{i_1}^{\left(M(l)\right)}\left(\frac{i_1}{T}\right)\right)-f\left(\s_2,\widetilde{\X}_{i_1}^{\left(M(l)\right)}\left(\frac{i_1}{T}\right)\right)\right)\,,\\
			\left.\left.f\left(\s_1,\widetilde{\X}_{i_2}^{(M)}\left(\frac{i_2}{T}\right)\right)-f\left(\s_2,\widetilde{\X}_{i_2}^{(M)}\left(\frac{i_2}{T}\right)\right)\right)\right|
		\end{multlined}\nonumber\\
		&\begin{multlined}[t][\linewidth]
			\leq \left|\cov\left(f\left(\s_1,\widetilde{\X}_{i_1}^{(M)}\left(\frac{i_1}{T}\right)\right)-f\left(\s_1,\widetilde{\X}_{i_1}^{\left(M(l)\right)}\left(\frac{i_1}{T}\right)\right)\,,\, f\left(\s_1,\widetilde{\X}_{i_2}^{(M)}\left(\frac{i_2}{T}\right)\right)-f\left(\s_2,\widetilde{\X}_{i_2}^{(M)}\left(\frac{i_2}{T}\right)\right)\right)\right|\\
			+\left|\cov\left(f\left(\s_2,\widetilde{\X}_{i_1}^{(M)}\left(\frac{i_1}{T}\right)\right)-f\left(\s_2,\widetilde{\X}_{i_1}^{\left(M(l)\right)}\left(\frac{i_1}{T}\right)\right)\,,\,f\left(\s_1,\widetilde{\X}_{i_2}^{(M)}\left(\frac{i_2}{T}\right)\right)-f\left(\s_2,\widetilde{\X}_{i_2}^{(M)}\left(\frac{i_2}{T}\right)\right)\right)\right|,
		\end{multlined}
	\end{align}
	and, again, we only examine the first covariance of (\ref{eq: FCLT4}) due to the same reason. Invoking Lipschitz continuity of $f$, we obtain similarly to Lemma~\ref{le: AbschM}
		\begin{align}\label{eq: FCLT5}
			&\left|\cov\left(f\left(\s_1,\widetilde{\X}_{i_1}^{(M)}\left(\frac{i_1}{T}\right)\right)-f\left(\s_1,\widetilde{\X}_{i_1}^{\left(M(l)\right)}\left(\frac{i_1}{T}\right)\right),\,f\left(\s_1,\widetilde{\X}_{i_2}^{(M)}\left(\frac{i_2}{T}\right)\right)-f\left(\s_2,\widetilde{\X}_{i_2}^{(M)}\left(\frac{i_2}{T}\right)\right)\right)\right|\nonumber\\
			&
				\leq E\left|\left(f\left(\s_1,\widetilde{\X}_{i_1}^{(M)}\left(\frac{i_1}{T}\right)\right)-f\left(\s_1,\widetilde{\X}_{i_1}^{\left(M(l)\right)}\left(\frac{i_1}{T}\right)\right)\right) \left(f\left(\s_1,\widetilde{\X}_{i_2}^{(M)}\left(\frac{i_2}{T}\right)\right)-f\left(\s_2,\widetilde{\X}_{i_2}^{(M)}\left(\frac{i_2}{T}\right)\right)\right)\right|\nonumber\\
				&\quad+C\,  \rho(\s_1,\s_2)\,  \sum_{M>|j|\geq M(l)}\vartheta^j.
		\end{align}
	In the later following calculations to bound the variance of $\nu_T$, we will need two suitable but different bounds. Therefore, we establish two alternative bounds for the first summand of (\ref{eq: FCLT5}). The first will make use of the closeness between the truncated and the two times truncated version of the companion process, whereas the second will consist of the difference between $\s_1$ and $\s_2$.
	\begin{enumerate}[i)]
		\item Using H\"older's inequality, we get  
		\begin{align*}
			&
				E\left|\left(f\left(\s_1,\widetilde{\X}_{i_1}^{(M)}\left(\frac{i_1}{T}\right)\right)-f\left(\s_1,\widetilde{\X}_{i_1}^{\left(M(l)\right)}\left(\frac{i_1}{T}\right)\right)\right)\left(f\left(\s_1,\widetilde{\X}_{i_2}^{(M)}\left(\frac{i_2}{T}\right)\right)-f\left(\s_2,\widetilde{\X}_{i_2}^{(M)}\left(\frac{i_2}{T}\right)\right)\right)\right|
		\\
			&\leq\left\|f\left(\s_1,\widetilde{\X}_{i_1}^{(M)}\left(\frac{i_1}{T}\right)\right)-f\left(\s_1,\widetilde{\X}_{i_1}^{\left(M(l)\right)}\left(\frac{i_1}{T}\right)\right)\right\|_\frac{2+\delta}{1+\delta}\,\sup_{\un s\in\mathcal S}\left\|f\left(\s,\widetilde{\X}_{i_2}^{(M)}\left(\frac{i_2}{T}\right)\right)\right\|_{2+\delta}\\
			&\leq C \left\|\un{\varepsilon}_0\right\|_\frac{2+\delta}{1+\delta}  \sum_{M>|j|\geq M(l)}\vartheta^j .
		\end{align*}
		\item On the other hand, the first summand of (\ref{eq: FCLT5}) can be bounded by	
		\begin{equation*}
			\sup_{M\in\N}\left\|f\left(\s_1,\widetilde{\X}_{i_1}^{(M)}\left(\frac{i_1}{T}\right)\right)\right\|_{2+\delta}\left\|f\left(\s_1,\widetilde{\X}_{i_2}^{(M)}\left(\frac{i_2}{T}\right)\right)-f\left(\s_2,\widetilde{\X}_{i_2}^{(M)}\left(\frac{i_2}{T}\right)\right)\right\|_\frac{2+\delta}{1+\delta}
			\leq C\,\rho\left(\s_1,\s_2\right)
		\end{equation*}
	using similar arguments as in the proof of Lemma~\ref{le: DiffProd}.
	\end{enumerate}
	The combination of these two bounds for (\ref{eq: FCLT5}) and similar arguments for (\ref{eq: FCLT4}) allow us to bound the covariance in (\ref{eq: FCLT3}) via
	\begin{align*}
		&\begin{multlined}[t][\linewidth]
			\left|\cov\left(\left(f\left(\s_1,\widetilde{\X}_{i_1}^{(M)}\left(\frac{i_1}{T}\right)\right)-f\left(\s_2,\widetilde{\X}_{i_1}^{(M)}\left(\frac{i_1}{T}\right)\right)\right)\right.\right.
			-\left(f\left(\s_1,\widetilde{\X}_{i_1}^{\left(M(l)\right)}\left(\frac{i_1}{T}\right)\right)-f\left(\s_2,\widetilde{\X}_{i_1}^{\left(M(l)\right)}\left(\frac{i_1}{T}\right)\right)\right),\\
			\left.\left.f\left(\s_1,\widetilde{\X}_{i_2}^{(M)}\left(\frac{i_2}{T}\right)\right)-f\left(\s_2,\widetilde{\X}_{i_2}^{(M)}\left(\frac{i_2}{T}\right)\right)\right)\right|
		\end{multlined}\\
		%
		%
		%
		&	\leq  C\min\Biggl\{\rho\left(\s_1,\s_2\right),\sum_{M>|j|\geq \left\lceil\frac{\left|i_1-i_2\right|}{2}\right\rceil}\vartheta^j\, \Biggr\}\\
		&	\leq C\, \min\left\{\rho\left(\s_1,\s_2\right),\left|i_1-i_2\right|^{-\frac{1+\delta}{\delta}}\right\}.
	\end{align*}
	%
	%
	%
	%
	Thus, we obtain for any $R_0\geq2$
	\begin{align*}
		\var\left(\nu_T\left(\s_1,\s_2\right)\right)
		&\leq\sum_{t=1}^{\mu_{T}}\sum_{i_1,i_2\in H_t}\,w_{i_1,T}\,w_{i_2,T}\,C\min\left\{\rho\left(\s_1,\s_2\right),\left|i_1-i_2\right|^{-\frac{1+\delta}{\delta}}\right\}\\
		&\leq C\, \mu_T\,\kappa_T\,d_T^{-1}\sum_{t=0}^{d_T}\min\left\{\rho\left(\s_1,\s_2\right),t^{-\frac{1+\delta}{\delta}}\right\}\\
		&\leq C\, \left(R_0\,\rho\left(\s_1,\s_2\right)+R_0^{-1/\delta}\right).
	\end{align*}
	With $R_0:=\left\lfloor\rho\left(\s_1,\s_2\right)^{-\frac{\delta}{1+\delta}}\right\rfloor$ and for any $r$ chosen sufficiently small, we get
	\begin{equation}\label{eq:var-bound}
		\var\left(\nu_T\left(\s_1,\s_2\right)\right)
		\leq C\, \rho\left(\s_1,\s_2\right)^{\frac{1}{1+\delta}}.
	\end{equation}
	Since we aim at the application of Bernstein's inequality to bound II in~\eqref{eq: EZD}, we first provide a suitable approximation of II by a sum of bounded random variables. To this end, we define
	\begin{equation*}
		\breve{f}\left(\s,\widetilde{\X}_t^{(M)}\left(\frac{t}{T}\right)\right) := 
		f\left(\s,\widetilde{\X}_t^{(M)}\left(\frac{t}{T}\right)\right)\I_{\Omega_{\sup,t} }
	\end{equation*}
	and
	\begin{equation*}
		\breve f^c\left(\s,\widetilde{\X}_t^{(M)}\left(\frac{t}{T}\right)\right)=\breve f\left(\s,\widetilde{\X}_t^{(M)}\left(\frac{t}{T}\right)\right)- E\breve f\left(\s,\widetilde{\X}_t^{(M)}\left(\frac{t}{T}\right)\right)
	\end{equation*}
	with 
	\begin{equation*}	\Omega_{\sup,t}:=\biggl\{\omega\in\Omega\,\bigg|\,\underset{ {\s\in\mathcal S,\; i\in H_t}}{\sup}w_{i,T}\left|f\left(\s,\widetilde{\X}_i^{(M)}\left(\frac{i}{T}\right)\right)\right|\leq d_T^{-\frac{\delta}{4+3\delta}}\biggr\}.
	\end{equation*}
	Note that in view of compactness of $(\mathcal S, \rho)$, we have
	\begin{align*}
		P \left(\Omega\backslash\Omega_{\sup,t}\right)
		&\begin{multlined}[t][0.87\linewidth]
			\leq  P\biggl(\underset{ {\s\in\mathcal S,\; i\in H_t}}{\sup}w_{t,T}\left|f\left(\s,\widetilde{\X}_i^{(M)}\left(\frac{i}{T}\right)\right)-f\left(\un{0},\widetilde{\X}_i^{(M)}\left(\frac{i}{T}\right)\right)\right|>\frac{1}{2}\,d_T^{-\frac{\delta}{4+3\delta}}\biggr)\\
			+P\biggl(\underset{ {\s\in\mathcal S,\; i\in H_t}}{\sup}w_{t,T}\left|f\left(\un{0},\widetilde{\X}_i^{(M)}\left(\frac{i}{T}\right)\right)\right|> \frac{1}{2}\,d_T^{-\frac{\delta}{4+3\delta}}\biggr)
		\end{multlined}\\
	&\leq C\,   \sum_{i\in H_t} w_{i,T}^{2+\delta}\,d_T^{ \frac{(2+\delta)\delta}{4+3\delta}}\left[
	E\left[g\left( \widetilde{\X}_i^{(M)}\left(\frac{i}{T}\right)\right)\right]^{2+\delta}+E\left|f\left(\un{0},\widetilde{\X}_i^{(M)}\left(\frac{i}{T}\right)\right)\right|^{2+\delta}\right]\\
	&\leq C\,   d_T^{-\frac{{\delta^2}}{2(4+3\delta)}-1+\frac{1}{2m}}\\
	&=\, o(1).
	\end{align*}
	Hence, II in~\eqref{eq: EZD} can be bounded from above by 
	\begin{align}\label{eq.absch-bern}
			& P\Biggl(2 \sum_{k=1}^{k_T} \underset{\substack{\s_1\in\mathcal F_k,\s_2\in\mathcal F_{k-1}\\ \rho\left(\s_1,\s_2\right)\leq 3\,r_k}}{\sup} \left|\nu_T\left(\s_1,\s_2\right)\right|>4\sum_{k=1}^{k_T}\lambda_k\Biggr)\nonumber\\
			&\begin{multlined}[t][\linewidth]
				\leq  P\Biggl(  \sum_{k=1}^{k_T} \underset{\substack{\s_1\in\mathcal F_k,\s_2\in\mathcal F_{k-1}\\ \rho\left(\s_1,\s_2\right)\leq 3\,r_k}}{\sup} \biggl|\sum_{t=1}^{\mu_T}\sum_{i\in H_t}w_{i,T}\left[\breve{f}^c\left(\s_1,\widetilde{\X}_i^{(M)}\left(\frac{i}{T}\right)\right)-\breve{f}^c\left(\s_2,\widetilde{\X}_i^{(M)}\left(\frac{i}{T}\right)\right)\right] \biggr|> \sum_{k=1}^{k_T}\lambda_k\Biggr)\\
				+  \sum_{k=1}^{k_T} \I_{C \, r_k\,\sum_{t=}^{\mu_T}\sum_{i\in H_t}\,w_{i,T}\, E\left(g\left( \widetilde{\X}_t^{(M)}\left(\frac{i}{T}\right)\right)\,\I_{\Omega\backslash\Omega_{\sup,t}}\right)\geq \frac{\lambda_k}{2}} +o(1).
			\end{multlined} 
	\end{align}
	Note that it holds
	\begin{equation*}
		\sum_{t=1}^{\mu_T}\sum_{i\in H_t}\, w_{i,T}\, E\left[g\left( \widetilde{\X}_t^{(M)}\left(\frac{i}{T}\right)\right)\,\I_{\Omega\backslash\Omega_{\sup,t}}\right] \leq C\, \,d_T^{1/2}\,\left( d_T^{-\frac{{\delta^2}}{2(4+3\delta)}-1+\frac{1}{2m}}\right)^{\frac{1+\delta}{2+\delta}}
	\end{equation*}
	which implies asymptotic negligibility of the middle term on the RHS of \eqref{eq.absch-bern} since
	$(r_k/\lambda_k)_k$ is uniformly bounded. 
	To bound the first summand on the RHS of~\eqref{eq.absch-bern}, we can apply Bernstein's inequality and get
	\begin{equation}\label{eq: FCLT7}
		\begin{aligned}
		%
		%
		%
		%
		%
		&	\sum_{k=1}^{k_T} \,\sum_{\substack{\s_1\in\mathcal F_k,\s_2\in\mathcal F_{k-1}\\ \rho\left(\s_1,\s_2\right)\leq 3\,r_k}}P\biggl(\sum_{t=1}^{\mu_{T}}\biggl|\sum_{i\in H_t} w_{i,T} 
			\left( \breve{f}^c\left(\s_1,\widetilde{\X}_i^{(M)}\left(\frac{i}{T}\right)\right)-\breve{f}^c\left(\s_2,\widetilde{\X}_i^{(M)}\left(\frac{i}{T}\right)\right)\right)\biggr|>\lambda_k\biggr)\\
		&\leq\, 2\sum_{k=1}^{k_T} D(k)\, D(k-1)\,\exp\left(-\frac{1}{2}\cdot\frac{\lambda_k^2}{V_{II,k}+\frac{\breve{M}\lambda_k}{3}}\right),
		\end{aligned}
	\end{equation}
	where, by definition of $H_t$ and $\breve f$,
	\begin{equation*}
		\breve{M}:=4\,d_T^{\frac{4+\delta(3-2m)}{2m(4+3\delta)}}\geq\sup_{\substack{\s_1, \s_2\in\mathcal S\\t\in\{1,\dots,\mu_T\}}}\biggl|\sum_{i\in H_t} w_{i,T}\left( \breve{f}\left(\s_1,\widetilde{\X}_i^{(M)}\left(\frac{i}{T}\right)\right)-\breve{f}\left(\s_2,\widetilde{\X}_i^{(M)}\left(\frac{i}{T}\right)\right)\right)\biggr|
	\end{equation*}
	and, in view of (\ref{eq:var-bound}), for some appropriately chosen $C<\infty$
	\begin{align*}
		V_{II,k} :=C\,r_k^{\frac{1}{1+\delta}} 
		&\geq 	\underset{\substack{\s_1\in\mathcal F_k,\s_2\in\mathcal F_{k-1}\\ \rho\left(\s_1,\s_2\right)\leq 3\,r_k}}{\sup}\var\left(\nu_T\left(\s_1,\s_2\right)\right)\\
		&\geq\underset{\substack{\s_1\in\mathcal F_k,\s_2\in\mathcal F_{k-1}\\ \rho\left(\s_1,\s_2\right)\leq 3\,r_k}}{\sup}\var\biggl(\sum_{t=1}^{\mu_{T}}\sum_{i\in H_t} w_{i,T}\left( \breve{f}\left(\s_1,\widetilde{\X}_i^{(M)}\left(\frac{i}{T}\right)\right)-\breve{f}\left(\s_2,\widetilde{\X}_i^{(M)}\left(\frac{i}{T}\right)\right)\right)\biggr).
	\end{align*}
	Having in mind that $r_{k_T}\leq r_k$ is fulfilled by construction, we take up on (\ref{eq: FCLT7}) to get  
	\begin{align}\label{eq: termII}
		\text{II}
		&\leq 2\sum_{k=1}^{k_T} D(k) D(k-1)\exp\left(-\frac{1}{2}\cdot\frac{\lambda_k^2}{V_{II,k}+\frac{\breve{M}\lambda_k}{3}}\right)+o(1)\nonumber\\
		%
		%
		%
		%
		%
		%
		&\leq 2\sum_{k=1}^{k_T}\exp\left(2\log\left(D(k)\right)-\frac{1}{2}\cdot\frac{\lambda_k^2}{C_1\,r_k^{\frac{1}{1+\delta}}+C_2\,r_{k_T}^\frac{1}{1+\delta}\lambda}\right)+o(1)\nonumber\\
		%
		%
		%
		%
		%
		%
		&\leq 2\sum_{k=1}^{k_T}\exp\left(-\frac{\bar{C}}{2}\,\lambda_k^2\,r_k^{-\frac{1}{1+\delta}}\right)+o(1)\nonumber\\
		%
		%
		%
		%
		%
		%
		&\leq  2\sum_{k\in\N}\exp\left(-C\,2^{\frac{k}{2(1+\delta)}}r^{-\frac{1}{2(1+\delta)}}\right)+o(1)\nonumber\\
		&\rull 0
	\end{align}
	with the help of (\ref{eq: logDk}) for suitably chosen constants $C_1, \;C_2,\;\bar C<\infty$.
	
	Concerning term III of (\ref{eq: EZD}), we can follow the same steps with $	V_{III}:= C r^{\frac{1}{1+\delta}}$ for some expediently chosen $C<\infty$ and obtain
	\begin{align}\label{eq: FCLT9}
			\text{III}
			&\leq \sum_{\substack{\s_1\in\mathcal F_k,\s_2\in\mathcal F_0\\ \rho\left(\s_1,\s_2\right)\leq 3\,r}}P\biggl(\biggl|\sum_{t=1}^{\mu_{T}}\sum_{i\in H_t} w_{i,T}\left( \breve{f}\left(\s_1,\widetilde{\X}_i^{(M)}\left(\frac{i}{T}\right)\right)-\breve{f}\left(\s_2,\widetilde{\X}_i^{(M)}\left(\frac{i}{T}\right)\right)\right)\biggr|>\frac{\lambda}{27}\biggr)+o(1)\nonumber\\
			&\leq 2\, D^2(0)\exp\left(-\frac{1}{2}\cdot\frac{C_1\lambda^2}{V_{III}+\frac{\breve{M}C_2\lambda}{3}}\right)+o(1)\nonumber\\
			&\leq 2\exp\left(2\log\left(D(0)\right)-\frac{1}{2}\cdot\frac{C_1\lambda^2}{V_{III}+\frac{\breve{M}C_2\lambda}{3}}\right)+o(1)\nonumber\\
			%
			%
			%
			%
			%
			%
			&\leq2\exp\left(2\log\left(D(0)\right)-C\,r^{-\frac{1}{1+\delta}}\right)+o(1)\nonumber\\
			&\rull 0
	\end{align}
	for suited constants $C_1, \;C_2,\;\bar C<\infty$.
	
	Now we move on with the remaining first summand in (\ref{eq: EZD}) and aim at verifying
	\begin{equation}\label{eq.ew-konv}
		\lim_{r\to 0}\limsup_{T\to\infty}\;E\Biggl(\underset{\substack{\s_1,\s_2\in\mathcal S\\ \rho\left(\s_1,\s_2\right)\leq r_{k_T}}}{\sup}\left|\nu_T\left(\s_1,\s_2\right)\right|\Biggr)=0.
	\end{equation}
	Once again, we need some further notation. For $\s\in\mathcal S$ let
	\begin{equation*}
		L_{t,T}\left(\s\right) := \sum_{i\in H_t} w_{i,T}\,f\left(\s,\widetilde{\X}_t^{(M)}\left(\frac{i}{T}\right)\right)
		\quad
		\text{and}
		\quad
		L_{t,T}^0\left(\s\right) :=\zeta_t\, L_{t,T}\left(\s\right),
	\end{equation*}
	where $\left(\zeta_t\right)_{t=1}^{\mu_T}$ are i.i.d.~Rademacher variables independent of $\left(\underline{\varepsilon}_t\right)_{t\in\Z}$. As $\left(L_{t,T}(\s)\right)_{t=1}^{\mu_T}$ consists of independent random variables by construction, we can apply a standard symmetrization lemma (see e.g.~Lemma 2.3.1 in \citet{VW00}) to get
	\begin{equation}\label{eq: FCLT10}
		E\Biggl(\underset{\substack{\s_1,\s_2\in\mathcal S\\ \rho\left(\s_1,\s_2\right)\leq r_{k_T}}}{\sup}\left|\nu_T\left(\s_1,\s_2\right)\right|\Biggr)
		\leq 2\,E\Biggl(\underset{\substack{\s_1,\s_2\in\mathcal S\\ \rho\left(\s_1,\s_2\right)\leq r_{k_T}}}{\sup}\biggl|\sum_{t=1}^{\mu_{T}}\left(L_{t,T}^0\left(\s_1\right)-L_{t,T}^0\left(\s_2\right)\right)\biggr|\Biggr).
	\end{equation}
	Note that $\sum_{t=1}^{\mu_{T}}L_{t,T}^0$ has sub-Gaussian increments conditionally on $L_{1,T},\dots,L_{\mu_T,T}$. This is the case since for $\s_1,\s_2\in\mathcal S$ and $\eta>0$, we get by applying Hoeffding's inequality
	\begin{align}\label{eq: Hoeff}
		P\biggl(\biggl|\sum_{t=1}^{\mu_{T}}L_{t,T}^0\left(\s_1\right)-L_{t,T}^0\left(\s_2\right)\biggr|>\hat{\rho}_{T,2}\left(\s_1,\s_2\right)\eta\,\bigg|\,L_{1,T},\dots,L_{\mu_T,T}\biggr)
		&\leq 2\exp\left(-\frac{\hat{\rho}_{ T,2}\left(\s_1,\s_2\right)^2\eta^2}{2\sum_{t=1}^{\mu_{ T}}\left(L_{t, T}\left(\s_1\right)-L_{t, T}\left(\s_2\right)\right)^2}\right)\nonumber\\
		& = 2\exp\left(-\frac{\eta^2}{2}\right)
	\end{align}
	with the random semimetric
	\begin{equation}\label{eq: DefSM}
		\hat{\rho}_{T,2}\left(\s_1,\s_2\right):=\biggl(\sum_{t=1}^{\mu_{T}}\left(L_{t,T}\left(\s_1\right)-L_{t,T}\left(\s_2\right)\right)^2\biggr)^{1/2}
	\end{equation}
	on $\mathcal S$. We aim at verifying (\ref{eq.ew-konv}) with the help of a maximal inequality for sub-Gaussian processes, which will be more convenient with a different semimetric. To obtain this new semimetric, we note that
	\begin{align}\label{normen}
		\left(L_{t,T}\left(\s_1\right)-L_{t,T}\left(\s_2\right)\right)^2
		& \leq \left(\left|L_{t,T}\left(\s_1\right)\right|+\left|L_{t,T}\left(\s_2\right)\right|\right)^{\frac{2+\delta}{3}} \, \left|L_{t,T}\left(\s_1\right)-L_{t,T}\left(\s_2\right)\right|^{\frac{4-\delta}{3}}\nonumber\\
		& \leq 2^\frac{2+\delta}{3}\; \left|L_{t,T}\right|_{\infty}^\frac{2+\delta}{3}\, \left|L_{t,T}\right|_{\mathrm{Lip}}^\frac{4-\delta}{3}\,\rho\left(\s_1,\s_2\right)^\frac{4-\delta}{3}
	\end{align}
	holds on $(\mathcal S,\rho)$ with $|L_{t,T}|_{\mathrm{Lip}}$ denoting the Lipschitz constant of $L_{t,T}$. 
	By defining
	\begin{equation*}
		Q_T:= 2^\frac{2+\delta}{6}\biggl(\sum_{t=1}^{\mu_{T}} \left|L_{t,T}\right|_{\infty}^\frac{2+\delta}{3}\, \left|L_{t,T}\right|_{\mathrm{Lip}}^\frac{4-\delta}{3}\biggr)^{1/2},
	\end{equation*}
	we get
	\begin{equation}\label{eq: FCLT11}
		EQ_T
		\leq C\biggl(\sum_{t=1}^{\mu_{T}} \left(E\left|L_{t,T}\right|_{\infty}^{2+\delta}\right)^{1/3}\, \left(E\left|L_{t,T}\right|_{\mathrm{Lip}}^{\frac{4-\delta}{2}}\right)^{2/3}\biggr)^{1/2}.
	\end{equation}
	Regarding the expectations on the RHS, we obtain  
	\begin{align*}
		E\left|L_{t,T}\right|_{\infty}^{2+\delta}
		&\leq\biggl\|\sum_{i\in H_t}w_{i,T}\left(C\,g\left(\widetilde{\X}^{(M)}_i\left(\frac{i}{T}\right)\right)+\left|f\left(\un 0,\widetilde{\X}^{(M)}_i\left(\frac{i}{T}\right)\right)\right|\right)\biggr\|_{2+\delta}^{2+\delta}\\
		%
		%
		%
		%
		%
		%
		&	\leq C \,\Biggl(\sum_{i\in H_t}d_T^{1/2}\left(\left\|  \,g\left(\widetilde{\X}^{(M)}_i\left(\frac{i}{T}\right)\right)\right\|_{2+\delta} +\left\|f\left(\un 0,\widetilde{\X}^{(M)}_i\left(\frac{i}{T}\right)\right)\right\|_{2+\delta}\right)\Biggr)^{2+\delta}
		\\
		%
		%
		%
		%
		%
		%
		&\leq C\,\kappa_T^{2+\delta}\,d_T^{-\frac{2+\delta}{2}}
	\end{align*}
	and
	\begin{equation*}
		E\left|L_{t,T}\right|_{\mathrm{Lip}}^{\frac{4-\delta}{2}}
		%
		%
		%
		%
		%
		%
		\leq\biggl(\sum_{i\in H_t}w_{i,T}\left\|g\left( \widetilde{\X}^{(M)}_i\left(\frac{i}{T}\right)\right)\right\|_{\frac{4-\delta}{2}}\biggr)^{\frac{4-\delta}{2}}
		\leq C\,\kappa_T^{\frac{4-\delta}{2}}\,d_T^\frac{4-\delta}{4}.
	\end{equation*}
	Together, we have for (\ref{eq: FCLT11})
	\begin{equation}\label{eq: EQT}
		EQ_T
		\leq C \biggl(\sum_{t=1}^{\mu_{T}} \kappa_T^{(2+\delta)\frac{1}{3}+\frac{4-\delta}{2}\cdot\frac{2}{3}}\, d_T^{\frac{1}{2}\left((2+\delta)\frac{1}{3}+\frac{4-\delta}{2}\cdot\frac{2}{3}\right)}\biggr)^{1/2}
		=C\left(\mu_T\,\kappa_T^2\,d_T^{-1}\right)^{1/2}\nonumber\\
		%
		%
		%
		\leq C\,d_T^{\frac{1}{4m}}.
	\end{equation}
	Having the definition of $Q_T$ in mind, it holds
	\begin{align*}
		\hat{\rho}_{T,2}\left(\s_1,\s_2\right)
		&\leq Q_T \;\rho\left(\s_1,\s_2\right)^{\frac{4-\delta}{6}}=:\breve{\rho}_{T}\left(\s_1,\s_2\right),
	\end{align*}
	and $\breve{\rho}_{T}$ is again a random semimetric as $\frac{4-\delta}{6}\in(0,1)$.
	Now we make use of Corollary~2.2.8 of \citet{VW00} to get
	\begin{align*}
		&E\Biggl(\underset{\substack{\s_1,\s_2\in\mathcal S\\ \rho\left(\s_1,\s_2\right)\leq r_{k_T}}}{\sup}\biggl|\sum_{t=1}^{\mu_{T}}\left(L_{t,T}^0\left(\s_1\right)-L_{t,T}^0\left(\s_2\right)\right)\biggr|\,\Bigg| \,L_{1,T},\dots,L_{\mu_T,T}\Biggr)\\
		%
		%
		%
		%
		%
		&\leq C_3\int_{0}^{Q_T r_{k_T}^{\frac{4-\delta}{6}}}\left(\log\left(D\left(u,\mathcal S,                {\rho}_T\right)\right)\right)^{1/2}du\nonumber\\
		&	\leq C_3\int_{0}^{Q_T r_{k_T}^{\frac{4-\delta}{6}}}\Biggl(\log\Biggl( D\biggl(\left(\frac{u}{Q_T}\right)^{\frac{6}{4-\delta}},\mathcal S, \rho\biggr)\Biggr)\Biggr)^{1/2}du\nonumber\\
		&\leq C\int_{0}^{Q_T r_{k_T}^{\frac{4-\delta}{6}}}\Biggl(\log\biggl( \left(\frac{u}{Q_T}\right)^{-\frac{6}{4-\delta}}+1\biggr)^d\Biggr)^{1/2}du\\
		%
		%
		%
		%
		%
		%
		&\leq C\,Q_T\int_{0}^{r_{k_T}^{\frac{4-\delta}{6}}}u^{-\frac{3}{4-\delta}}\,du\\
		&\leq  C\,Q_T\,r_{k_T}^{\frac{1-\delta}{6}}.
	\end{align*}
	Returning to (\ref{eq: FCLT10}), we get from (\ref{eq: defkn}) and (\ref{eq: EQT}) 
	\begin{equation*}
		E\Biggl(\underset{\substack{\s_1,\s_2\in\mathcal S\\ \rho\left(\s_1,\s_2\right)\leq r_{k_T}}}{\sup}\biggl|\sum_{t=1}^{\mu_{T}}\left(L_{t,T}^0\left(\s_1\right)-L_{t,T}^0\left(\s_2\right)\right)\biggr|\Biggr) 
		\leq   C\,EQ_T\, r_{k_T}^{\frac{1-\delta}{6}}
		%
		%
		%
		%
		%
		%
		\leq C\,d_T^{-\frac{1}{12m}},
	\end{equation*}
	which tends to 0 as $T\to\infty$, and the proof is completed.
\end{proof}\medskip
 
As already indicated in the Bootstrap Algorithm \ref{al: BS}, some indices can be responsible for a change in the sign of $k_i$. The following definition groups those indices for easier handling:
	\begin{defi}[Endpoints]\ \\
		Considering Algorithm \ref{al: BS}, all indices $t\in{1,\dots,T}$ which might cause a sign switch are called \emph{endpoints} $EP=EP_1\cup EP_2$ with
		\begin{equation*}
			EP_1:=\left\{t\in\{1,\dots,T\}\,\middle|\,1\leq t\leq TD_T\right\}
	\quad\text{and}\quad
			EP_2:=\left\{t\in\{1,\dots,T\}\,\middle|\,T-TD_T< t\leq T\right\}.
		\end{equation*}
\end{defi}
\begin{proof}[Proof of Theorem~\ref{th: CLTBS}]	
	  This proof is inspired by \citet{DPP13} but enhanced to meet the demands imposed by our assumptions.
	   For sake of notational simplicity, we consider $J=1$ only. We split the sum up in one sum containing all indices belonging to whole bootstrap blocks without including endpoints and two with the remaining indices in the way that we have
	\begin{align}\label{eq: T4BS0}
		\sum_{t=1}^{T} w_{t,T}\, f\left(\s, \X_{t,T}^\star\right)
			&\begin{multlined}[t][0.81\linewidth]
				=\sum_{t=1}^{L_T\left\lceil\left(TD_T+1\right)/L_T\right\rceil} w_{t,T}\, f\left(\s, \X_{t,T}^\star\right)+\sum_{t=L_T\left\lceil\left(TD_T+1\right)/L_T\right\rceil+1}^{L_T\left\lfloor\left(T-TD_T\right)/L_T\right\rfloor} w_{t,T}\, f\left(\s, \X_{t,T}^\star\right)\\
				+\sum_{t=L_T\left\lfloor\left(T-TD_T\right)/L_T\right\rfloor+1}^{T} w_{t,T}\, f\left(\s, \X_{t,T}^\star\right)
			\end{multlined}
	\nonumber\\
		&=:\text{I}+\text{II}+\text{III}.
	\end{align}
Recall that $L_T=o(TD_T)$. In view of Lemma~\ref{le: KovEP} and Assumption~\ref{as: Ra}, we can bound 
	\begin{equation*}
		E^\star|\text{I}|\leq 	E^\star\biggl|\sum_{t=1}^{L_T\left\lceil\left(TD_T+1\right)/L_T\right\rceil} w_{t,T}\, f\left(\s, \X_{t,T}^\star\right)\biggr|
		=\go_P\left(TD_T\,d_T^{-1/2}\right)\\
		=\go_P\left(d_T^{-\frac{\delta}{2\left(2+\delta\right)}}\right).
	\end{equation*}
	 Due to the similar structure of III, we obtain $\text{I}+\text{III}=o_{P^\star}\left(1\right)$. The remaining term II of (\ref{eq: T4BS0})  can be rewritten to make the single blocks visible, to wit
	\begin{equation*}
		\text{II}
		=\sum_{t=\left\lceil\left(TD_T+1\right)/L_T\right\rceil}^{\left\lfloor\left(T-TD_T\right)/L_T\right\rfloor-1}\sum_{j=1}^{L_T} w_{tL_T+j,T}\,\f^\star\left(\s,\X_{tL_T+j,T}^\star\right)
		=:\sum_{t=\left\lceil\left(TD_T+1\right)/L_T\right\rceil}^{\left\lfloor\left(T-TD_T\right)/L_T\right\rfloor-1}\xi_{t,T}^\star.
	\end{equation*}
	Since the bootstrap blocks are independent, the newly defined $\left(\xi^\star_{t,T}\right)$ are independent as well. Lemma~\ref{th: VarBS} gives us
	\begin{equation}\label{eq: T4BS1}
		P-\lim_{T\to\infty}\sum_{t=\left\lceil\left(TD_T+1\right)/L_T\right\rceil}^{\left\lfloor\left(T-TD_T\right)/L_T\right\rfloor-1}\var^\star\left(\xi_{t,T}^\star\right)=\mathbf V\left(\s,\s\right).
	\end{equation}
	Thus, it suffices to consider the case $\mathbf V\left(\s,\s\right)>0$ only. 
	 We aim for applying the classical central limit theorem for independent data using Lyapunov's condition. To this end, note that
		\begin{equation}\label{eq: T4BS3}
		E\left(E^\star\left|\xi_{t,T}^\star\right|^{2+\delta}\right)
		\leq C\,d_T^{-\frac{2+\delta}{2}}\biggl(\sum_{j=1}^{L_T}\ind{\left\{w_{tL_T+j,T}>0\right\}}\biggr)^{2+\delta}.
	\end{equation}
	Furthermore, it holds
	\begin{equation}\label{eq: sumA}
		\sum_{t=\left\lceil\left(TD_T+1\right)/L_T\right\rceil}^{\left\lfloor\left(T-TD_T\right)/L_T\right\rfloor-1}\biggl(\sum_{j=1}^{L_T}\ind{\left\{w_{tL_T+j,T}>0\right\}}\biggr)^{2+\delta}
		\leq d_T\,L_T^{1+\delta}.
	\end{equation}
	Now, combining (\ref{eq: T4BS1}) and (\ref{eq: T4BS3}) and using (\ref{eq: sumA}), we obtain
	\begin{equation*}
		\frac{\sum_{t=\left\lceil\left(TD_T+1\right)/L_T\right\rceil}^{\left\lfloor\left(T-TD_T\right)/L_T\right\rfloor-1}E^\star\left|\xi_{t,T}^\star\right|^{2+\delta}}{\left(\sum_{t=\left\lceil\left(TD_T+1\right)/L_T\right\rceil}^{\left\lfloor\left(T-TD_T\right)/L_T\right\rfloor-1}\var^\star\left(\xi_{t,T}^\star\right)\right)^{\frac{2+\delta}{2}}}
		=\frac{d_T\,L_T^{1+\delta}\go_P\left(d_T^{-\frac{2+\delta}{2}}\right)}{\left(V\left(\s,\s\right)\right)^{\frac{2+\delta}{2}}+o_P(1)}
		=o_P(1).
	\end{equation*}
 Thereby, Lyapunov's condition is fulfilled, which implies asymptotic normality as desired. The second assertion of the theorem is an immediate consequence of the first one due to continuity of the Gaussian distribution function, see Lemma 2.11 of \citet{VV98}.
\end{proof}\medskip

Within the next proof and the corresponding auxiliary results in the Appendix, we will use some additional notation. In particular, we abbreviate $\mathbb X_T:=(X_{1,T}',\dots, X_{T,T}')'$ and $P(\cdot\mid \mathbb X_T=\un x_T)$ by  $P^\star_T(\cdot)$. Corresponding expectations and variances are denoted by $E^\star_T$ and $\var^\star_T$.

\begin{proof}[Proof of Theorem~\ref{th: FCLTBS}]
	With the same arguments as in the proof of Theorem~\ref{th: FCLT}, it is sufficient  to show
that there exist sets $\left(\Omega_T\right)_{T\in\N}$ with $P\left(\mathbb X_T\in\Omega_T\right)\to 1$ as $T\to\infty$ such that for any $(\un x_T)_T$ with $x_T\in\Omega_T$ for all $T$ it holds
	\begin{equation}\label{eq: BS-T}
		\underset{r\to 0}{\lim}\,\limsup_{T\to\infty}\, P\biggl(\sup_{\rho\left(\s_1,\s_2\right)<r}\biggl|\sum_{t=1}^T w_{t,T}\left(\f^\star\left(\s_1, \X^\star_{t,T}\right)-\f^\star\left(\s_2, \X^\star_{t,T}\right)\right)\biggr|>\lambda\,\bigg|\,\mathbb X_T=\un x_T\biggr)=0
	\end{equation}
	with $\lambda>0$, whereas $\f^\star$ is defined in (\ref{eq: CentBS}).

Here, we only consider the case (a) in Assumption~\ref{as: Fu2} and~\ref{as: G-BS}. Part (b) is deferred to Lemma~\ref{le: BS-FCLT-b} in the Appendix. First, we define 
\begin{equation*}
	G_T=\biggl\{\un x_T\,\bigg|\, \sum_{t=L_T \left\lfloor T/L_T\right\rfloor+1}^{T}  \,E^\star_T[ g\left( \X^\star_{t,T}\right)]\leq   L_T \, d_T^{\delta/(2(1+\delta))}\biggr\}.
\end{equation*}
From $\sup_{1\leq t\leq T}E|E^\star g\left( \X^\star_{t,T}\right)|\leq C$, we get $P(G_T)\to 1$ as $T\to\infty$ applying Markov's inequality. 
	With the sets $A_T$, $B_T$ and $K_T$ established in Lemmata \ref{le: UB}, \ref{le: BT} and \ref{le: K}, respectively,  
	we set
	\begin{equation*}
		\Omega_T:=A_T\cap B_T\cap G_T\cap K_T.
	\end{equation*} 
	In view of the above mentioned results, it holds $\lim_{T\to\infty}P\left(\mathbb X_T\in\Omega_T\right)=1$.
	
	After these preparations, we start by splitting the LHS of (\ref{eq: BS-T}) such  that we get one sum containing the indices of whole independent blocks and a second one containing the remaining indices:
	\begin{align}\label{eq: BS-T1}
		&P^\star_T\biggl(\sup_{{\rho\left(\s_1,\s_2\right)<r}}\biggl|\sum_{t=1}^T w_{t,T}\left(\f^\star\left(\s_1, \X^\star_{t,T}\right)-\f^\star\left(\s_2, \X^\star_{t,T}\right)\right)\biggr|>\lambda \biggr)\nonumber\\
		&\begin{multlined}[b][0.88\linewidth]
			\leq P^\star_T\biggl(\sup_{{\rho\left(\s_1,\s_2\right)<r}}\biggl|\sum_{t=1}^{L_T \left\lfloor T/L_T\right\rfloor} w_{t,T}\left(\f^\star\left(\s_1, \X^\star_{t,T}\right)-\f^\star\left(\s_2, \X^\star_{t,T}\right)\right)\biggr|>\frac{\lambda}{2} \biggr)\\
			+P^\star_T\biggl(\sup_{{\rho\left(\s_1,\s_2\right)<r}}\biggl|\sum_{t=L_T \left\lfloor T/L_T\right\rfloor+1}^{T} w_{t,T}\left(\f^\star\left(\s_1, \X^\star_{t,T}\right)-\f^\star\left(\s_2, \X^\star_{t,T}\right)\right)\biggr|>\frac{\lambda}{2} \biggr).
		\end{multlined}
	\end{align}
	By construction of $G_T$, the second sum on the RHS of (\ref{eq: BS-T1}) can be bounded by 
	\begin{align*}
		\frac{4}{\lambda} \sum_{t=L_T \left\lfloor T/L_T\right\rfloor+1}^{T} w_{t,T}\,E^\star_T\left[ \sup_{{\rho\left(\s_1,\s_2\right)<r}}\left|f\left(\s_1, \X^\star_{t,T}\right)-f\left(\s_2, \X^\star_{t,T}\right)\right|  \right]
		&\leq{C\,r\, d_T^{-1/2}}\,\sum_{t=L_T \left\lfloor T/L_T\right\rfloor+1}^{T}   E^\star_T g\left( \X^\star_{t,T}\right)\\
		&=\go\left(d_T^{-\frac{1-\delta}{2(1+\delta)}}\right),
	\end{align*}
	which tends to 0 as $T\to\infty$ leaving us with the first sum on the RHS of (\ref{eq: BS-T1}) to deal with. Similarly to the proof of Theorem~\ref{th: FCLT}, we define $\nu_{T}^\star(\s)$ as well as $\nu_{T}^\star\left(\s_1,\s_2\right)$ by 
	\begin{equation*}
		\nu_{T}^\star(\s):= \sum_{t=0}^{\left\lfloor T/L_T\right\rfloor-1}\sum_{j=1}^{L_T} w_{tL_T+j,T}\,\f^\star\left(\s_1,\X_{tL_T+j,T}^\star\right)
		\quad\text{and}\quad
		\nu_{T}^\star\left(\s_1,\s_2\right):= \nu_{T}^\star\left(\s_1\right) - \nu_{T}^\star\left(\s_2\right),
	\end{equation*}
	respectively.  We will use the same notation regarding the sequence $(r_k)_k$, the index sets $\left(\mathcal F_k\right)_{k=0}^{k_T}$ and the maps $\left(\pi_k\right)_{k=0}^{k_T}$ as in the proof of Theorem~\ref{th: FCLT}. As before, we split the left over sum of (\ref{eq: BS-T1}) up into
	\begin{align}\label{eq: BS-T1a}
		&P^\star_T\Biggl(\sup_{\overset{\s_1,\s_2\in\mathcal S}{\rho\left(\s_1,\s_2\right)<r}}\biggl|\sum_{t=1}^{L_T \left\lfloor T/L_T\right\rfloor} w_{t,T}\left(\f^\star\left(\s_1, \X^\star_{t,T}\right)-\f^\star\left(\s_2, \X^\star_{t,T}\right)\right)\biggr|>\frac{\lambda}{2} \Biggr)\nonumber\\
		&\begin{multlined}[t][\linewidth]
				\leq P^\star_T\Biggl(2\sup_{\overset{\s_1,\s_2\in\mathcal S}{\rho\left(\s_1,\s_2\right)\leq r_{k_T}}}\left|\nu_T^\star\left(\s_1,\s_2\right)\right|>\frac{\lambda}{6} \Biggr)
			+P^\star_T\Biggl(2\sum_{k=1}^{k_T}\sup_{\overset{\s_1 \in \mathcal F_k, \s_2\in\mathcal F_{k-1}}{\rho\left(\s_1,\s_2\right)\leq 3\,r_k}}\left|\nu_T^\star\left(\s_1,\s_2\right)\right|>\frac{\lambda}{6} \Biggr)\\
			+ P^\star_T\Biggl(\sup_{\overset{\s_1,\s_2\in\mathcal F_0}{\rho\left(\s_1,\s_2\right)\leq 3\,r}}\left|\nu_T^\star\left(\s_1,\s_2\right)\right|>\frac{\lambda}{6} \Biggr)
		\end{multlined}
		\nonumber\\
		&=:\text{I}+\text{II}+\text{III}
	\end{align}
	postponing the magnitude of $k_T$. One can use Bernstein's inequality for the discussion of terms II and~III. Here, we will carry out the details for~II, only. To this end, we consider the bootstrap variance of $\nu_T^\star\left(\s_1,\s_2\right)$ first. With the help of Lemmata~\ref{le: BT} and \ref{le: DiffVar}, we obtain (with $(c_T)_T$ as in Lemma~\ref{le: BT})
	\begin{equation*}
		\var^\star_T\left(\nu_T^\star\left(\s_1,\s_2\right)\right)
		%
		%
		%
		\leq C\left\{  \left|\s_1-\s_2\right|_1^{1/2}\,+\, L_T^{-1}  \,+\,c_T\,L_T\left(T^{-\frac{3+\delta}{2(4+\delta)}}\,+\, D_T^{\frac{3+\delta}{2(4+\delta)}}\,+\, (TD_T)^{-\frac{1}{2}}\right) \right\}.
\end{equation*}
Next, Lemma~\ref{le: UB} allows us to bound
\begin{equation}\label{eq: UB}
	\biggl|\sum_{j=1}^{L_T} w_{tL_T+j,T}\left(\f^\star\left(\s_1,\X_{tL_T+j,T}^\star\right)-\f^\star\left(\s_2,\X_{tL_T+j,T}^\star\right)\right)\biggr|
	\leq C\,L_T\,d_T^{-\frac{2+\delta^2}{2(4+\delta)}}
	%
	%
	%
	%
	%
	%
	<C\,d_T^{-\frac{2-2\delta+\delta^3}{2(1+\delta)(4+\delta)}}
\end{equation}
due to Assumption~\ref{as: Ra}.	
Abbreviating  $y_T:=L_T^{-1}  \,+\,c_T\,L_T\left(T^{-\frac{3+\delta}{2(4+\delta)}}\,+\, D_T^{\frac{3+\delta}{2(4+\delta)}}\,+\, (TD_T)^{-1/2}\right)$, we specify $k_T$ such that
\begin{equation}\label{eq: BS-T2}
	\left(\max\left\{y_T,d_T^{-\frac{2-2\delta+\delta^3}{2(1+\delta)(4+\delta)}}\right\}\right)^{\frac{(1+\delta)(2+\delta)(4+\delta)}{\delta^2\left(\delta-\delta^2\right)}}\leq r_{k_T}\leq L_T^{-\frac{1+\delta}{\delta^2}}.
\end{equation}
Tedious straightforward calculations similar to  the proof of Theorem~3.13 in \citet{B21} show that the upper bound is strictly larger than the lower bound for all sufficiently large $T$.
The specific choice of the upper bound in (\ref{eq: BS-T2}) will become relevant during the examination of term I of equation (\ref{eq: BS-T1a}), while the lower bound is essential for a successful application of Bernstein's inequality. In particular for~II in \eqref{eq: BS-T2}, we will make us of $\rho\left(\s_1,\s_2\right)\leq 3\,r_k$ to obtain 
\begin{equation}\label{eq: BS-T4}
	\var^\star\left(\nu_T^\star\left(\s_1,\s_2\right)\right)\leq C\Biggl(r_k^{1/2}+r_{k_T}^{\frac{\delta^2\left(\delta-\delta^2\right)}{(1+\delta)(2+\delta)(4+\delta)}}\Biggr)
	\leq C\,r_k^{\frac{\delta^2\left(\delta-\delta^2\right)}{(1+\delta)(2+\delta)(4+\delta)}}
\end{equation}
for $T>\bar{T}$ as it holds $r_{k_T}<r_k$ for all $k=0,1,\dots,k_T$ and because of
\begin{equation*}
	\frac{1}{2}-\frac{\delta^2\left(\delta-\delta^2\right)}{(1+\delta)(2+\delta)(4+\delta)}>0\quad\forall \delta \in (0,1).
\end{equation*}
After having determined the bounds to be used in Bernstein's inequality, we turn our attention to $\lambda$ and, as in the proof of Theorem~\ref{th: FCLT}, define a sequence $\left(\lambda_k\right)_{k\in\N}$ satisfying
\begin{equation*}
	\sum_{k\in\N}\lambda_k\leq \frac{\lambda}{3} 
\end{equation*}
for $r$ sufficiently small by setting
\begin{equation*}
	\lambda_k := r_k^\frac{\delta^2\left(\delta-\delta^2\right)}{4(1+\delta)(2+\delta)(4+\delta)}\lor \left(\frac{4}{\bar{C}}\log\left(D(k)\right) r_k^\frac{\delta^2\left(\delta-\delta^2\right)}{(1+\delta)(2+\delta)(4+\delta)}\right)^{1/2},\quad k\in \N,
\end{equation*}
for a finite constant $\bar{C}>0$, which will be specified further during the upcoming calculations. Again, $D(k)=\mathcal O\left(r_k^{-d}\right)$ ensures the summability of $\left(\lambda_k\right)_{k}.$

At this point, we return to II in (\ref{eq: BS-T1a}) and use both the definition of $\left(\lambda_k\right)_{k=1}^{k_T}$ and Bernstein's inequality with the previously established upper bounds in (\ref{eq: UB}) and (\ref{eq: BS-T4}) to get
\begin{align*}
	&\limsup_{T\to\infty}P^\star_T\Biggl(2\sum_{k=1}^{k_T}\sup_{\overset{\s_1 \in \mathcal F_k, \s_2\in\mathcal F_{k-1}}{\rho\left(\s_1,\s_2\right)\leq 3\,r_k}}\left|\nu_T^\star\left(\s_1,\s_2\right)\right|>\frac{\lambda}{6} \Biggr)\\
	%
	%
	%
	%
	%
	%
	%
	&\leq\limsup_{T\to\infty} 2 \sum_{k=1}^{k_T} D(k) D(k-1) \exp\left(-\frac{1}{2}\cdot\frac{\lambda_k^2}{C_1\,r_k^\frac{\delta^2\left(\delta-\delta^2\right)}{(1+\delta)(2+\delta)(4+\delta)}+C_2\,d_T^{\frac{2-2\delta+\delta^3}{2(1+\delta)(4+\delta)}}\lambda_k/3}\right)\\
	%
	%
	%
	%
	&\leq\limsup_{T\to\infty} 2 \sum_{k=1}^{k_T}\exp\left(2 \log(D(k))-\bar{C}\frac{\lambda_k^2}{r_k^\frac{\delta^2\left(\delta-\delta^2\right)}{(1+\delta)(2+\delta)(4+\delta)}}\right)\\
	%
	%
	%
	%
	%
	%
	%
	%
	%
	&\leq 2 \sum_{k\in\N}\exp\left(-C\, 2^{\frac{k\delta^2\left(\delta-\delta^2\right)}{2(1+\delta)(2+\delta)(4+\delta)}}r^{-\frac{\delta^2\left(\delta-\delta^2\right)}{2(1+\delta)(2+\delta)(4+\delta)}}\right)\\
	&\underset{r\to 0}{\longrightarrow} 0.
\end{align*}
%
Now term I of (\ref{eq: BS-T1a}) is left, and, as in the proof of Lemma~\ref{th: FCLT}, we aim for making use of a symmetrization lemma.  To this end, let
\begin{equation}\label{eq: Lstern}
	L_{t,T}^\star\left(\s\right):= \sum_{j=1}^{L_T} w_{tL_T+j,T}\,f\left(\s,\X_{tL_T+j,T}^\star\right)\quad\text{and}\quad L_{t,T}^{\star,0}\left(\s\right):= \zeta_t\,L_{t,T}^\star\left(\s\right),\quad \s\in\mathcal S,
\end{equation}
%
%
%
where $\left(\zeta_t\right)_{t=0}^{\left\lfloor T/L_T\right\rfloor-1}$ are again i.i.d.~Rademacher variables, but this time, they are independent of the series of uniformly distributed random variables used in Algorithm \ref{al: BS} and $\mathbb X_T$. This gives
\begin{align}\label{eq: BS-T4a}
		\text{I}&=	P^\star_T\biggl(\sup_{{\s_1,\s_2\in\mathcal S,\;}{\rho\left(\s_1,\s_2\right)\leq r_{k_T}}}\left|\nu_T^\star\left(\s_1,\s_2\right)\right|>\frac{\lambda}{12} \biggr)\nonumber\\
		&\leq \frac{12}{\lambda}\,E^\star_T\biggl(\sup_{{\s_1,\s_2\in\mathcal S,\;\rho\left(\s_1,\s_2\right)\leq r_{k_T}}}\left|\nu_T^\star\left(\s_1,\s_2\right)\right|\biggr)\nonumber\\	
		&\leq \frac{24}{\lambda}\,E^\star_T\biggl(\sup_{ {\s_1,\s_2\in\mathcal S,\;\rho\left(\s_1,\s_2\right)\leq r_{k_T}}}\biggl|\sum_{t=0}^{\left\lfloor T/L_T\right\rfloor-1}\left(L_{t,T}^{\star,0}\left(\s_1\right)-L_{t,T}^{\star,0}\left(\s_2\right)\right)\biggr|\biggr).  
\end{align}
Next, we define
\begin{equation}\label{eq: BS-T5}
	\hat{\rho}_{T,2}\left(\s_1,\s_2\right):=\biggl(\sum_{t=0}^{\left\lfloor T/L_T\right\rfloor-1}\left(L_{t,T}^\star\left(\s_1\right)-L_{t,T}^\star\left(\s_2\right)\right)^2\biggr)^{1/2},\quad \s_1,\s_2\in\mathcal S.
\end{equation}
For any $\eta>0$, Hoeffding's inequality gives
$$\begin{aligned}
	&P^\star\Biggl(\sup_{\overset{\s_1,\s_2\in\mathcal S}{\rho\left(\s_1,\s_2\right)\leq r_{k_T}}}\biggl|\sum_{t=0}^{\left\lfloor T/L_T\right\rfloor-1}L_{t,T}^{\star,0}\left(\s_1\right)-L_{t,T}^{\star,0}\left(\s_2\right)\biggr|>\hat{\rho}_{T,2}\left(\s_1,\s_2\right)\eta\,\Bigg |\,L_{1,T}^\star,\dots,L_{L_T \left\lfloor T/L_T\right\rfloor,T}^\star\Biggr)\\
	&\leq 2\exp\left(-\frac{\hat{\rho}_{T,2}\left(\s_1,\s_2\right)^2\eta^2}{2\sum_{t=1}^{L_T \left\lfloor T/L_T\right\rfloor}\left(L_{t,T}^\star\left(\s_1\right)-L_{t,T}^\star\left(\s_2\right)\right)^2}\right)\\& =2\exp\left(-\frac{\eta^2}{2}\right).
\end{aligned}$$
Thus, $\sum_{t=0}^{\left\lfloor T/L_T\right\rfloor-1}L_{t,T}^{\star,0}$ possesses sub-Gaussian increments conditionally on $L_{0,T}^\star,\dots,L_{\left\lfloor T/L_T\right\rfloor-1,T}^\star$. We continue by establishing an upper bound for the difference in (\ref{eq: BS-T5}) as follows:
\begin{equation}\label{normenBS}
	\left(L_{t,T}^\star\left(\s_1\right)-L_{t,T}^\star\left(\s_2\right)\right)^2
	\leq 2^{1-\delta}\; \left|L_{t,T}^\star\right|_{\infty}^{1-\delta} \left|L_{t,T}^\star\right|_{\mathrm{Lip}}^{1+\delta}\rho\left(\s_1,\s_2\right)^{1+\delta},\quad \s_1,\s_2\in\mathcal S.
\end{equation}
In order to create a more suitable semimetric during the further course of this proof, we define
\begin{equation*}
	Q_T:= 2^\frac{1-\delta}{2}\biggl(\sum_{t=0}^{\left\lfloor T/L_T\right\rfloor-1} \left|L_{t,T}^\star\right|_{\infty}^{1-\delta}\, \left|L_{t,T}^\star\right|_{\mathrm{Lip}}^{1+\delta}\biggr)^{1/2}
\end{equation*}
and introduce a new semimetric using
\begin{equation*}
	\hat{\rho}_{T,2}\left(\s_1,\s_2\right)\leq Q_T\rho\left(\s_1,\s_2\right)^{\frac{1+\delta}{2}}=:\breve{\rho}_T\left(\s_1,\s_2\right).
\end{equation*}
Now we come back to (\ref{eq: BS-T4a}). Using Corollary~2.2.8 of \citet{VW00}, we obtain
\begin{align*}
	\text{I} &\leq \frac{24}{\lambda}\,E^\star_T\left[E^*\biggl(\sup_{ {\s_1,\s_2\in\mathcal S,\;\rho\left(\s_1,\s_2\right)\leq r_{k_T}}}\biggl|\sum_{t=0}^{\left\lfloor T/L_T\right\rfloor-1}\left(L_{t,T}^{\star,0}\left(\s_1\right)-L_{t,T}^{\star,0}\left(\s_2\right)\right)\biggr|\,\bigg|\,L_{1,T}^\star,\dots,L_{L_T \left\lfloor T/L_T\right\rfloor,T}^\star\biggr)\right]\\
	&= \frac{24}{\lambda}\,E^\star_T\left[E^*\biggl(\sup_{ {\s_1,\s_2\in\mathcal S,\,\breve{\rho}\left(\s_1,\s_2\right)\leq Q_T r_{k_T}^{\frac{1+\delta}{2}}}}\biggl|\sum_{t=0}^{\left\lfloor T/L_T\right\rfloor-1}\left(L_{t,T}^{\star,0}\left(\s_1\right)-L_{t,T}^{\star,0}\left(\s_2\right)\right)\biggr|\,\bigg|\,L_{1,T}^\star,\dots,L_{L_T \left\lfloor T/L_T\right\rfloor,T}^\star\biggr)\right]\\
	&\leq C\,E^\star_T \Biggl(\int_{0}^{Q_T r_{k_T}^{\frac{1+\delta}{2}}}\biggl(\log D\left(u,\mathcal S,\breve{\rho}_T\right)\biggr)^{1/2}du\Biggr)\\
	&\leq C \,E^\star_T\Biggl(\int_{0}^{Q_T r_{k_T}^{\frac{1+\delta}{2}}}\biggl(\log \biggl(\left(\frac{u}{Q_T}\right)^{-\frac{2}{1+\delta}}+1\biggr)^d\biggr)^{1/2}du\Biggr)\nonumber\\
	&=C\, E^\star_T[ Q_T]\, \int_{0}^{r_{k_T}^{\frac{1+\delta}{2}}}\left(\log \left(u^{-\frac{2}{1+\delta}}+1\right)^d\right)^{1/2}du.
\end{align*}
Moreover, Jensen's inequality yields
\begin{equation*}
	E^\star_T[ Q_T]\leq C_1\biggl(\sum_{t=0}^{\left\lfloor T/L_T\right\rfloor-1} E^\star\left(\left|L_{t,T}^\star\right|_{\infty}^{1-\delta}\, \left|L_{t,T}^\star\right|_{\mathrm{Lip}}^{1+\delta}\right)\biggr)^{1/2}\leq  C\,L_T^{\frac{1+\delta}{2}},
\end{equation*}
where the latter inequality follows from Lemma~\ref{le: K}. Plugging this result into the previous calculations and remembering the upper bound of $r_{k_T}$ in (\ref{eq: BS-T2}) gives 
\begin{equation*}
	\text{I}
	\leq\limsup_{T\to\infty}\, C\,L_T^{\frac{1+\delta}{2}}\,\int_{0}^{r_{k_T}^{\frac{1+\delta}{2}}}\left(\log \left(u^{-\frac{2}{1+\delta}}+1\right)^d\right)^{1/2}du \leq  \limsup_{T\to\infty}\, C\,L_T^{\frac{1+\delta}{2}}\,L_T^{-\frac{1+\delta}{2\delta}}=0.
\end{equation*}
This terminates the proof. 
\end{proof}\medskip

\noindent
{\bf Acknowledgement} This research is partly funded by the Volkswagen Foundation, Germany (Professorinnen für Niedersachsen
des Niedersächsischen Vorab).
\bigskip

\noindent
{\bf Disclosure Statements}
The authors report there are no competing interests to declare.
\bigskip
\linespread{1.0}

\bibliographystyle{harvard}

\bigskip

\appendix
\section{Appendix/Online Supplement}

Throughout this section again, $C$ denotes a generic constant that may change its value from line to line.

\subsection{Auxiliary results required for Section~2}

\begin{lemma}\label{le: AbschM}\ 
	Under Assumptions \ref{as: Pro} valid for $k=0$ and Assumption~\ref{as: Fu}, we have for every $M\in\N$ and	some $C_{trunc}<\infty$
	\begin{equation*}
	\sup_{u\in[0,1]}	\left\|f\left(\s,\widetilde{\X}_t(u)\right)-f\left(\s,\widetilde{\X}_t^{(M)}(u)\right)\right\|_1\leq   \,C_{trunc}\,\vartheta^M\,.
	\end{equation*}
\end{lemma}

\begin{proof}

		With the use of \eqref{eq.au} and the Lipschitz condition (\ref{eq: Lip}), it holds for some $C_{trunc}<\infty$
	\begin{align*}
	\sup_{u\in[0,1]}	\left\|f\left(\s,\widetilde{\X}_t(u)\right)-f\left(\s,\widetilde{\X}_t^{(M)}(u)\right)\right\|_1\,
	&	\leq C_{Lip}\, \sup_{u\in[0,1]}  \left\|\widetilde{\X}_t(u)-\widetilde{\X}_t^{(M)}(u)\right\|_1 \\
	&\leq  C_{Lip}\,\sup_{u\in[0,1]}  \left\|\sum_{|j|\geq M}A(u,j)\,\un\varepsilon_{t-j}\right\|_1 \\
	&\leq   C_{trunc}\,\vartheta^M.
	\end{align*}
\end{proof}\medskip

\begin{lemma}\label{le: DiffProd}\ 
	Suppose Assumptions \ref{as: Pro} for $k=0$ and Assumption~\ref{as: Fu} are fulfilled. Then, it holds uniformly for all $t_1,t_2\in \Z$, $\left(\s_1,u_1\right), \left(\s_2,u_2\right)\in\mathcal S\times[0,1]$
	\begin{equation*}
		\left\|\left(f\left(\s_1,\widetilde{\X}_{t_1}\left(u_1\right)\right)-f\left(\s_1,\widetilde{\X}_{t_1}^{(M)}\left(u_1\right)\right)\right)f\left(\s_2,\widetilde{\X}_{t_2}\left(u_2\right)\right)\right\|_1\leq C_{DP}\,\vartheta^\frac{M\,\delta}{1+\delta}
	\end{equation*}	
for any $M\in\N$ with some positive constant $C_{DP}<\infty$ being independent of $M$.
\end{lemma}

\begin{proof}
By an iterative application of H\"older's inequality, we obtain
	\begin{align*}
		&\left\|\left(f\left(\s_1,\widetilde{\X}_{t_1}\left(u_1\right)\right)-f\left(\s_1,\widetilde{\X}_{t_1}^{(M)}\left(u_1\right)\right)\right)f\left(\s_2,\widetilde{\X}_{t_2}\left(u_2\right)\right)\right\|_1\\
		&\leq C\, \left\|f\left(\s_1,\widetilde{\X}_{t_1}\left(u_1\right)\right)-f\left(\s_1,\widetilde{\X}_{t_1}^{(M)}\left(u_1\right)\right)\right\|_{\frac{2+\delta}{1+\delta}}\\
	&	\leq C\,  \left(E\left|f\left(\s_1,\widetilde{\X}_{t_1}\left(u_1\right)\right)-f\left(\s_1,\widetilde{\X}_{t_1}^{(M)}\left(u_1\right)\right)\right|\right)^\frac{\delta}{1+\delta}\nonumber\\
		&\leq C \,\vartheta^\frac{M\,\delta}{1+\delta},
	\end{align*}
where the latter inequality follows from 
 Lemma~\ref{le: AbschM}.
\end{proof}\medskip

\begin{lemma}[Covariance Bounds]\label{le: KovAb}\ \\
	Suppose Assumption~\ref{as: Pro} holds true for $k=0$ and Assumption~\ref{as: Fu} is valid. Then, we have for some $C_{Cov}<\infty$ and every $h\in\Z$
	\begin{equation*}
\sup_{\s_1,\s_2\in\mathcal S, \,u_1,u_2\in [0,1]}	\left|\cov\left(f\left(\s_1,\widetilde{\X}_0\left(u_1\right)\right),f\left(\s_2,\widetilde{\X}_h\left(u_2\right)\right)\right)\right|\leq C_{Cov}\,\vartheta^\frac{|h|\delta}{2(1+\delta)}.
	\end{equation*}
\end{lemma}

\begin{proof}
	W.l.o.g.~consider $h\geq 1$ and let $M:=\left\lceil h/2\right\rceil$. We split the covariance up	and get
	\begin{align*}
		&\left|\cov\left(f\left(\s_1,\widetilde{\X}_0\left(u_1\right)\right),f\left(\s_2,\widetilde{\X}_h\left(u_2\right)\right)\right)\right|\\
		&\begin{multlined}[t][\linewidth]
			\leq \left|\cov\left(f\left(\s_1,\widetilde{\X}_0\left(u_1\right)\right)-f\left(\s_1,\widetilde{\X}_0^{(M)}\left(u_1\right)\right),f\left(\s_2,\widetilde{\X}_h\left(u_2\right)\right)\right)\right|\\
			+\left|\cov\left(f\left(\s_1,\widetilde{\X}_0^{(M)}\left(u_1\right)\right),f\left(\s_2,\widetilde{\X}_h\left(u_2\right)\right)-f\left(\s_2,\widetilde{\X}_h^{(M)}\left(u_2\right)\right)\right)\right|
		\end{multlined}
\\
		&=:\text{I}+\text{II}.
	\end{align*}
Obviously, it suffices to verify the result for the first summand.
 We start by splitting term I up further using Lemmata~\ref{le: AbschM}	and \ref{le: DiffProd}                
		\begin{align*}
		I&\begin{multlined}[t][0.99\linewidth]
			\leq
			\left|E\left(\left[f\left(\s_1,\widetilde{\X}_0\left(u_1\right)\right)-f\left(\s_1,\widetilde{\X}_0^{(M)}\left(u_1\right)\right)\right]f\left(\s_2,\widetilde{\X}_h\left(u_2\right)\right)\right)\right|\\
				+E\left|f\left(\s_1,\widetilde{\X}_0\left(u_1\right)\right)-f\left(\s_1,\widetilde{\X}_0^{(M)}\left(u_1\right)\right)\right|E\left|f\left(\s_2,\widetilde{\X}_h\left(u_2\right)\right)\right|
		\end{multlined}
		\\
			&\leq C_{DP}\,\vartheta^\frac{M\,\delta}{1+\delta}\, \,+\, C_{trunc}\vartheta^M,  
		\end{align*}
	which can obviously be bounded by $C\,\vartheta^\frac{h\delta}{2(1+\delta)}$.
\end{proof} \medskip

\begin{lemma}[Convergence of covariances]\label{le: Kov}\ \\
	Suppose Assumption~\ref{as: Pro} is satisfied for $k=1$ and Assumptions  \ref{as: We} as well as \ref{as: Fu} hold. Then, for any $\s_1,\s_2\in\mathcal S$ we have
	\begin{equation*}
		\cov\biggl(\sum_{t=1}^Tw_{t,T}\,f\left(\s_1,\X_{t,T}\right),\sum_{t=1}^Tw_{t,T}\,f\left(\s_2,\X_{t,T}\right)\biggr)\Tinfty \mathbf V \left(\s_1,\s_2\right), 
	\end{equation*}
where $\mathbf V(\s_1,\s_2)$ is defined in Theorem~\ref{th: CLT}.
\end{lemma}

\begin{proof} 
First, we rewrite the covariance and obtain 
		\begin{multline*}
			\lim_{T\to\infty} \sum_{t_1,t_2=1}^T \cov\left(w_{t_1,T}\,f\left(\s_1,\X_{t_1,T}\right),w_{t_2,T}\,f\left(\s_2,\X_{t_2,T}\right)\right)\\
			%
			%
			%
			%
			=\lim_{T\to\infty} \sum_{h=-(T-1)}^{T-1}\sum_{t=\max\{1,1-h\}}^{\min\{T,T-h\}}w_{t,T}\,w_{t+h,T}\, \cov\left(f\left(\s_1,\X_{t,T}\right),f\left(\s_2,\X_{t+h,T}\right)\right).
		\end{multline*}
		Note that 
		\begin{equation*}
		\sup_t	\left|\cov\left(f\left(\s_1, \X_{t,T}\right),f\left(\s_2, \X_{t+h,T}\right)\right)\right|\leq C\,\rho^{\frac{|h|\delta}{2+2\delta}}	
		\end{equation*}
	can be verified along the lines of the proof of Lemma~\ref{le: KovAb}. Hence, these covariances are absolutely summable w.r.t.~$h$. Due to the assumptions regarding the weights, Lebesgue's dominated convergence theorem is applicable, and we obtain
		\begin{multline}\label{eq: Kov2}
			\lim_{T\to\infty} \sum_{h=-(T-1)}^{T-1}\sum_{t=\max\{1,1-h\}}^{\min\{T,T-h\}}w_{t,T}\,w_{t+h,T}\, \cov\left(f\left(\s_1,\X_{t,T}\right),f\left(\s_2,\X_{t+h,T}\right)\right)\\
			=\sum_{h\in\Z}\lim_{T\to\infty} \sum_{t=\max\{1,1-h\}}^{\min\{T,T-h\}}w_{t,T}\,w_{t+h,T}\, \cov\left(f\left(\s_1,\X_{t,T}\right),f\left(\s_2,\X_{t+h,T}\right)\right)\I_{|h|<T}.
		\end{multline}
	Next, we incorporate the companion process into the covariance terms
		\begin{align*}
			&\cov\left(f\left(\s_1,\X_{t,T}\right),f\left(\s_2,\X_{t+h,T}\right)\right)\\
		&\begin{multlined}[t][\linewidth]
			=\cov\left(f\left(\s_1,\X_{t,T}\right)-f\left(\s_1,\widetilde{\X}_t\left(\frac{t}{T}\right)\right),f\left(\s_2,\X_{t+h,T}\right)\right)\\
			+\cov\left(f\left(\s_1,\widetilde{\X}_t\left(\frac{t}{T}\right)\right),f\left(\s_2,\X_{t+h,T}\right)-f\left(\s_2,\widetilde{\X}_{t+h}\left(\frac{t+h}{T}\right)\right)\right)\\
			+\cov\left(f\left(\s_1,\widetilde{\X}_t\left(\frac{t}{T}\right)\right),f\left(\s_2,\widetilde{\X}_{t+h}\left(\frac{t+h}{T}\right)\right)\right)
		\end{multlined}
			\\
			&=:\text{I}+\text{II}+\text{III}
		\end{align*}
		and show asymptotic negligibility of I and II. Using Remark~\ref{le: AF}, we get 
		\begin{align*}
		\text{I}	&
				\leq E\left(\left|f\left(\s_1,\X_{t,T}\right)-f\left(\s_1,\widetilde{\X}_t\left(\frac{t}{T}\right)\right)\right|\ \left|f\left(\s_2,\X_{t+h,T}\right)\right|\right)
				+O(T^{-1}).
		\end{align*}
		Inspired by the proof of Lemma~\ref{le: KovAb}, we obtain for the first summand from above with double use of H\"older's inequality
		
		\begin{align*}
			&E\left(\left|f\left(\s_1,\X_{t,T}\right)-f\left(\s_1,\widetilde{\X}_t\left(\frac{t}{T}\right)\right)\right| \left|f\left(\s_2,\X_{t+h,T}\right)\right|\right)\\
			&\begin{multlined}[t][\linewidth]
				\leq C\Biggl(\left(E\left|f\left(\s_1,\X_{t,T}\right)-f\left(\s_1,\widetilde{\X}_t\left(\frac{t}{T}\right)\right)\right|\right)^\frac{\delta}{1+\delta}\;\left(E\left|f\left(\s_1,\X_{t,T}\right)-f\left(\s_1,\widetilde{\X}_t\left(\frac{t}{T}\right)\right)\right|^{{2+\delta}}\right)^\frac{1}{(1+\delta)(2+\delta)}\Biggr)
			\end{multlined}			
		\\
			&\leq C\,T^{-\frac{\delta}{1+\delta}}
		\end{align*}
and, hence, $\text{I}=\go\left(T^{-\frac{\delta}{1+\delta}}\right)$. Analogously, we have $\text{II}=\go\left(T^{-\frac{\delta}{1+\delta}}\right)$. 		
		Thus, we can write (\ref{eq: Kov2}) as follows:
		\begin{multline*}
			\sum_{h\in\Z}\lim_{T\to\infty} \sum_{t=\max\{1,1-h\}}^{\min\{T,T-h\}}w_{t,T}\,w_{t+h,T}\, \cov\left(f\left(\s_1,\X_{t,T}\right),f\left(\s_2,\X_{t+h,T}\right)\right)\I_{|h|<T}\\
			=\sum_{h\in\Z}\lim_{T\to\infty} \sum_{t=\max\{1,1-h\}}^{\min\{T,T-h\}}w_{t,T}\,w_{t+h,T}\, \cov\left(f\left(\s_1,\widetilde{\X}_0\left(\frac{t}{T}\right)\right),f\left(\s_2,\widetilde{\X}_h\left(\frac{t+h}{T}\right)\right)\right)\I_{|h|<T}.
		\end{multline*}
	Now we introduce some notational adaptions and set
		\begin{equation*}
			\widetilde{\X}_t(z)=\begin{cases}
				\widetilde{\X}_t(1), & z>1,\\
				\widetilde{\X}_t(0), & z<0,
			\end{cases}\quad\text{for all }t\in\Z\text{ and}\quad	w_{t,T}=\begin{cases}
			w_{T,T}, & t>T,\\
			w_{1,T}, & t<1.
		\end{cases}
		\end{equation*}
This allows us to consider
		\begin{align*}
			&\begin{multlined}[t][\linewidth]
				\biggl|\sum_{h\in\Z}\lim_{T\to\infty} \biggl\{\sum_{t=\max\{1,1-h\}}^{\min\{T,T-h\}}w_{t,T}\,w_{t+h,T}\, \cov\left(f\left(\s_1,\widetilde{\X}_0\left(\frac{t}{T}\right)\right),f\left(\s_2,\widetilde{\X}_h\left(\frac{t+h}{T}\right)\right)\right)\\
				-\sum_{t=1}^T w_{t,T}\,w_{t+h,T}\, \cov\left(f\left(\s_1,\widetilde{\X}_0\left(\frac{t}{T}\right)\right),f\left(\s_2,\widetilde{\X}_h\left(\frac{t+h}{T}\right)\right)\right)\biggr\}\I_{|h|<T}\biggr|.
			\end{multlined}
		\end{align*}
	Using Lemma~\ref{le: KovAb}, the absolute value of the term in curly brackets can be bounded by
	$ C\, d_T^{-1}\, h \vartheta^{h\delta/(2\delta+2)}$,
		which converges to 0 as $T\to\infty$ for fixed $h$. The next step will be to change the argument of $\widetilde{\X}_h(\cdot)$ to eliminate $h$. To this end, we show asymptotic negligibility of 		
		\begin{align*}
			\cov\left(f\left(\s_1,\widetilde{\X}_0\left(\frac{t}{T}\right)\right),f\left(\s_2,\widetilde{\X}_h\left(\frac{t+h}{T}\right)\right)-f\left(\s_2,\widetilde{\X}_h\left(\frac{t}{T}\right)\right)\right).
		\end{align*}
		 With the help of the Lipschitz condition \eqref{eq: Lip} in Assumption~\ref{as: Fu} and Remark~\ref{le: AP}(i), it holds
		\begin{equation*}
			\left|E\left(f\left(\s_2,\widetilde{\X}_h\left(\frac{t+h}{T}\right)\right)-f\left(\s_2,\widetilde{\X}_h\left(\frac{t}{T}\right)\right)\right)\right|
			\leq C_1\frac{|h|}{T}.
		\end{equation*}
		 Next, remembering the use of H\"older's inequality in the proof of Lemma~\ref{le: DiffProd}, we obtain 
		\begin{equation*}
			E\left(\left|f\left(\s_1,\widetilde{\X}_0\left(\frac{t}{T}\right)\right)\right|\left|f\left(\s_2,\widetilde{\X}_h\left(\frac{t+h}{T}\right)\right)-f\left(\s_2,\widetilde{\X}_h\left(\frac{t}{T}\right)\right)\right|\right)
			\leq C_2\left(\frac{|h|}{T}\right)^\frac{\delta}{1+\delta}.
		\end{equation*}
		Thus, we have
		\begin{equation*}
			\left|\cov\left(f\left(\s_1,\widetilde{\X}_0\left(\frac{t}{T}\right)\right),f\left(\s_2,\widetilde{\X}_h\left(\frac{t+h}{T}\right)\right)-f\left(\s_2,\widetilde{\X}_h\left(\frac{t}{T}\right)\right)\right)\right|
			\leq \go\left(\left(\frac{|h|}{T}\right)^\frac{\delta}{1+\delta}\right)+\go\left(\frac{|h|}{T}\right).
		\end{equation*}
	 For fixed $h$, this converges to 0 as $T\to\infty$. Summing up, we get
		\begin{multline*}
			\lim_{T\to\infty} \sum_{h=-(T-1)}^{T-1}\sum_{t=\max\{1,1-h\}}^{\min\{T,T-h\}}w_{t,T}\,w_{t+h,T}\, \cov\left(f\left(\s_1,\X_{t,T}\right),f\left(\s_2,\X_{t+h,T}\right)\right)\\
			=\sum_{h\in\Z}\lim_{T\to\infty} \sum_{t=1}^{T}w_{t,T}\,w_{t+h,T}\, \cov\left(f\left(\s_1,\widetilde{\X}_0\left(\frac{t}{T}\right)\right),f\left(\s_2,\widetilde{\X}_h\left(\frac{t}{T}\right)\right)\right).
		\end{multline*}
		Finally, we obtain existence of the RHS using Lemma~\ref{le: KovAb} again.
\end{proof}\medskip

\begin{lemma}\label{L: rest-FCLT}
Under the prerequisites of Theorem~\ref{th: FCLT} combined with Assumption~\ref{as: Fu2}(ii)(b), it holds 
\begin{equation*}
	\lim_{r\to 0}\limsup_{T\to\infty}\,P\Biggl(\underset{\substack{\s_1,\s_2\in\mathcal S\\ \rho\left(\s_1,\s_2\right)<r}}{\sup}\sum_{t=1}^{\mu_{T}}\sum_{i\in H_t} w_{i,T} \left|\f\left(\s_1,\widetilde{\X}_i^{(M)}\left(\frac{i}{T}\right)\right)-\f\left(\s_2,\widetilde{\X}_i^{(M)}\left(\frac{i}{T}\right)\right)\right|>\frac{\lambda}{9}\Biggr)=0
\end{equation*}
using the notation of the proof of Theorem~\ref{th: FCLT}.
\end{lemma}

\begin{proof}
	We bound the quantity of interest as in~\eqref{eq: EZD} and start by specifying our choice for $m$ and demanding both a lower and an upper bound for $r_{k_T}$. Here, we assume $m>1+\frac{2}{\delta(1+\delta)}$, and that the bounds for $r_{k_T}$ are of the following form: 
\begin{equation}\label{eq: defknB}
	d_T^{-\frac{(m-1)(1+\delta)}{2m}}\leq r_{k_{T}}\leq d_T^{-\frac{1}{m\delta}}.
\end{equation}
Again, by our choice of $m$, we guarantee for the existence of a $k_T\in\N$ fulfilling the aforesaid requirements.
%
%
%
%
%
%
%
%
%
At this point, we return to the three summands on the RHS of~(\ref{eq: EZD}). The treatment of the individual terms will be carried out analogously to part (a). On account of this, we start by establishing the necessary bounds for the application of Bernstein's inequality on term II. First, we constitute an upper bound for the variance of the inner sum of $\nu_T$. With $l:=\left|i_1-i_2\right|$ and the very same arguments as in (a), we obtain 
\begin{align}\label{eq: FCLT3B}
	 \var\left(\nu_T\left(\s_1,\s_2\right)\right) 
	%
	%
	%
	\leq C\,\rho\left(\s_1,\s_2\right)^{\frac{1}{1+\delta}}
\end{align}
 using the boundedness of $f$ appropriately.
Going back to term II of equation (\ref{eq: EZD}), we use Bernstein's inequality and the notation we introduced in (\ref{eq: Sumlambdak}) to get
\begin{equation*}
	\text{II}\nonumber
	\leq 2\sum_{k=1}^{k_T} D(k)\, D(k-1)\,\exp\left(-\frac{1}{2}\cdot\frac{\lambda_k^2}{V_{II,k}+\frac{\breve{M}\lambda_k}{3}}\right).
\end{equation*}
Here, using again the boundedness of $f$,  we can set
\begin{equation*}
	\breve{M}:=C\,d_T^{-\frac{m-1}{2m}}\geq\sup_{\substack{\s_1, \s_2\in\mathcal S\\t\in\{1,\dots,\mu_T\}}}\biggl|\sum_{i\in H_t} w_{i,T}\left( \f\left(\s_1,\widetilde{\X}_i^{(M)}\left(\frac{i}{T}\right)\right)-\f\left(\s_2,\widetilde{\X}_i^{(M)}\left(\frac{i}{T}\right)\right)\right)\biggr|
\end{equation*}
and, in view of \eqref{eq: FCLT3B},
\begin{align*}
		V_{II,k}&:= C\,r_k^{\frac{1}{1+\delta}} 
		\geq
		\underset{\substack{\s_1\in\mathcal F_k,\s_2\in\mathcal F_{k-1}\\ \rho\left(\s_1,\s_2\right)\leq 3\,r_k}}{\sup}\var\left(\nu_T\left(\s_1,\s_2\right)\right).
\end{align*}
Similarly to (\ref{eq: termII}), we obtain
\begin{align*}
	\text{II}
	\leq 2\sum_{k=1}^{k_T}\exp\left(2\log\left(D(k)\right)-\frac{1}{2}\cdot\frac{\lambda_k^2}{C\,r_k^{\frac{1}{1+\delta}}}\right)
	\leq  2\sum_{k\in\N}\,\exp\left(-\frac{\overline{C}}{2}\,2^{\frac{k}{2(1+\delta)}}r^{-\frac{1}{2(1+\delta)}}\right)
	\rull 0.
\end{align*}
By the same arguments as before, we get
\begin{equation*}
	V_{III}:= C\, r^{\frac{1}{1+\delta}}
	\geq \underset{\substack{\s_1,\s_2\in\mathcal F_0\\ \rho\left(\s_1,\s_2\right)\leq 3\,r}}{\sup}\var\biggl(\sum_{t=1}^{\mu_{T}}\sum_{i\in H_t} w_{i,T}\left( \f\left(\s_1,\widetilde{\X}_i^{(M)}\left(\frac{i}{T}\right)\right)-\f\left(\s_2,\widetilde{\X}_i^{(M)}\left(\frac{i}{T}\right)\right)\right)\biggr).
\end{equation*}
Using the same $\breve{M}$ as while treating term II, Bernstein's inequality applied to equation (\ref{eq: FCLT9}) gives us 
\begin{align*}
	&P\Biggl(\underset{\substack{\s_1,\s_2\in\mathcal F_0\\ \rho\left(\s_1,\s_2\right)\leq 3\,r}}{\sup} \biggl|\sum_{t=1}^{\mu_{T}}\sum_{i\in H_t} w_{i,T}\left( \f\left(\s_1,\widetilde{\X}_i^{(M)}\left(\frac{i}{T}\right)\right)-\f\left(\s_2,\widetilde{\X}_i^{(M)}\left(\frac{i}{T}\right)\right)\right)\biggr|>\frac{\lambda}{27}\Biggr)\\
	%
	%
	%
	%
	%
	%
	&\leq 2\,D^2(0) \exp\left(-\frac{1}{2}\cdot\frac{C\lambda^2}{V_{III}+\frac{\breve{M}C\lambda}{3}}\right)\\
	%
	%
	%
	%
	%
	%
	&\leq 2\exp\left(2\log\left(D(0)\right)-\overline{C}\,r^{-\frac{1}{1+\delta}}\right)\\
	&\rull 0.
\end{align*}
Now term I is left. As in part (a), by Markov's inequality it is sufficient to show
\begin{equation}\label{eq.ew-konvB}
	\lim_{r\to 0}\limsup_{T\to\infty}\;E\Biggl(\underset{\substack{\s_1,\s_2\in\mathcal S\\ \rho\left(\s_1,\s_2\right)\leq r_{k_T}}}{\sup}\left|\nu_T\left(\s_1,\s_2\right)\right|\Biggr)=0.
\end{equation}
Again, using the notation of the proof of Theorem~\ref{th: FCLT} (case a), we get
\begin{equation*}
	E\Biggl(\underset{\substack{\s_1,\s_2\in\mathcal S\\ \rho\left(\s_1,\s_2\right)\leq r_{k_T}}}{\sup}\left|\nu_T\left(\s_1,\s_2\right)\right|\Biggr)
	\leq 2\,E\Biggl(\underset{\substack{\s_1,\s_2\in\mathcal S\\ \rho\left(\s_1,\s_2\right)\leq r_{k_T}}}{\sup}\biggl|\sum_{t=1}^{\mu_{T}}\left(L_{t,T}^0\left(\s_1\right)-L_{t,T}^0\left(\s_2\right)\right)\biggr|\Biggr)
\end{equation*}
similar to (\ref{eq: FCLT10}). By  Hoeffding's inequality we obtain 
\begin{equation*}
	P\biggl(\biggl|\sum_{t=1}^{\mu_{T}}L_{t,T}^0\left(\s_1\right)-L_{t,T}^0\left(\s_2\right)\biggr|>\hat{\rho}_{T,2}\left(\s_1,\s_2\right)\eta\,\bigg|\,L_{1,T},\dots,L_{\mu_T,T}\biggr)
	\leq 2\exp\left(-\frac{\eta^2}{2}\right)
\end{equation*}
 for any $\s_1,\s_2\in\mathcal S$ and $\eta>0$ comparable to (\ref{eq: Hoeff}) with $\hat{\rho}_{T,2}$ as in (\ref{eq: DefSM}). Thus, we can check for (\ref{eq.ew-konvB}) with the use of a maximal inequality for sub-Gaussian processes again. To ease the following part, we use a slightly different semimetric as above. In order to get the said new semimetric, we notice that
\begin{equation*}
	\left(L_{t,T}\left(\s_1\right)-L_{t,T}\left(\s_2\right)\right)^2
	\leq 2^\frac{2-\delta}{2} \left|L_{t,T}\right|_{\infty}^\frac{2-\delta}{2} \left|L_{t,T}\right|_{\mathrm{Lip}}^\frac{2+\delta}{2}\;\rho\left(\s_1,\s_2\right)^\frac{2+\delta}{2}
\end{equation*}
holds for $\s_1, \s_2\in\mathcal S$, comparable to (\ref{normen}). Next, set
\begin{equation}\label{eq: FCLT10Ba}
	Q_T:= 2^\frac{2-\delta}{4}\biggl(\sum_{t=1}^{\mu_{T}} \left|L_{t,T}\right|_{\infty}^\frac{2-\delta}{2}\, \left|L_{t,T}\right|_{\mathrm{Lip}}^\frac{2+\delta}{2}\biggr)^{1/2}.
\end{equation}
We can use the bounds
\begin{equation*}
	\left|L_{t,T}\right|_{\infty}^{\frac{2-\delta}{2}}
	=\biggl(\sup_{\s\in\mathcal S}\biggl|\sum_{i\in H_t}w_{i,T}\,f\left(\s,\widetilde{\X}^{(M)}_i\left(\frac{i}{T}\right)\right)\biggr|\biggr)^{\frac{2-\delta}{2}}
	\leq C_1\,\kappa_T^{\frac{2-\delta}{2}}d_T^{-\frac{2-\delta}{4}}
\end{equation*}
and
\begin{align*}
	\left|L_{t,T}\right|_{\mathrm{Lip}}^{\frac{2+\delta}{2}}
	%
	%
	%
	%
	%
	%
	&\leq\biggl(\sum_{i\in H_t}w_{i,T}\,g\left( \widetilde{\X}^{(M)}_i\left(\frac{i}{T}\right)\right)\biggr)^{{\frac{2+\delta}{2}}} 
	\leq C_2\,d_T^{\frac{2+\delta}{4}}\biggl(\sum_{i\in H_t}g\left( \widetilde{\X}^{(M)}_i\left(\frac{i}{T}\right)\right)\biggr)^{{\frac{2+\delta}{2}}}
\end{align*}
to obtain  
\begin{equation}\label{eq: QT}
	Q_T
	%
	%
	%
	%
	%
	%
	\leq C_3 \,\Biggl(\kappa_T^{\frac{2-\delta}{2}}d_T\sum_{t=1}^{\mu_{T}} \left(\sum_{i\in H_t}g\left( \widetilde{\X}^{(M)}_i\left(\frac{i}{T}\right)\right)\right)^{{\frac{2+\delta}{2}}}\Biggr)^{1/2}.
\end{equation}
Now we can define the new semimetric using the definition of $Q_T$:
\begin{align*}
	\hat{\rho}_{T,2}\left(\s_1,\s_2\right)
	&\leq Q_T \;\rho\left(\s_1,\s_2\right)^{\frac{2+\delta}{4}}=:\breve{\rho}_{T}\left(\s_1,\s_2\right).
\end{align*}
 As in part (a), we have
\begin{equation}\label{eq: FCLT11B}
	E\Biggl(\underset{\substack{\s_1,\s_2\in[-S,S]^d\\ \rho\left(\s_1,\s_2\right)\leq r_{k_T}}}{\sup}\biggl|\sum_{t=1}^{\mu_{T}}\left(L_{t,T}^0\left(\s_1\right)-L_{t,T}^0\left(\s_2\right)\right)\biggr|\,\Bigg|\,L_{1,T},\dots,L_{\mu_T,T}\Biggr)
	\leq C\int_{0}^{Q_T r_{k_T}^{\frac{2+\delta}{4}}}\left(\log\left(D\left(u,\mathcal S,\breve{\rho}_T\right)\right)\right)^{1/2}du.
\end{equation}
By Assumption~\ref{as: Fu2}, it holds
\begin{equation*}
	D\left(u,\mathcal S, \breve{\rho}_{T}\right)
	=D\left(\left(\frac{u}{Q_T}\right)^{\frac{4}{2+\delta}},\mathcal S, \rho\right)
	\leq C\,\left({\left(\frac{u}{Q_T}\right)^{-\frac{4}{2+\delta}}}+1\right)^d.
\end{equation*}
Now we can insert this bound into (\ref{eq: FCLT11B}) and obtain
\begin{equation}\label{eq: int}
	E\Biggl(\underset{\substack{\s_1,\s_2\in\mathcal S\\ \rho\left(\s_1,\s_2\right)\leq r_{k_T}}}{\sup}\biggl|\sum_{t=1}^{\mu_{T}}\left(L_{t,T}^0\left(\s_1\right)-L_{t,T}^0\left(\s_2\right)\right)\biggr|\Biggr)
	\leq C\,E[Q_T]\int_{0}^{r_{k_T}^{\frac{2+\delta}{4}}}u^{-\frac{2}{2+\delta}}\,du\,\leq\,  C\,E[Q_T]\,r_{k_T}^{\delta/4} .
\end{equation}
Next, we focus on $E[Q_T]$ and get with the use of (\ref{eq: QT})
\begin{equation}\label{eq: EQTB}
	E[Q_T]
	\leq C_3\left(\kappa_T^{\frac{2-\delta}{2}}d_T^{-1}\sum_{t=1}^{\mu_{T}} \biggl\|\sum_{i\in H_t}g\left( \widetilde{\X}^{(M)}_i\left(\frac{i}{T}\right)\right)\biggr\|_{{\frac{2+\delta}{2}}}^{{\frac{2+\delta}{2}}}\right)^{1/2}
	%
	%
	%
	%
	%
	\leq C_3\,d_T^{\frac{1}{4m}}.
\end{equation}
Finally, combining (\ref{eq: int})  and (\ref{eq: EQTB}) and using the upper bound of $r_{k_T}$ in (\ref{eq: defknB}) we have
\begin{equation*}
	E\Biggl(\underset{\substack{\s_1,\s_2\in\mathcal S\\ \rho\left(\s_1,\s_2\right)\leq r_{k_T}}}{\sup}\biggl|\sum_{t=1}^{\mu_{T}}\left(L_{t,T}^0\left(\s_1\right)-L_{t,T}^0\left(\s_2\right)\right)\biggr|\Biggr)
	%
	%
	%
	%
	%
	%
	=C\,d_T^{-\frac{1}{4m}},
\end{equation*}
which tends to 0 as $T\to\infty$, and the proof is completed.
\end{proof}\medskip

\subsection{Auxiliary results related to Section 3}

\begin{lemma}\label{le: KovEP}\ \\
	Assume the validity of Assumption~\ref{as: Ge}. Then, it holds  
 	\begin{equation*}
		\sup_{s\in\mathcal S}\sup_{t\leq T}E\left[E^\star |f\left(\s,\X_{t,T}^\star\right)|^{2+\delta}\right]<\infty.
	\end{equation*}
\end{lemma}

\begin{proof}
	We distinguish whether $t$ is an endpoint or not. 
	Starting with $t\notin EP$, we have
	\begin{equation*}
	\sup_{s\in\mathcal S}\sup_{t\leq T}	E\left[E^\star |f\left(\s,\X_{t,T}^\star\right)|^{2+\delta}\right]=\sup_{s\in\mathcal S}\sup_{t\leq T}\frac{1}{2\,TD_T+1}\sum_{r=-TD_T}^{TD_T}E|f\left(\s,\X_{t+r,T}\right)|^{2+\delta}<\infty	\end{equation*}
	by Assumption~\ref{as: Fu}. The cases $t\in EP_1$ and $t\in EP_2$ can be treated similarly, and we only examine the latter case further. There, we obtain
	\begin{equation*}
		E^\star |f\left(\s,\X_{t,T}^\star\right)|^{2+\delta}=\frac{1}{2\,TD_T+1}\biggl(\sum_{r=-TD_T}^{T-t}|f\left(\s,\X_{t+r,T}\right)|^{2+\delta}+\sum_{r=T-t+1}^{TD_T}|f\left(\s,\X_{t-r,T}\right)|^{2+\delta}\biggr).
	\end{equation*}
	Consequently, the desired result follows with the same arguments as above.        
\end{proof}\medskip

\begin{lemma}\label{le: ProdFu}\ \\
	Under Assumption~\ref{as: Ge}, we have for all $t_1,t_2\in\{1,\dots,T\}$ and $\s\in\mathcal S$
	\begin{equation*}
		f\left(\s,\X_{t_1,T}\right)f\left(\s,\X_{t_2,T}\right)=f\left(\s,\widetilde{\X}_{t_1}\left(\frac{t_1}{T}\right)\right)f\left(\s,\widetilde{\X}_{t_2}\left(\frac{t_2}{T}\right)\right)+\go_P\left(T^{-\frac{3+\delta}{2(4+\delta)}}\right).
	\end{equation*}
	Here, the $\go_P$-term does not depend on the choices for $t_1$ and $t_2$ or $\s$.  
\end{lemma}

\begin{proof}
	We want to make use of the closeness between the process $\left(\X_{t,T}\right)_{t=1}^T$ and its companion process $\left(\widetilde{\X}_t\left(\frac{t}{T}\right)\right)_{t\in\Z}$. So the first step is to rewrite the difference between the two products:
	\begin{align*}
		&E\left|f\left(\s,\X_{t_1,T}\right)f\left(\s,\X_{t_2,T}\right)-f\left(\s,\widetilde{\X}_{t_1}\left(\frac{t_1}{T}\right)\right)f\left(\s,\widetilde{\X}_{t_2}\left(\frac{t_2}{T}\right)\right)\right|\\
		%
		%
		%
		&	\leq E\left|\left(f\left(\s,\X_{t_1,T}\right)-f\left(\s,\widetilde{\X}_{t_1}\left(\frac{t_1}{T}\right)\right)\right)f\left(\s,\X_{t_2,T}\right)\right|
			+E\left|f\left(\s,\widetilde{\X}_{t_1}\left(\frac{t_1}{T}\right)\right)\left(f\left(\s,\X_{t_2,T}\right)-f\left(\s,\widetilde{\X}_{t_2}\left(\frac{t_2}{T}\right)\right)\right)\right|\\
		%
		%
		&=:\text{I}+\text{II}.
	\end{align*}
	Starting with I, we get using Remark~\ref{le: AF}
	\begin{align*}
		\text{I}
		%
		%
		%
		%
		%
		%
		&\leq C\left(E\left|f\left(\s,\X_{t_1,T}\right)-f\left(\s,\widetilde{\X}_{t_1}\left(\frac{t_1}{T}\right)\right)\right|\right)^\frac{3+\delta}{2(4+\delta)}\;\left(E\left|f\left(\s,\X_{t_1,T}\right)-f\left(\s,\widetilde{\X}_{t_1}\left(\frac{t_1}{T}\right)\right)\right|^\frac{2(5+\delta)}{6+\delta}\right)^\frac{3+\delta}{2(4+\delta)}\\
		&\leq   C\,{T^{-\frac{3+\delta}{2(4+\delta)}}}
	\end{align*}
	and, in complete analogy,   $\text{II}=\go\left({T^{\frac{3+\delta}{2(4+\delta)}}}\right)$. 
\end{proof}\medskip

\begin{lemma}\label{le: ChangeArg}\ \\
	Suppose Assumptions \ref{as: Ge} and \ref{as: Ra} are fulfilled. Then, for all $h,k\in\{1,\dots,T\}$, $-TD_T\leq r,l\leq TD_T$ and $\s\in\mathcal S$ it holds
	\begin{equation*}
		f\left(\s,\widetilde{\X}_{h+r}\left(\frac{h+r}{T}\right)\right)=f\left( \s,\widetilde{\X}_{h+r}\left(\frac{h}{T}\right)\right)+\go_P\left(D_T\right)
	\end{equation*}
	and
	$$
		f\left(\s,\widetilde{\X}_{h+r}\left(\frac{h+r}{T}\right)\right)f\left(\s,\widetilde{\X}_{k+l}\left(\frac{k+l}{T}\right)\right) 
		=f\left(\s,\widetilde{\X}_{h+r}\left(\frac{h}{T}\right)\right)f\left(\s,\widetilde{\X}_{k+l}\left(\frac{k}{T}\right)\right)+\go_P\left(D_T^\frac{3+\delta}{2(4+\delta)}\right).
	$$
	Both $\go_P$-terms are independent of the choices for $h,k,r,l$ and $\s$.
\end{lemma}

\begin{proof}
	\begin{enumerate}[(i)]
		\item Applying  Remark \ref{le: AP}(ii) after using the Lipschitz condition (\ref{eq: Lip}), we get part (i) straightforwardly as it holds for $-TD_T\leq r\leq TD_T$
		\begin{align*}
			&E\left|f\left( \s,\widetilde{\X}_{h+r}\left(\frac{h+r}{T}\right)\right)-f\left( \s,\widetilde{\X}_{h+r}\left(\frac{h}{T}\right)\right)\right| \leq C_{Lip}\frac{|r|}{T}
			\leq C\,D_T.
		\end{align*}
		\item We split up
		\begin{align*}
			&E\left|f\left(\s,\widetilde{\X}_{h+r}\left(\frac{h+r}{T}\right)\right)f\left(\s,\widetilde{\X}_{k+l}\left(\frac{k+l}{T}\right)\right)-f\left(\s,\widetilde{\X}_{h+r}\left(\frac{h}{T}\right)\right)f\left(\s,\widetilde{\X}_{k+l}\left(\frac{k}{T}\right)\right)\right|\\
			&\begin{multlined}[t][\linewidth]
				\leq E\left|\left(f\left( \s,\widetilde{\X}_{h+r}\left(\frac{h+r}{T}\right)\right)-f\left( \s,\widetilde{\X}_{h+r}\left(\frac{h}{T}\right)\right)\right)f\left(\s,\widetilde{\X}_{k+l}\left(\frac{k}{T}\right)\right)\right|\\
				+E\left|f\left( \s,\widetilde{\X}_{h+r}\left(\frac{h+r}{T}\right)\right)\left(f\left( \s,\widetilde{\X}_{k+l}\left(\frac{k+l}{T}\right)\right)-f\left(\s,\widetilde{\X}_{k+l}\left(\frac{k}{T}\right)\right)\right)\right|
			\end{multlined}\\
			&=:\text{I}+\text{II}.
		\end{align*}
		It suffices to investigate I. Using H\"older's inequality, we get 
		\begin{align*}
				\text{I}&\leq\left\|f\left(\s,\widetilde{\X}_{h+r}\left(\frac{h+r}{T}\right)\right)-f\left( \s,\widetilde{\X}_{h+r}\left(\frac{h}{T}\right)\right)\right\|_\frac{4+\delta}{3+\delta}\left\|f\left(\s,\widetilde{\X}_{k+l}\left(\frac{k}{T}\right)\right)\right\|_{4+\delta}
				\\
				&\begin{multlined}[t][0.99\linewidth]
					\leq C\,\left(E\left(\left|f\left(\s,\widetilde{\X}_{h+r}\left(\frac{h+r}{T}\right)\right)-f\left( \s,\widetilde{\X}_{h+r}\left(\frac{h}{T}\right)\right)\right|^\frac{1}{2}\right.\right.\\
					\left.\left.\cdot\left|f\left(\s,\widetilde{\X}_{h+r}\left(\frac{h+r}{T}\right)\right)-f\left( \s,\widetilde{\X}_{h+r}\left(\frac{h}{T}\right)\right)\right|^{\frac{4+\delta}{3+\delta}-\frac{1}{2}}\right)\right)^\frac{3+\delta}{4+\delta}
				\end{multlined} 
				\\
				&				\leq C\left(\frac{|r|}{T}\right)^\frac{3+\delta}{2(4+\delta)}\,\left(E\left|f\left(\s,\widetilde{\X}_{h+r}\left(\frac{h+r}{T}\right)\right)-f\left( \s,\widetilde{\X}_{h+r}\left(\frac{h}{T}\right)\right)\right|^\frac{2(5+\delta)}{6+\delta}\right)^\frac{3+\delta}{2(4+\delta)}\\
				&	=\go\left(D_T^\frac{3+\delta}{2(4+\delta)}\right).
		\end{align*}
	\end{enumerate}
\end{proof}\medskip

\begin{lemma}[Product Covariance Bound I]\label{le: KovVier}\ \\
	Suppose Assumption~\ref{as: Ge} is satisfied. Then, for all $u_1,\dots, u_4\in[0,1]$, $\s\in\mathcal S$
	\begin{enumerate}[(i)]
		\item and $t_1\in\N$ and $t_2,r\in\N_0$ fulfilling $t_1>t_2$ we have for some $\widetilde\rho\in (0,1)$
		\begin{equation*}
			\left|\cov\left(\f\left(\s,\widetilde\X_{t_1+r}(u_1)\right)\f\left(\s,\widetilde\X_{t_1}(u_2)\right),\f\left(\s,\widetilde\X_{t_2}(u_3)\right)\f\left(\s,\widetilde\X_0(u_4)\right)\right)\right| 
			\leq {C_{Cov,2i}}\;\widetilde \rho^{t_1-t_2}.
		\end{equation*}
		\item and $t_1,t_2\in\N$ with $t_1<t_2$ it holds for some $\widetilde\rho\in (0,1)$
		\begin{equation*}
			\left|\cov\left(\f\left(\s,\widetilde\X_{t_1+t_2}(u_1)\right)\f\left(\s,\widetilde\X_{t_1}(u_2)\right),\f\left(\s,\widetilde\X_{t_2}(u_3)\right)\f\left(\s,\widetilde\X_0(u_4)\right)\right)\right|
		\leq C_{Cov,2ii}\;\max \{\widetilde \rho^{t_2-t_1},\widetilde \rho^{t_1}\}
			.
		\end{equation*}
	\end{enumerate}
Here, the constants $C_{Cov,2i},\;C_{Cov,2ii}<\infty$ are independent of $\s$, $u_1,\dots,u_4$ as well as of $t_1,t_2$ and $r$.
\end{lemma}

{This lemma's proof uses the abbreviation
\begin{equation}\label{eq: CentM}
	\f_M\left(\s,\widetilde\X_t(u)\right):=f\left(\s,\widetilde\X_t^{(M)}(u)\right)-Ef\left(\s,\widetilde\X_t(u)\right)
\end{equation}
for $t\in\Z$, $\s\in\R^d$, $u\in[0,1]$ and some $M\in\N$.} 
\medskip
\begin{proof}
	\begin{enumerate}[(i)]
		\item Since $t_1+r\geq t_1>t_2>0$, we set $M:=\left\lceil\frac{t_1-t_2}{2}\right\rceil$. Now we are going to follow the proof of Lemma~\ref{le: KovAb}  but modified for the covariance of products. Hence, we start by inserting the truncated version of the companion process as introduced in (\ref{eq: Trun1}) with the above defined truncation parameter $M$ and obtain
		\begin{align}\label{eq: KovVier1}
			&\left|\cov\left(\f\left(\s,\widetilde\X_{t_1+r}(u_1)\right)\f\left(\s,\widetilde\X_{t_1}(u_2)\right),\f\left(\s,\widetilde\X_{t_2}(u_3)\right)\f\left(\s,\widetilde\X_0(u_4)\right)\right)\right| \nonumber\\
			&\begin{multlined}[t][\linewidth]
				\leq\left|\cov\left(\f\left(\s,\widetilde\X_{t_1+r}(u_1)\right)\f\left(\s,\widetilde\X_{t_1}(u_2)\right),\right.\right.\\
				\left.\left.\f\left(\s,\widetilde\X_{t_2}(u_3)\right)\f\left(\s,\widetilde\X_0(u_4)\right)-\f_M\left(\s,\widetilde\X_{t_2}(u_3)\right)\f_M\left(\s,\widetilde\X_0(u_4)\right)\right)\right|\\
				+\left|\cov\left(\f\left(\s,\widetilde\X_{t_1+r}(u_1)\right)\f\left(\s,\widetilde\X_{t_1}(u_2)\right)-\f_M\left(\s,\widetilde\X_{t_1+r}(u_1)\right)\f_M\left(\s,\widetilde\X_{t_1}(u_2)\right),\right.\right.\\
				\left.\left.
				\f\left(\s,\widetilde\X_{t_2}(u_3)\right)\f\left(\s,\widetilde\X_0(u_4)\right)\right)\right|
				%
				%
			\end{multlined}\nonumber\\
			&=:\text{I}+\text{II}
		\end{align}
		%
	 using (\ref{eq: CentM}).	As the terms I and II have a similar structure, we focus on term I and split up further
		\begin{align}\label{eq: KovVier2}
			\text{I}
			&\begin{multlined}[t][0.99\linewidth]
				\leq\left|\cov\left(\f\left(\s,\widetilde\X_{t_1+r}(u_1)\right)\f\left(\s,\widetilde\X_{t_1}(u_2)\right)\,,\,\f\left(\s,\widetilde\X_{t_2}(u_3)\right)\left(\f\left(\s,\widetilde\X_0(u_4)\right)-\f_M\left(\s,\widetilde\X_0(u_4)\right)\right)\right)\right|\\
				+\left|\cov\left(\f\left(\s,\widetilde\X_{t_1+r}(u_1)\right)\f\left(\s,\widetilde\X_{t_1}(u_2)\right)\,\,\f_M\left(\s,\widetilde\X_0(u_4)\right)\left(\f\left(\s,\widetilde\X_{t_2}(u_3)\right)-\f_M\left(\s,\widetilde\X_{t_2}(u_3)\right)\right)\right)\right|
			\end{multlined}\nonumber\\
			&=:\text{Ia}+\text{Ib}.
		\end{align}
		For the same reasons as above, we limit ourselves to the investigation of the first subterm of (\ref{eq: KovVier2}) and obtain by an iterative application of H\"older's inequality and Lemma~\ref{le: AbschM}
		\begin{align*}
			\text{Ia}
			&\begin{multlined}[t][0.98\linewidth]
				\leq E\left|\f\left(\s,\widetilde\X_{t_2+r}(u_1)\right)\f\left(\s,\widetilde\X_{t_1}(u_2)\right)\f\left(\s,\widetilde\X_{t_2}(u_3)\right) \,\left[\f\left(\s,\widetilde\X_0(u_4)\right)-\f_M\left(\s,\widetilde\X_0(u_4)\right)\right]\right|\\
				+C\,  E\left|\f\left(\s,\widetilde\X_{t_2}(u_3)\right)\left(\f\left(\s,\widetilde\X_0(u_4)\right)-\f_M\left(\s,\widetilde\X_0(u_4)\right)\right)\right|
			\end{multlined}\\
			&\leq C\, 
			\left\|\f\left(\s,\widetilde\X_0(u_4)\right)-\f_M\left(\s,\widetilde\X_0(u_4)\right)\right\|_{\frac{4+\delta}{1+\delta}}
		 \\
		&\leq C\,\left(E\left|f\left(\s,\widetilde\X_0(u_4)\right)-f\left(\s,\widetilde\X_0^{(M)}(u_4)\right)\right|\right)^\frac{\delta}{(3+\delta)}
		 \\
			&\leq C\left(\rho^M\right)^\frac{\delta}{3+\delta}.
		\end{align*}
		In complete analogy, we get the same result as upper bound for Ib. In conclusion, this upper bound is valid for term I and II in (\ref{eq: KovVier1}) as well. Setting $\widetilde\rho= \rho^\frac{\delta}{ 2(3+\delta)}$ finishes the proof. 
		%
	
		\item Applying part (i) and Lemma~\ref{le: KovAb} leads to
		\begin{align*}
		&	\left|\cov\left(\f\left(\s,\widetilde\X_{t_1+t_2}(u_1)\right)\f\left(\s,\widetilde\X_{t_1}(u_2)\right),\f\left(\s,\widetilde\X_{t_2}(u_3)\right)\f\left(\s,\widetilde\X_0(u_4)\right)\right)\right|
		\\
		&\begin{multlined}[t][\linewidth]
			\leq \left|\cov\left(\f\left(\s,\widetilde\X_{t_1+t_2}(u_1)\right)\f\left(\s,\widetilde\X_{t_2}(u_3)\right),\f\left(\s,\widetilde\X_{t_1}(u_2)\right)\f\left(\s,\widetilde\X_0(u_4)\right)\right)\right|\\
			+\left|\cov\left(\f\left(\s,\widetilde\X_{t_1+t_2}(u_1)\right)\f\left(\s,\widetilde\X_{t_1}(u_2)\right)\right)\, \cov\left(\f\left(\s,\widetilde\X_{t_2}(u_3)\right),\,\f\left(\s,\widetilde\X_0(u_4)\right)\right)\right|\\
			+\left|\cov\left(\f\left(\s,\widetilde\X_{t_1+t_2}(u_1)\right),\f\left(\s,\widetilde\X_{t_2}(u_3)\right)\right)\, \cov\left(\f\left(\s,\widetilde\X_{t_1}(u_2)\right),\,\f\left(\s,\widetilde\X_0(u_4)\right)\right)\right|
		\end{multlined}\\
	&\leq C\max\{\widetilde \rho^{t_2-t_1},\widetilde \rho^{t_1}\}
		\end{align*}
		for some $\widetilde \rho\in(0,1)$.
	\end{enumerate}
\end{proof}\medskip
\begin{lemma}\label{le: KovBSKov}\ \\
	Suppose Assumptions \ref{as: Ge} and \ref{as: Ra} are valid. Then, for all indices $t_1,t_2\in\{1,\dots,T\}\backslash EP$ 
	as well as for all $\s\in\mathcal S$ and $u_1,\, u_2\in[0,1]$ it holds
	\begin{align*}
		&\begin{multlined}[t][\linewidth]
			\frac{1}{2\,TD_T+1}\sum_{r=-TD_T}^{TD_T}\biggl(\biggl(f\left(\s,\widetilde{\X}_{t_1+r}(u_1)\right)-\frac{1}{2\,TD_T+1}\sum_{l=-TD_T}^{TD_T}f\left(\s,\widetilde{\X}_{t_1+l}(u_1)\right)\biggr)\\
			\cdot \biggl(f\left(\s,\widetilde{\X}_{t_2+r}(u_2)\right)-\frac{1}{2\,TD_T+1}\sum_{l=-TD_T}^{TD_T}f\left(\s,\widetilde{\X}_{t_2+l}(u_2)\right)\biggr)\biggr)
		\end{multlined}\\	
		&=\cov\left(f\left(\s,\widetilde{X}_{t_1}(u_1)\right),f\left(\s,\widetilde{X}_{t_2}(u_2)\right)\right)+\go_P\left(\left(TD_T\right)^{-1/2}\right).
	\end{align*}
	At that, the occurring $\go_P$-term is independent of $t_1,\,t_2,\,\s$ and $u_1,\, u_2$.
\end{lemma}	
\begin{proof}
	First, note that by Lemma~\ref{le: KovAb} we have
	\begin{equation*}
		(2\,TD_T+1)^{-1}\sum_{l=-TD_T}^{TD_T}\bar f (\s,\widetilde{\X}_{t_2+l}(u))=\go_P((TD_T)^{-1/2})
	\end{equation*}
	for all $u\in[0,1]$. Thus, we obtain
	\begin{align*}
		&\begin{multlined}[t][\linewidth]
			\frac{1}{2\,TD_T+1}\sum_{r=-TD_T}^{TD_T}\biggl(\biggl(f\left(\s,\widetilde{\X}_{t_1+r}(u_1)\right)-\frac{1}{2\,TD_T+1}\sum_{l=-TD_T}^{TD_T}f\left(\s,\widetilde{\X}_{t_1+l}(u_1)\right)\biggr)\\
			\cdot \biggl(f\left(\s,\widetilde{\X}_{t_2+r}(u_2)\right)-\frac{1}{2\,TD_T+1}\sum_{l=-TD_T}^{TD_T}f\left(\s,\widetilde{\X}_{t_2+l}(u_2)\right)\biggr)\biggr)
		\end{multlined}\\
		&=\frac{1}{2\,TD_T+1}\sum_{r=-TD_T}^{TD_T}\f\left(\s,\widetilde{\X}_{t_1+r}(u_1)\right)\f\left(\s,\widetilde{\X}_{t_2+r}(u_2)\right)+\go_P((TD_T)^{-1}).
	\end{align*}
	W.l.o.g.~we consider $v:=t_1-t_2>0$ only. It remains to bound
	\begin{align*}
			&E\biggl(\frac{1}{2\,TD_T+1}\sum_{r=-TD_T}^{TD_T}\f\left(\s,\widetilde{\X}_{t_1+r}(u_1)\right)\f\left(\s,\widetilde{\X}_{t_2+r}(u_2)\right)-\cov\left(f\left(\s,\widetilde{X}_{t_1}(u_1)\right),f\left(\s,\widetilde{X}_{t_2}(u_2)\right)\right)\biggr)^2\\
			&\leq \frac{1}{(2\,TD_T+1)}\sum_{t=-{2\, TD_T}}^{2\,TD_T}
			\cdot\left|\cov\left(\f\left(\s,\widetilde{\X}_{t+v}(u_1)\right)\f\left(\s,\widetilde{\X}_{t}(u_2)\right),\f\left(\s,\widetilde{\X}_{v}(u_1)\right)\f\left(\s,\widetilde{\X}_{0}(u_2)\right)\right)\right|\\
			&\begin{multlined}[t][\linewidth]
				\leq\frac{2}{2\,TD_T+1}\sum_{t=1}^{v-1}\left|\cov\left(\f\left(\s,\widetilde{\X}_{t+v}(u_1)\right)\f\left(\s,\widetilde{\X}_{t}(u_2)\right),\f\left(\s,\widetilde{\X}_{v}(u_1)\right)\f\left(\s,\widetilde{\X}_{0}(u_2)\right)\right)\right|\\
				+\frac{2}{2\,TD_T+1}\sum_{t=v+1}^{2\,TD_T}\left|\cov\left(\f\left(\s,\widetilde{\X}_{t+v}(u_1)\right)\f\left(\s,\widetilde{\X}_{t}(u_2)\right),\f\left(\s,\widetilde{\X}_{v}(u_1)\right)\f\left(\s,\widetilde{\X}_{0}(u_2)\right)\right)\right|\\
				+\go((TD_T)^{-1})
			\end{multlined} 	\\
			&=:\text{I}+\text{II}+\go((TD_T)^{-1}).
	\end{align*}
	Applying Lemma~\ref{le: KovVier}(i) to term I and Lemma~\ref{le: KovVier}(ii) to II, we obtain the desired result, namely $I+II=\go ((TD_T)^{-1})$.
\end{proof}\medskip
\begin{lemma}[$P$-Convergence of the Bootstrap Variance]\label{th: VarBS}\ \\
	Suppose Assumptions \ref{as: Ge} and \ref{as: Ra} are true. Then, we have for all $\s\in\mathcal S$
	\begin{equation*}
		\var^\star\biggl(\sum_{t=1}^{T} w_{t,T}\, f\left(\s, \X_{t,T}^\star\right)\biggr)\pto\mathbf V\left(\s,\s\right)
	\end{equation*}
	as $T\to\infty$, where $\mathbf V \left(\s,\s\right)$ has its roots in Theorem~\ref{th: CLT}.
\end{lemma}
\begin{proof}
	We start by dividing the quantity of interest  into three sums such that the middle one does not contain endpoints and only full bootstrap blocks:
	\begin{align}\label{eq: T3BS1}
		&\var^\star\biggl(\sum_{t=1}^{T} w_{t,T}\, f\left(\s, \X_{t,T}^\star\right)\biggr)\nonumber\\
		&\begin{multlined}[t][\linewidth]
		=\var^\star\biggl(\sum_{t=1}^{L_T\left\lceil\left(TD_T+1\right)/L_T\right\rceil} w_{t,T}\, f\left(\s, \X_{t,T}^\star\right)\biggr)
			+\var^\star\biggl(\sum_{t=L_T\left\lceil\left(TD_T+1\right)/L_T\right\rceil+1}^{L_T\left\lfloor\left(T-TD_T\right)/L_T\right\rfloor} w_{t,T}\, f\left(\s, \X_{t,T}^\star\right)\biggr)\\
			+\var^\star\biggl(\sum_{t=L_T\left\lfloor\left(T-TD_T\right)/L_T\right\rfloor+1}^{T} w_{t,T}\, f\left(\s, \X_{t,T}^\star\right)\biggr)
		\end{multlined}
		\nonumber\\
		&=:\text{I}+\text{II}+\text{III}.
	\end{align}
	Note that Assumption~\ref{as: Ra} and independence of the bootstrap blocks guarantee validity of this partition for any sufficiently large $T$. Since terms I and III from (\ref{eq: T3BS1}) are of a similar type, we focus on the first. From Lemma~\ref{le: KovEP} and the rates in Assumption~\ref{as: Ra}, we obtain
	\begin{equation}\label{eq: T3BS1a}
		\text{I}
		=\sum_{t_1,t_2=1}^{L_T\left\lceil\left(TD_T+1\right)/L_T\right\rceil}w_{t_1,T}\,w_{t_2,T}\,\cov^\star\left(f\left(\s, \X_{t_1,T}^\star\right),f\left(\s, \X_{t_2,T}^\star\right)\right)
		%
		%
		%
		%
		%
		%
		=\go_P\left(d_T^{-\frac{\delta}{2+\delta}}\right)= o_P(1).
	\end{equation}
 Thus, we only need to set our focus on the second bootstrap variance term of (\ref{eq: T3BS1}):
	\begin{multline*}
		\text{II}
		=\sum_{t=\left\lceil\left(TD_T+1\right)/L_T\right\rceil}^{\left\lfloor\left(T-TD_T\right)/L_T\right\rfloor-1}\sum_{j=1}^{L_T}\sum_{l=1}^{L_T} w_{tL_T+j,T}\,w_{tL_T+l,T}
	\biggl(\frac{1}{2\,TD_T+1}\sum_{r=-TD_T}^{TD_T}f\left(\s, \X_{tL_T+j+r,T}\right)f\left(\s, \X_{tL_T+l+r,T}\right)\\
		-\frac{1}{\left(2\,TD_T+1\right)^2}\sum_{r=-TD_T}^{TD_T}f\left(\s, \X_{tL_T+j+r,T}\right)\sum_{k=-TD_T}^{TD_T}f\left(\s, \X_{tL_T+l+k,T}\right)\biggr).
	\end{multline*}
	We aim for transforming the bootstrap covariance into the real world one with negligible error. To do so, the first step is to change the process $\left(\X_{t,T}\right)$ to the companion process with the aid of Lemma~\ref{le: ProdFu}:
	\begin{multline}\label{eq: T3BS3}
	\text{II}=	
		\sum_{t=\left\lceil\left(TD_T+1\right)/L_T\right\rceil}^{\left\lfloor\left(T-TD_T\right)/L_T\right\rfloor-1}\sum_{j=1}^{L_T}\sum_{l=1}^{L_T} w_{tL_T+j,T}\,w_{tL_T+l,T}\\
		\cdot\left(\frac{1}{2\,TD_T+1}\sum_{r=-TD_T}^{TD_T}f\left(\s, \widetilde{\X}_{tL_T+j+r}\left(\frac{tL_T+j+r}{T}\right)\right)f\left(\s, \widetilde{\X}_{tL_T+l+r}\left(\frac{tL_T+l+r}{T}\right)\right)\right.\\
			-\frac{1}{\left(2\,TD_T+1\right)^2}\biggl(\sum_{r=-TD_T}^{TD_T}f\left(\s, \widetilde{\X}_{tL_T+j+r}\left(\frac{tL_T+j+r}{T}\right)\right)\\
			\left.\cdot\sum_{k=-TD_T}^{TD_T}f\left(\s, \widetilde{\X}_{tL_T+l+k}\left(\frac{tL_T+l+k}{T}\right)\right)\biggr)\right) 
			+\mathcal{O}_P\left({L_T}{T^{-\frac{3+\delta}{2(4+\delta)}}}\right).
	\end{multline}
	Recalling that $L_T=o(d_T^{\frac{\delta}{2(1+\delta)}})$, straight-forward calculations give $\mathcal{O}_P\left({L_T}{T^{-\frac{3+\delta}{2(4+\delta)}}}\right)=o_P(1)$. Next, we change the argument of $\left(\widetilde{\X}_t(u)\right)$ so that it loses the dependence of the inner summation index. With the aid of Lemma~\ref{le: ChangeArg}, the first summand of (\ref{eq: T3BS3}) becomes
	\begin{align}\label{eq: T3BS3a}
		&\begin{multlined}[t][\linewidth]
			\sum_{t=\left\lceil\left(TD_T+1\right)/L_T\right\rceil}^{\left\lfloor\left(T-TD_T\right)/L_T\right\rfloor-1}\sum_{h=-(L_T-1)}^{L_T-1}\sum_{j=\max\{1,1-h\}}^{\min\{L_T,L_T-h\}} w_{tL_T+j+h,T}\,w_{tL_T+j,T}\\
			\cdot\left(\frac{1}{2\,TD_T+1}\sum_{r=-TD_T}^{TD_T}f\left(\s, \widetilde{\X}_{tL_T+j+h+r}\left(\frac{tL_T+j+h+r}{T}\right)\right)\right.f\left(\s, \widetilde{\X}_{tL_T+j+r}\left(\frac{tL_T+j+r}{T}\right)\right)\\
			-\frac{1}{\left(2\,TD_T+1\right)^2}\biggl(\sum_{r=-TD_T}^{TD_T}f\left(\s, \widetilde{\X}_{tL_T+j+h+r}\left(\frac{tL_T+j+h+r}{T}\right)\right)\\
			\left.\cdot\sum_{k=-TD_T}^{TD_T}f\left(\s, \widetilde{\X}_{tL_T+j+k}\left(\frac{tL_T+j+k}{T}\right)\right)\biggr)\right)
		\end{multlined}\nonumber\\
	&\begin{multlined}[t][\linewidth]
		=\sum_{t=\left\lceil\left(TD_T+1\right)/L_T\right\rceil}^{\left\lfloor\left(T-TD_T\right)/L_T\right\rfloor-1}\sum_{h=-(L_T-1)}^{L_T-1}\sum_{j=\max\{1,1-h\}}^{\min\{L_T,L_T-h\}} w_{tL_T+j+h,T}\,w_{tL_T+j,T}\\
		\cdot\left(\frac{1}{2\,TD_T+1}\sum_{r=-TD_T}^{TD_T}f\left(\s, \widetilde{\X}_{tL_T+j+h+r}\left(\frac{tL_T+j+h}{T}\right)\right)f\left(\s, \widetilde{\X}_{tL_T+j+r}\left(\frac{tL_T+j}{T}\right)\right)\right.\\
		-\frac{1}{\left(2\,TD_T+1\right)^2}\biggl(\sum_{r=-TD_T}^{TD_T}f\left(\s, \widetilde{\X}_{tL_T+j+h+r}\left(\frac{tL_T+j+h}{T}\right)\right)\\
		\left.\cdot\sum_{k=-TD_T}^{TD_T}f\left(\s, \widetilde{\X}_{tL_T+j+k}\left(\frac{tL_T+j}{T}\right)\right)\biggr)\right)
		+\go_P\left(L_TD_T^\frac{3+\delta}{2(4+\delta)}\right).
	\end{multlined}
\end{align}
Note that under Assumption~\ref{as: Ra}, it holds $\go_P\left(L_TD_T^\frac{3+\delta}{2(4+\delta)}\right)=o_P(1)$. Hence, we focus again on the first summand on the RHS of (\ref{eq: T3BS3a}) and obtain its equivalence to
	\begin{multline}\label{eq: T3BS5}
			\sum_{t=\left\lceil\left(TD_T+1\right)/L_T\right\rceil}^{\left\lfloor\left(T-TD_T\right)/L_T\right\rfloor-1}\sum_{h=-(L_T-1)}^{L_T-1}\sum_{j=\max\{1,1-h\}}^{\min\{L_T,L_T-h\}} w_{tL_T+j+h,T}\,w_{tL_T+j,T}\\
			\cdot \cov\left(f\left(\s, \widetilde{\X}_{h}\left(\frac{tL_T+j+h}{T}\right)\right),f\left(\s, \widetilde{\X}_{0}\left(\frac{tL_T+j}{T}\right)\right)\right)
			+\go_P\left(\frac{L_T}{\sqrt{TD_T }}\right)
	\end{multline}
	by Lemma~\ref{le: KovBSKov}. Comparably to the proof of Lemma~\ref{le: Kov}, we would like to rewrite the inner sum in order to eliminate the minimum and maximum determining the index bounds. Therefore, we show asymptotic negligibility of the difference between the first summand of (\ref{eq: T3BS5}) and
	\begin{multline}\label{eq: T3BS6}
		\sum_{t=\left\lceil\left(TD_T+1\right)/L_T\right\rceil}^{\left\lfloor\left(T-TD_T\right)/L_T\right\rfloor-1}\sum_{h=-(L_T-1)}^{L_T-1}\sum_{j=1}^{L_T} w_{tL_T+j+h,T}\,w_{tL_T+j,T}\\
		\cdot\cov\left(f\left(\s, \widetilde{\X}_{h}\left(\frac{tL_T+j+h}{T}\right)\right),f\left(\s, \widetilde{\X}_{0}\left(\frac{tL_T+j}{T}\right)\right)\right)
	\end{multline}
	$T$ tends to infinity. Note that (\ref{eq: T3BS6}) is well-defined for sufficiently large $T$ as we skipped endpoints within the summation. Regarding the difference in question, we get with the very same arguments as in the proof of Lemma~\ref{le: Kov}
	\begin{align*}
		&\begin{multlined}[t][\linewidth]
			\left|\sum_{t=\left\lceil\left(TD_T+1\right)/L_T\right\rceil}^{\left\lfloor\left(T-TD_T\right)/L_T\right\rfloor-1}\sum_{h=-(L_T-1)}^{L_T-1}\sum_{j=\max\{1,1-h\}}^{\min\{L_T,L_T-h\}} w_{tL_T+j+h,T}\,w_{tL_T+j,T}\right.\\
			\cdot\cov\left(f\left(\s, \widetilde{\X}_{h}\left(\frac{tL_T+j+h}{T}\right)\right),f\left(\s, \widetilde{\X}_{0}\left(\frac{tL_T+j}{T}\right)\right)\right)\\
			-\sum_{t=\left\lceil\left(TD_T+1\right)/L_T\right\rceil}^{\left\lfloor\left(T-TD_T\right)/L_T\right\rfloor-1}\sum_{h=-(L_T-1)}^{L_T-1}\sum_{j=1}^{L_T} w_{tL_T+j+h,T}\,w_{tL_T+j,T}\\
			\cdot\left.\cov\left(f\left(\s, \widetilde{\X}_{h}\left(\frac{tL_T+j+h}{T}\right)\right),f\left(\s, \widetilde{\X}_{0}\left(\frac{tL_T+j}{T}\right)\right)\right)\right|
		\end{multlined}\\
		&=\go\left(L_T^{-1}\right).
	\end{align*}
Thus, we can proceed with (\ref{eq: T3BS6}). The next step will be to incorporate the sum over $j$ back into the sum over $t$. Thereby, we obtain equivalence of (\ref{eq: T3BS6}) and
	\begin{equation}\label{eq: T3BS8}
		\sum_{h=-(L_T-1)}^{L_T-1}\sum_{t=L_T\left\lceil\left(TD_T+1\right)/L_T\right\rceil+1}^{L_T\left\lfloor\left(T-TD_T\right)/L_T\right\rfloor}w_{t+h,T}\,w_{t,T}\,\cov\left(f\left(\s, \widetilde{\X}_{h}\left(\frac{t+h}{T}\right)\right),f\left(\s, \widetilde{\X}_{0}\left(\frac{t}{T}\right)\right)\right)
	\end{equation}
	 for all sufficiently large $T$. By following exactly the lines of the proof of Lemma~\ref{le: Kov}, we can rewrite (\ref{eq: T3BS8}) as
	\begin{multline}\label{eq: T3BS9}
			\sum_{h=-(L_T-1)}^{L_T-1}\sum_{t=L_T\left\lceil\left(TD_T+1\right)/L_T\right\rceil+1}^{L_T\left\lfloor\left(T-TD_T\right)/L_T\right\rfloor}w_{t+h,T}\,w_{t,T}
			\left(\cov\left(f\left(\s, \widetilde{\X}_{h}\left(\frac{t}{T}\right)\right),f\left(\s, \widetilde{\X}_{0}\left(\frac{t}{T}\right)\right)\right)+\go\left(\left(\frac{|h|}{T}\right)^\frac{\delta}{1+\delta}\right)\right)\\
		%
		%
		%
			=\sum_{h=-(L_T-1)}^{L_T-1}\sum_{t=L_T\left\lceil\left(TD_T+1\right)/L_T\right\rceil+1}^{L_T\left\lfloor\left(T-TD_T\right)/L_T\right\rfloor}w_{t+h,T}\,w_{t,T}\,\cov\left(f\left(\s, \widetilde{\X}_{h}\left(\frac{t}{T}\right)\right),f\left(\s, \widetilde{\X}_{0}\left(\frac{t}{T}\right)\right)\right)
			+\go\left(\left(\frac{L_T}{T}\right)^\frac{\delta}{1+\delta}\right).
	\end{multline}
	Since the last term of (\ref{eq: T3BS9}) tends to 0 as $T\to \infty$, it remains to show
	\begin{equation}\label{eq: T3BS10}
		\lim_{T\to\infty}\sum_{h=-(L_T-1)}^{L_T-1}\sum_{t=L_T\left\lceil\left(TD_T+1\right)/L_T\right\rceil+1}^{L_T\left\lfloor\left(T-TD_T\right)/L_T\right\rfloor}w_{t+h,T}\,w_{t,T}\,\cov\left(f\left(\s, \widetilde{\X}_{h}\left(\frac{t}{T}\right)\right),f\left(\s, \widetilde{\X}_{0}\left(\frac{t}{T}\right)\right)\right) 
		=\mathbf V \left(\s,\s\right).
	\end{equation}
	 In view of Lebesgue's dominated convergence theorem, we obtain that the LHS of \eqref{eq: T3BS10} equals
	\begin{equation*}
		\sum_{h\in\Z}\lim_{T\to\infty}\sum_{t=L_T\left\lceil\left(TD_T+1\right)/L_T\right\rceil+1}^{L_T\left\lfloor\left(T-TD_T\right)/L_T\right\rfloor}w_{t+h,T}\,w_{t,T}\,\cov\left(f\left(\s, \widetilde{\X}_{h}\left(\frac{t}{T}\right)\right),f\left(\s, \widetilde{\X}_{0}\left(\frac{t}{T}\right)\right)\right).
	\end{equation*}
By definition of $V_h(\s,\s)$, it suffices to show asymptotic negligibility of 	
	\begin{align*}
	&\begin{multlined}[t][\linewidth]
		\sum_{h\in\Z}	\biggl|\sum_{t=1}^{T}w_{t+h,T}\,w_{t,T}\,\cov\left(f\left(\s, \widetilde{\X}_{h}\left(\frac{t}{T}\right)\right),f\left(\s, \widetilde{\X}_{0}\left(\frac{t}{T}\right)\right)\right)\\
		-\sum_{h\in\Z}\sum_{t=L_T\left\lceil\left(TD_T+1\right)/L_T\right\rceil+1}^{L_T\left\lfloor\left(T-TD_T\right)/L_T\right\rfloor}w_{t+h,T}\,w_{t,T}\,\cov\left(f\left(\s, \widetilde{\X}_{h}\left(\frac{t}{T}\right)\right),f\left(\s, \widetilde{\X}_{0}\left(\frac{t}{T}\right)\right)\right)\biggr|
	\end{multlined}	\\
&\begin{multlined}[t][\linewidth]
	\leq
	\sum_{h\in\Z}\sum_{t=1}^{L_T\left\lceil\left(TD_T+1\right)/L_T\right\rceil}w_{t+h,T}\,w_{t,T}\left|\cov\left(f\left(\s, \widetilde{\X}_{h}\left(\frac{t}{T}\right)\right),f\left(\s, \widetilde{\X}_{0}\left(\frac{t}{T}\right)\right)\right)\right|\\
	+\sum_{h\in\Z}\sum_{t=L_T\left\lfloor\left(T-TD_T\right)/L_T\right\rfloor+1}^{T}w_{t+h,T}\,w_{t,T}\left|\cov\left(f\left(\s, \widetilde{\X}_{h}\left(\frac{t}{T}\right)\right),f\left(\s, \widetilde{\X}_{0}\left(\frac{t}{T}\right)\right)\right)\right|,
\end{multlined}
		\end{align*}
which is, indeed, of order $\go\left(d_T^{-\frac{\delta}{2+\delta}}\right)$.
\end{proof}\medskip

\begin{lemma}\label{le: UB}\ \\
	Suppose Assumptions \ref{as: Ge} to \ref{as: G-BS} hold true. Then, $\lim_{T\to\infty}P(\mathbb X_T\in A_T)=1$, where
	\begin{equation*}
		A_T:=\Biggl\{\un x_T\in \R^{dT}\,\Bigg|\,\underset{\substack{\s\in\mathcal S\\ t=1,\dots,T}}{\sup}w_{t,T}\left|f\left(\s,\X^\star_{t,T}\right)\right|\leq d_T^{-\frac{2+\delta^2}{2(4+\delta)}}~\text{ for }~\mathbb X_T=\un x_T\Biggr\}.
	\end{equation*}
\end{lemma}
\begin{proof}[Proof] 
	Recalling the Lipschitz condition in Assumption~\ref{as: Fu2}, Assumption~\ref{as: We}, and the moment conditions on~$g$, we obtain asymptotic negligibility of
	\begin{equation*}
		P\left(\mathbb X_T\in A_T^C\right)
		%
		%
		%
		\leq C\,d_T^{\frac{2+\delta^2}{2}}\, \sum_{t=1}^{T}w_{t,T}^{4+\delta}\,\left[E\left[g\left( \X^\star_{t,T}\right)\right]^{4+\delta}+E\left|f\left(\un{0},\X^\star_{t,T}\right)\right|^{4+\delta}\right]  =\go\left(d_T^{-\frac{\delta-\delta^2}{2}}\right).
	\end{equation*}
\end{proof}\medskip

\begin{lemma}\label{le: BT}\ \\
	Suppose that Assumptions \ref{as: Ge} and \ref{as: Ra} are fulfilled. Then, there exist sets $\left(B_T\right)_{T\in\N}$ with $P\left(\mathbb X_T \in B_T\right)\to 1$ as $T\to\infty$ such that for any sequence $(\un x_T)_T$ with $\un x_T\in B_T$ for all $T$ and for any deterministic sequence $(c_T)_T$ with $c_T\Tinfty \infty$
	\begin{align*}
		&\var\biggl(\sum_{t=1}^{L_T\left\lfloor T/L_T\right\rfloor} w_{t,T} \left(f\left(\s_1, \X_{t,T}^\star\right)-f\left(\s_2, \X_{t,T}^\star\right)\right)\,\bigg|\,\mathbb X_T=\un x_T\biggr)\\
		&\begin{multlined}[t][\linewidth]
			=\sum_{h=-(L_T-1)}^{L_T-1}\sum_{t=L_T\left\lceil\left(TD_T+1\right)/L_T\right\rceil+1}^{L_T\left\lfloor\left(T-TD_T\right)/L_T\right\rfloor}w_{t+h,T}\,w_{t,T}\\
			\cdot\cov\left(f\left(\s_1, \widetilde{\X}_{h}\left(\frac{t}{T}\right)\right)-f\left(\s_2, \widetilde{\X}_{h}\left(\frac{t}{T}\right)\right),f\left(\s_1, \widetilde{\X}_{0}\left(\frac{t}{T}\right)\right)-f\left(\s_2, \widetilde{\X}_{0}\left(\frac{t}{T}\right)\right)\right)\\
			+\go\left( L_T^{-1}   \,+\,c_T\,L_T\left(T^{-\frac{3+\delta}{2(4+\delta)}}\,+\, D_T^{\frac{3+\delta}{2(4+\delta)}}\,+\, (TD_T)^{-1/2}\right) \right)
		\end{multlined}
	\end{align*}
	holds on $B_T$ for $\s_1,\s_2\in\mathcal S$, where both $\s_1$ and $\s_2$ have no influence on the $\go$-term.
\end{lemma}
\begin{proof}[Proof]
	The proof follows the lines of the proof of Lemma~\ref{th: VarBS} aiming to replace the $\go_P$-terms in the former proof by $\go$-terms with slightly modified rates to assure them to hold uniformly on $(B_T)_T$. Let $(c_T)_T$ be a deterministic, monotonically non-decreasing sequence tending to infinity at an arbitrary (slow) rate. Starting from \eqref{eq: T3BS1}, we first argue that the effect of the endpoints is of negligible order.  To this end, define
	\begin{equation*}
		B_{T,1}:=\biggl\{\un x_T\in\R^{dT}\,\bigg|\,[\text{ I+III  in \eqref{eq: T3BS1} with }\var^\star \text{ substituted by }\var^\star_T] \leq d_T^{-\frac{\delta}{2+\delta}}\, c_T\biggr\}.
	\end{equation*}
	From the proof of~\eqref{eq: T3BS1a}, we obtain $P(\mathbb X_T\in B_{T,1})\to 1$ as $T\to \infty$ by Markov's inequality. At this point, we come back to the remaining term~II of~(\ref{eq: T3BS1}) and define II$_T$ as II in~(\ref{eq: T3BS1}) with $\var^\star$ substituted by $\var^\star_T$. Further, we set
	 \begin{multline*}
 		R_T:=\text{II}_T-
 		\sum_{h=-(L_T-1)}^{L_T-1}\sum_{t=L_T\left\lceil\left(TD_T+1\right)/L_T\right\rceil+1}^{L_T\left\lfloor\left(T-TD_T\right)/L_T\right\rfloor}w_{t+h,T}\,w_{t,T}\\
 		\cdot\cov\left(f\left(\s_1, \widetilde{\X}_{h}\left(\frac{t}{T}\right)\right)-f\left(\s_2, \widetilde{\X}_{h}\left(\frac{t}{T}\right)\right),f\left(\s_1, \widetilde{\X}_{0}\left(\frac{t}{T}\right)\right)-f\left(\s_2, \widetilde{\X}_{0}\left(\frac{t}{T}\right)\right)\right)
	 \end{multline*}
and
\begin{equation*}
	  B_{T,2}:=\left\{\un x_T\,\middle |\, \left|R_T\right|\leq C\left( L_T^{-1}+ \left(\frac{L_T}{T}\right)^\frac{\delta}{1+\delta}\,+\,c_T\,L_T\left(T^{-\frac{3+\delta}{2(4+\delta)}}\,+\, D_T^{\frac{3+\delta}{2(4+\delta)}}\,+\, (TD_T)^{-1/2}\right)\right) \right\}.
\end{equation*}
 With the arguments in the proof of Lemma~\ref{th: VarBS} and Markov's inequality, we get 
 \begin{equation*}
 	 P(\mathbb X_T\in B_{T,2})\Tinfty 1.
 \end{equation*}
 Additionally, Assumption~\ref{as: Ra} yields (together with $T\geq d_T$)
 \begin{equation*}
 	\go\left(d_T^{\frac{\delta}{2+\delta}}\right)\subseteq\go\left(L_T^{-1}\right)\quad\text{and}\quad \go\left(\left(\frac{L_T}{T}\right)^\frac{\delta}{1+\delta}\right)\subseteq\go\left(L_T^{-1}\right) .
 \end{equation*}
Thus, the desired result follows from $P(\mathbb X_T\in B_{T,1}\cap B_{T,2})\to 1$ as $T\to \infty$.
\end{proof}\medskip
\begin{lemma}\label{le: K}\ \\
	Provided the validity of Assumptions \ref{as: Ra} and \ref{as: G-BS} (case (a)), there exists a sequence of sets $\left(K_T\right)_{T\in\N}$ with $P(\mathbb X_T \in K_T)\Tinfty 1$ such that for any sequence $(\un x_T)_T$ with $\un x_T\in K_T$ for all $T$ 
	\begin{equation*}
		\sum_{t=0}^{\left\lfloor T/L_T\right\rfloor-1} E\left(\left|L_{t,T}^\star\right|_{\infty}^{1-\delta}\, \left|L_{t,T}^\star\right|_{\mathrm{Lip}}^{1+\delta}\,\middle |\, \mathbb X_t=x_t\right)
		\leq L_T^{1+\delta}
	\end{equation*}
 with~{\neu}$	L_{t,T}^\star\left(\s\right)$ being defined as in~\eqref{eq: Lstern}.
\end{lemma}
\begin{proof}[Proof]
	 As $P(\mathbb X_T\in K_T)\leq L_T^{-1-\delta}\sum_{t=0}^{\left\lfloor T/L_T\right\rfloor-1} E\left[E^\star\left(\left|L_{t,T}^\star\right|_{\infty}^{1-\delta}\, \left|L_{t,T}^\star\right|_{\mathrm{Lip}}^{1+\delta}\right)\right]$,  it is enough to show that the sum on the RHS is of lower order than $(L_T^{1+\delta})_T$.
	An iterative application of both H\"older's and Jensen's inequality  yields
	\begin{align}\label{eq: K1}
		&\sum_{t=0}^{\left\lfloor T/L_T\right\rfloor-1} E\left[E^\star\left(\left|L_{t,T}^\star\right|_{\infty}^{1-\delta}\, \left|L_{t,T}^\star\right|_{\mathrm{Lip}}^{1+\delta}\right)\right]\nonumber\\
		&\leq\sum_{t=0}^{\left\lfloor T/L_T\right\rfloor-1} \left(EE^\star\left|L_{t,T}^\star\right|_{\infty}^{4+\delta}\right)^{\frac{1-\delta}{4+\delta}}\left(EE^\star\left|L_{t,T}^\star\right|_{\mathrm{Lip}}^{\frac{2(1+\delta)(4+\delta)}{2+3\delta}}\right)^{\frac{2+3\delta}{2(4+\delta)}}\nonumber\\
		&\begin{multlined}[t][\linewidth]
			\leq C\sum_{t=0}^{\left\lfloor T/L_T\right\rfloor-1}\left(E\biggl(\sum_{j=1}^{L_T}w_{tL_T+j,T}\left\|g\left(\X_{tL_T+j,T}^\star\right)\right\|_{4+\delta,\star}\biggr)^{4+\delta}\right.\\
			\left.+E\biggl(\sum_{j=1}^{L_T}w_{tL_T+j,T}\left\|f\left(\un 0,\X_{tL_T+j,T}^\star\right)\right\|_{4+\delta,\star}\biggr)^{4+\delta}\right)^{\frac{1-\delta}{4+\delta}}\\
			\cdot\left(E\biggl(\sum_{j=1}^{L_T}w_{tL_T+j,T}\left\|g\left( \X_{tL_T+j,T}^\star\right)\right\|_{\frac{2(1+\delta)(4+\delta)}{2+3\delta},\star}\biggr)^{\frac{2(1+\delta)(4+\delta)}{2+3\delta}}\right)^{\frac{2+3\delta}{2(4+\delta)}}
		\end{multlined}\nonumber\\
		&\begin{multlined}[b][\linewidth]
			\leq C\sum_{t=0}^{\left\lfloor T/L_T\right\rfloor-1}\left(\left(E\biggl(\sum_{j=1}^{L_T}w_{tL_T+j,T}\left\|g\left( \X_{tL_T+j,T}^\star\right)\right\|_{4+\delta,\star}\biggr)^{4+\delta}\right)^{\frac{1-\delta}{4+\delta}}\right.\\
			\left.+\left(E\biggl(\sum_{j=1}^{L_T}w_{tL_T+j,T}\left\|f\left(\un 0,\X_{tL_T+j,T}^\star\right)\right\|_{4+\delta,\star}\biggr)^{4+\delta}\right)^{\frac{1-\delta}{4+\delta}}\right)\\
			\cdot\left(E\biggl(\sum_{j=1}^{L_T}w_{tL_T+j,T}\left\|g\left( \X_{tL_T+j,T}^\star\right)\right\|_{\frac{2(1+\delta)(4+\delta)}{2+3\delta},\star}\biggr)^{\frac{2(1+\delta)(4+\delta)}{2+3\delta}}\right)^{\frac{2+3\delta}{2(4+\delta)}}
		\end{multlined}\nonumber\\
		&=: C\sum_{t=0}^{\left\lfloor T/L_T\right\rfloor-1}\left(\text{I}+\text{II}\right)\cdot\text{III}
\end{align}
	with an obvious definition of I, II, and III. In order to show that the RHS is of order $o(L_T^{1+\delta})$, we investigate the three newly defined terms one by one beginning with the second. From Assumption~\ref{as: Ge}, we obtain
	\begin{equation}\label{eq: K7}
		\text{II}
		%
		%
		%
		%
		%
		%
		\leq \biggl(\sum_{j=1}^{L_T}w_{tL_T+j,T}\left\|\left\|f\left(\un 0,\X_{tL_T+j,T}^\star\right)\right\|_{4+\delta,\star}\right\|_{4+\delta}\biggr)^{1-\delta}\,
\leq \,C\,\biggl(\sum_{j=1}^{L_T}w_{tL_T+j,T}\biggr)^{1-\delta}.
	\end{equation}
	Proceeding with the first subterm on the RHS of (\ref{eq: K1}), we play on Assumption~\ref{as: G-BS} and obtain for $tL_T+j\notin EP$ (and similarly for endpoints)
	\begin{align*}
		&\left\|\left\|g\left(\X_{tL_T+j,T}^\star\right)\right\|_{4+\delta,\star}\right\|_{4+\delta}\nonumber\\
		&\begin{multlined}[t][\linewidth]
			\leq \left\|\left\|g\left(\X_{tL_T+j,T}^\star\right)\right\|_{4+\delta,\star}  
			-\left(\frac{1}{2\,TD_T+1}\sum_{r=-TD_T}^{TD_T}
			\left|g\left( 
			\widetilde{\X}_{tL_T+j+r}\left(\frac{tL_T+j+r}{T}\right)\right)\right|^{4+\delta}\right)^{\frac{1}{4+\delta}}\right\|_{4+\delta}\\
			+\left\|\left(\frac{1}{2\,TD_T+1}\sum_{r=-TD_T}^{TD_T}\left|g\left( \,\widetilde{\X}_{tL_T+j+r}\left(\frac{tL_T+j+r}{T}\right)\right)\right|^{4+\delta}\right)^{\frac{1}{4+\delta}}\right\|_{4+\delta}
		\end{multlined}\\
		&\leq C \left(\frac{1}{T}+1\right).
	\end{align*}
	  This leads to
	\begin{equation}\label{eq: K15}
		\text{I}\,\leq \,C\,\biggl(\sum_{j=1}^{L_T}w_{tL_T+j,T}\biggr)^{1-\delta}.
	\end{equation}
	At this point, only term III of equation (\ref{eq: K1}) is left to be examined. Because said term has the same building type as term I, we can repeat the belonging procedure to get 
	\begin{equation}\label{eq: K16}
		\text{III}\,
		\leq C\,\biggl(\sum_{j=1}^{L_T}w_{tL_T+j,T}\biggr)^{1+\delta}.
	\end{equation}
	Joining equations (\ref{eq: K7}), (\ref{eq: K15}) and (\ref{eq: K16}), we obtain
	\begin{equation*}
		C\sum_{t=0}^{\left\lfloor T/L_T\right\rfloor-1} \biggl(\sum_{j=1}^{L_T}w_{tL_T+j,T} \biggr)^{1-\delta}
		\cdot\biggl(\sum_{j=1}^{L_T}w_{tL_T+j,T} \biggr)^{1+\delta}	\leq C\sum_{t=0}^{\left\lfloor T/L_T\right\rfloor-1}\biggl(\sum_{j=1}^{L_T}w_{tL_T+j,T}\biggr)^2=C\,L_T
	\end{equation*}
	as an upper bound for equation (\ref{eq: K1}), which is of order $o( L_T^{1+\delta}).$
\end{proof}\medskip
\begin{lemma}\label{le: DiffVar}\ \\
	Suppose Assumptions \ref{as: Ge}, \ref{as: Ra} and \ref{as: G-BS} are satisfied. Then, for all $\s_1,\s_2\in\mathcal S$ it holds
	\begin{align}\label{eq: DiffVar1}
		&\begin{multlined}[t][\linewidth]
			\sum_{h=-(L_T-1)}^{L_T-1}\sum_{t=L_T\left\lceil\left(TD_T+1\right)/L_T\right\rceil+1}^{L_T\left\lfloor\left(T-TD_T\right)/L_T\right\rfloor}w_{t+h,T}\,w_{t,T}\\
			\cdot\left|\cov\left(f\left(\s_1, \widetilde{\X}_{h}\left(\frac{t}{T}\right)\right)-f\left(\s_2, \widetilde{\X}_{h}\left(\frac{t}{T}\right)\right),f\left(\s_1, \widetilde{\X}_{0}\left(\frac{t}{T}\right)\right)-f\left(\s_2, \widetilde{\X}_{0}\left(\frac{t}{T}\right)\right)\right)\right|
		\end{multlined}\nonumber\\
		&\leq C_{DC}\left|\s_1-\s_2\right|_1^{1/2}
	\end{align}
	for some positive constant $C_{DC}<\infty$ neither depending on $\s_1,\s_2$ nor on $t$.
\end{lemma}
\begin{proof}
	First, we split the covariance up inserting the truncated version of the companion process with truncation parameter $M:=\left\lceil|h|/2\right\rceil$ like in (\ref{eq: Trun1}) and obtain
	\begin{align*}
		&\left|\cov\left(f\left(\s_1,\widetilde{\X}_0\left(\frac{t}{T}\right)\right)-f\left(\s_2,\widetilde{\X}_0\left(\frac{t}{T}\right)\right),f\left(\s_1,\widetilde{\X}_h\left(\frac{t}{T}\right)\right)-f\left(\s_2,\widetilde{\X}_h\left(\frac{t}{T}\right)\right)\right)\right|\\
		&\begin{multlined}[t][\linewidth]
			\leq\left|\cov\left(f\left(\s_1,\widetilde{\X}_0\left(\frac{t}{T}\right)\right)-f\left(\s_2,\widetilde{\X}_0\left(\frac{t}{T}\right)\right)\right.\right.
			-f\left(\s_1,\widetilde{\X}_0^{(M)}\left(\frac{t}{T}\right)\right)+f\left(\s_2,\widetilde{\X}_0^{(M)}\left(\frac{t}{T}\right)\right),\\
			\left.\left.f\left(\s_1,\widetilde{\X}_h\left(\frac{t}{T}\right)\right)-f\left(\s_2,\widetilde{\X}_h\left(\frac{t}{T}\right)\right)\right)\right|\\
			+\left|\cov\left(f\left(\s_1,\widetilde{\X}_0^{(M)}\left(\frac{t}{T}\right)\right)-f\left(\s_2,\widetilde{\X}_0^{(M)}\left(\frac{t}{T}\right)\right),\right.\right.\\
			f\left(\s_1,\widetilde{\X}_h\left(\frac{t}{T}\right)\right)-f\left(\s_2,\widetilde{\X}_h\left(\frac{t}{T}\right)\right)
			\left.\left.-f\left(\s_1,\widetilde{\X}_h^{(M)}\left(\frac{t}{T}\right)\right)+f\left(\s_2,\widetilde{\X}_h^{(M)}\left(\frac{t}{T}\right)\right)\right)\right|
		\end{multlined}\\
		&=:\text{I}+\text{II}.
	\end{align*}
	For symmetry reasons, we examine only term I. With the use of the Cauchy-Schwarz inequality, we get
	\begin{align}\label{eq: DiffVar4}
		\text{I}
		&\begin{multlined}[t][0.98\linewidth]
			\leq \left(E\left|f\left(\s_1,\widetilde{\X}_0\left(\frac{t}{T}\right)\right)-f\left(\s_2,\widetilde{\X}_0\left(\frac{t}{T}\right)\right)
			-\left(f\left(\s_1,\widetilde{\X}_0^{(M)}\left(\frac{t}{T}\right)\right)-f\left(\s_2,\widetilde{\X}_0^{(M)}\left(\frac{t}{T}\right)\right)\right)\right|^2\right)^{1/2}\\
			\cdot\left(E\left|f\left(\s_1,\widetilde{\X}_h\left(\frac{t}{T}\right)\right)-f\left(\s_2,\widetilde{\X}_h\left(\frac{t}{T}\right)\right)\right|^2\right)^{1/2}
		\end{multlined}\nonumber\\
		&=:\text{Ia}\cdot\text{Ib}.
	\end{align}
	Continuing again with the first factor,  we apply H\"older's inequality iteratively  and similar arguments as in the proof of Lemma~\ref{le: DiffProd} to obtain
	\begin{align*}
		\text{Ia}
		&\begin{multlined}[t][0.96\linewidth]
			\leq\left(E\left|f\left(\s_1,\widetilde{\X}_0\left(\frac{t}{T}\right)\right)-f\left(\s_2,\widetilde{\X}_0\left(\frac{t}{T}\right)\right)\right.\right.\\
			\left.\left.-\left(f\left(\s_1,\widetilde{\X}_0^{(M)}\left(\frac{t}{T}\right)\right)-f\left(\s_2,\widetilde{\X}_0^{(M)}\left(\frac{t}{T}\right)\right)\right)\right|^\frac{2+\delta}{1+\delta}\right)^{\frac{1+\delta}{2+\delta}\cdot\frac{1}{2}}\\
			\cdot\left(E\left|f\left(\s_1,\widetilde{\X}_0\left(\frac{t}{T}\right)\right)-f\left(\s_2,\widetilde{\X}_0\left(\frac{t}{T}\right)\right)\right.\right.\\
			\left.\left.-\left(f\left(\s_1,\widetilde{\X}_0^{(M)}\left(\frac{t}{T}\right)\right)-f\left(\s_2,\widetilde{\X}_0^{(M)}\left(\frac{t}{T}\right)\right)\right)\right|^{2+\delta}\right)^{\frac{1}{2+\delta}\cdot\frac{1}{2}}
		\end{multlined}\\
		&\leq C\, \rho^{\frac{M\delta}{2(1+\delta)}}.
	\end{align*}
	For Ib in equation (\ref{eq: DiffVar4}), we use H\"older's inequality anew and get in case (a) of Assumption~\ref{as: G-BS}
	\begin{align*}
		\text{Ib}
		&\begin{multlined}[t][0.96\linewidth]
			\leq C\, \left(E\left|f\left(\s_1,\widetilde{\X}_h\left(\frac{t}{T}\right)\right)-f\left(\s_2,\widetilde{\X}_h\left(\frac{t}{T}\right)\right)\right|^\frac{2+\delta}{1+\delta}\right)^{\frac{1+\delta}{2+\delta}\cdot\frac{1}{2}}\leq C\left|\s_1-\s_2\right|_1^{1/2}.
		\end{multlined}
	\end{align*}
	Part (b) can be treated in a similar manner using boundedness of $f$ instead of  H\"older's inequality.
	Consequently, we have
	\begin{equation*}
		\text{I}\leq C\,\rho^{\frac{M\delta}{2(1+\delta)}}\left|\s_1-\s_2\right|_1^{1/2}.
	\end{equation*}
	Hence, we can bound (\ref{eq: DiffVar1}) by
	\begin{equation*}
		C\,	\sum_{h=-(L_T-1)}^{L_T-1}\sum_{t=L_T\left\lceil\left(TD_T+1\right)/L_T\right\rceil+1}^{L_T\left\lfloor\left(T-TD_T\right)/L_T\right\rfloor}w_{t+h,T}\,w_{t,T}\,\rho^{\frac{M\delta}{2(1+\delta)}} \left|\s_1-\s_2\right|_1^{1/2}
		\leq C_{DC}\left|\s_1-\s_2\right|_1^{1/2}.
	\end{equation*}
	%
\end{proof}\medskip
\begin{lemma}\label{le: KT}\ \\
	Under Assumptions \ref{as: Ra} and \ref{as: G-BS} (case (b)),
		there exists a sequence of sets $\left(\bar K_T\right)_{T\in\N}$ with $P(\mathbb X_T \in \bar K_T)\Tinfty 1$, such that for any sequence $(\un x_T)_T$ with $\un x_T\in\bar K_T$ for all $T$ 
		\begin{equation*}
			\sum_{t=0}^{\left\lfloor T/L_T\right\rfloor-1} E\left(\left|L_{t,T}^\star\right|_{\infty}^{\frac{2-\delta}{2}}\, \left|L_{t,T}^\star\right|_{\mathrm{Lip}}^{\frac{2+\delta}{2}}\,\middle |\, \mathbb X_t=x_t\right)
			\leq C\, L_T^{1+\delta}
		\end{equation*}
		holds with $L_{t,T}^\star$ being defined in \eqref{eq: Lstern}.
	%
	%
\end{lemma}
\begin{proof}[Proof]
	This proof models itself on the proof of Lemma~\ref{le: K}, i.e.~we show that 
		\begin{equation*}
			E\biggl[	\sum_{t=0}^{\left\lfloor T/L_T\right\rfloor-1} E^\star\left(\left|L_{t,T}^\star\right|_{\infty}^{\frac{2-\delta}{2}}\, \left|L_{t,T}^\star\right|_{\mathrm{Lip}}^{\frac{2+\delta}{2}}\,\right)\biggr]=O(L_T).
		\end{equation*}
		First, note that
		\begin{equation}\label{eq: BS-T15B}
			\left|L_{t,T}^\star\right|_{\infty}
			\leq C\,\sum_{j=1}^{L_T}w_{tL_T+j,T}
		\end{equation}
		and
		\begin{align}\label{eq: BS-T16B}
			\left|L_{t,T}^\star\right|_{\mathrm{Lip}}
			%
			%
			%
			%
			&=\sum_{j=1}^{L_T}w_{tL_T+j,T}\,g\left(\X_{tL_T+j,T}^\star\right).
		\end{align}
		Hence, we have
		\begin{equation*}
			E^\star\biggl( \left|L_{t,T}^\star\right|_{\mathrm{Lip}}\biggr)^{\frac{2+\delta}{2}}
			\leq \biggl(\sum_{j=1}^{L_T}w_{tL_T+j,T}\left\|g\left(\X_{tL_T+j,T}^\star\right)\right\|_{\frac{2+\delta}{2},\star}\biggr)^{\frac{2+\delta}{2}}.
		\end{equation*}
		For sake of notational simplicity, we consider non-endpoints only in the sequel. We get 
		\begin{align}\label{eq: KT2b}
			&\left\|\left\|g\left(\X_{tL_T+j,T}^\star\right)\right\|_{\frac{2+\delta}{2},\star}\right\|_{\frac{2+\delta}{2}}\nonumber\\
			&\begin{multlined}[t][\linewidth]
				\leq \left\|\left\|g\left(\X_{tL_T+j,T}^\star\right)\right\|_{\frac{2+\delta}{2},\star}-\biggl(\frac{1}{2\,TD_T+1}\sum_{r=-TD_T}^{TD_T}\left|g\left(\widetilde{\X}_{tL_T+j+r}\left(\frac{tL_T+j+r}{T}\right)\right)\right|^{\frac{2+\delta}{2}}\biggr)^{\frac{2}{2+\delta}}\right\|_{\frac{2+\delta}{2}}\\
				+\left\|\biggl(\frac{1}{2\,TD_T+1}\sum_{r=-TD_T}^{TD_T}\left|g\left(\widetilde{\X}_{tL_T+j+r}\left(\frac{tL_T+j+r}{T}\right)\right)\right|^{\frac{2+\delta}{2}}\biggr)^{\frac{2}{2+\delta}}\right\|_{\frac{2+\delta}{2}}
				%
				%
				%
				%
				%
			\end{multlined}\nonumber\\
			&=:\text{I}+\text{II}.
		\end{align}
		Term~I can be bounded by $C/T$ similarly to the respective term in the aforementioned proof. Straightforward arguments yield II$=O(1)$. Taken all together, from \eqref{eq: BS-T15B}, \eqref{eq: BS-T16B} and \eqref{eq: KT2b} we obtain
		\begin{equation*}
			E\biggl[\sum_{t=0}^{\left\lfloor T/L_T\right\rfloor-1} E^\star\left(\left|L_{t,T}^\star\right|_{\infty}^{\frac{2-\delta}{2}}\, \left|L_{t,T}^\star\right|_{\mathrm{Lip}}^{\frac{2+\delta}{2}} \right)\biggr]
			\leq C\,
			\sum_{t=0}^{\left\lfloor T/L_T\right\rfloor-1} \biggl(\sum_{j=1}^{L_T}w_{tL_T+j,T}\biggr)^{\frac{2-\delta}{2}}\,\biggl(\sum_{j=1}^{L_T}w_{tL_T+j,T}\biggr)^{\frac{2+\delta}{2}}
			\leq \, C\, L_T.
				\end{equation*}
\end{proof}\medskip
 \begin{lemma}\label{le: BS-FCLT-b}
		Suppose that the set of assumptions of Theorem~\ref{th: FCLTBS} holds with part (b) in Assumption~\ref{as: Fu2} and~\ref{as: G-BS}. Then, there exist sets $\left(\bar\Omega_T\right)_{T\in\N}$ with $P\left(\mathbb X_T\in\bar \Omega_T\right)\to 1$ as $T\to\infty$ such that for any $(\un x_T)_T$ with $x_T\in\bar\Omega_T$ for all $T$ 
		\begin{equation*}
			\underset{r\to 0}{\lim}\,\limsup_{T\to\infty}\, P\biggl(\sup_{\rho\left(\s_1,\s_2\right)<r}\biggl|\sum_{t=1}^T w_{t,T}\left(\f^\star\left(\s_1, \X^\star_{t,T}\right)-\f^\star\left(\s_2, \X^\star_{t,T}\right)\right)\biggr|>\lambda\,\bigg|\,\mathbb X_T=\un x_T\biggr)=0
		\end{equation*}
		with $\lambda>0$, whereas $\f^\star$ is defined in (\ref{eq: CentBS}).
\end{lemma}

\begin{proof}
	First, we define subsets $\left(\bar\Omega_T\right)_{T\in\N}$ of $\Omega$ as 
		\begin{equation*}
		\bar\Omega_T=A_T\cap B_T\cap \bar K_T
		\end{equation*}
	with $A_T$ as in Lemma~\ref{le: UB},  $B_T$  as in Lemma~\ref{le: BT}  and~$\bar{K}_T$ being defined in  Lemma~\ref{le: KT}.
	Then, we obtain $\lim_{T\to\infty}P\left(\bar{\Omega}_T\right)=1$.	
	As in the proof of part (a), we split up 
		\begin{align}\label{eq: BS-T1B}
			&P^\star_T\Biggl(\sup_{ {\rho\left(\s_1,\s_2\right)<r}}\biggl|\sum_{t=1}^T w_{t,T}\left(\f^\star\left(\s_1, \X^\star_{t,T}\right)-\f^\star\left(\s_2, \X^\star_{t,T}\right)\right)\biggr|>\lambda \Biggr)\nonumber\\
			&\begin{multlined}[t][\linewidth]
				\leq P^\star_T\Biggl(\sup_{ {\rho\left(\s_1,\s_2\right)<r}}\biggl|\sum_{t=1}^{L_T \left\lfloor T/L_T\right\rfloor} w_{t,T}\left(\f^\star\left(\s_1, \X^\star_{t,T}\right)-\f^\star\left(\s_2, \X^\star_{t,T}\right)\right)\biggr|>\frac{\lambda}{2} \Biggr)\\
				+\,P^\star_T\Biggl(\sup_{{\rho\left(\s_1,\s_2\right)<r}}\biggl|\sum_{t=L_T \left\lfloor T/L_T\right\rfloor+1}^{T} w_{t,T}\left(\f^\star\left(\s_1, \X^\star_{t,T}\right)-\f^\star\left(\s_2, \X^\star_{t,T}\right)\right)\biggr|>\frac{\lambda}{2} \Biggr).
			\end{multlined}		
		\end{align}
		With the use of Markov's inequality, asymptotic negligibility of the  second sum of the RHS of~(\ref{eq: BS-T1B}) can be verified since it can be bounded by
		\begin{align*}
			\frac{2}{\lambda}\,E^\star_T\Biggl(\sup_{{\rho\left(\s_1,\s_2\right)<r}}\biggl|\sum_{t=L_T \left\lfloor T/L_T\right\rfloor+1}^{T} w_{t,T}\left(\f^\star\left(\s_1, \X^\star_{t,T}\right)-\f^\star\left(\s_2, \X^\star_{t,T}\right)\right)\biggr|\Biggr)
			&=o\left(d_T^{-\frac{1}{2(1+\delta)}}\right).
		\end{align*}
		It remains to consider
		\begin{align}\label{eq: BS-T1aB}
			&P^\star_T\Biggl(\sup_{ {\rho\left(\s_1,\s_2\right)<r}}\biggl|\sum_{t=1}^{L_T \left\lfloor T/L_T\right\rfloor} w_{t,T}\left(\f^\star\left(\s_1, \X^\star_{t,T}\right)-\f^\star\left(\s_2, \X^\star_{t,T}\right)\right)\biggr|>\frac{\lambda}{2} \Biggr)\nonumber\\
			&\begin{multlined}[t][\linewidth]
				\leq P^\star_T\Biggl(2\sup_{\overset{\s_1,\s_2\in\mathcal S}{\rho\left(\s_1,\s_2\right)\leq r_{k_T}}}\left|\nu_T^\star\left(\s_1,\s_2\right)\right|>\frac{\lambda}{6} \Biggr)
				+ P^\star_T\Biggl(\sup_{\overset{\s_1,\s_2\in\mathcal F_0}{\rho\left(\s_1,\s_2\right)\leq 3\,r}}\left|\nu_T^\star\left(\s_1,\s_2\right)\right|>\frac{\lambda}{6} \Biggr)\\
				+P^\star_T\Biggl(2\sum_{k=1}^{k_T}\sup_{\overset{\s_1 \in \mathcal F_k, \s_2\in\mathcal F_{k-1}}{\rho\left(\s_1,\s_2\right)\leq 3\,r_k}}\left|\nu_T^\star\left(\s_1,\s_2\right)\right|>\frac{\lambda}{6} \Biggr)
			\end{multlined}\nonumber\\
			&=:\text{I}+\text{II}+\text{III}.
	\end{align}
	While I and II can be treated as in the proof of Theorem~\ref{th: FCLTBS},
	we have to adapt the investigation of term~III. Using the notation of the proof of Theorem~\ref{th: FCLTBS}, it is enough to verify asymptotic negligibility of 
		\begin{equation}\label{eq: BS-T12aB}
			\,E^\star_T\Biggl(\sup_{\overset{\s_1,\s_2\in\mathcal S}{\rho\left(\s_1,\s_2\right)\leq r_{k_T}}}\biggl|\sum_{t=0}^{\left\lfloor T/L_T\right\rfloor-1}\left(L_{t,T}^{\star,0}\left(\s_1\right)-L_{t,T}^{\star,0}\left(\s_2\right)\right)\biggr|\Biggr).
		\end{equation}
		To this end, first recall that $\sum_{t=0}^{\left\lfloor T/L_T\right\rfloor-1}L_{t,T}^{\star,0}$ comes with sub-Gaussian increments w.r.t.~$\hat\rho_{T,2}$ defined in \eqref{eq: BS-T5} conditionally on $L_{0,T}^\star,\dots,L_{\left\lfloor T/L_T\right\rfloor-1,T}^\star$.
		%
		%
		Second, it holds
		\begin{equation*}
			\left(L_{t,T}^\star\left(\s_1\right)-L_{t,T}^\star\left(\s_2\right)\right)^2 
			\leq 2^{1-\delta} \left|L_{t,T}^\star\right|_{\infty}^{\frac{2-\delta}{2}} \;\left|L_{t,T}^\star\right|_{\mathrm{Lip}}^{\frac{2+\delta}{2}}\rho\left(\s_1,\s_2\right)^{\frac{2+\delta}{2}},\quad \s_1,\s_2\in\mathcal S,
		\end{equation*}
		as seen comparably in (\ref{normenBS}). The next step is again to establish a semimetric which suits us more than $\rho$ and $\hat\rho_{T,2}$, respectively. To this end, let 
	\begin{equation*}
		Q_T:= 2^\frac{2-\delta}{4}\biggl(\sum_{t=0}^{\left\lfloor T/L_T\right\rfloor-1} \left|L_{t,T}^\star\right|_{\infty}^{\frac{2-\delta}{2}}\, \left|L_{t,T}^\star\right|_{\mathrm{Lip}}^{\frac{2+\delta}{2}}\biggr)^{1/2}
	\end{equation*}
	%
	and define
	\begin{equation*}
		\hat{\rho}_{T,2}\left(\s_1,\s_2\right)\leq Q_T\rho\left(\s_1,\s_2\right)^{\frac{2+\delta}{4}}=:\breve{\rho}_T\left(\s_1,\s_2\right).
	\end{equation*}
	%
	%
	Returning to (\ref{eq: BS-T12aB}), we obtain from the maximal inequality for sub-Gaussian processes in Corollary~2.2.8 of \citet{VW00} that
	\begin{align*}
	&	E^\star_T\left[E^\star\Biggl(\sup_{{\rho\left(\s_1,\s_2\right)\leq r_{k_T}}}\biggl|\sum_{t=0}^{\left\lfloor T/L_T\right\rfloor-1}\left(L_{t,T}^{\star,0}\left(\s_1\right)-L_{t,T}^{\star,0}\left(\s_2\right)\right)\biggr|\,\Bigg|\,L_{1,T}^\star,\dots,L_{L_T \left\lfloor T/L_T\right\rfloor,T}^\star\Biggr)\right]\\
		&\leq \,C\, E^\star\Biggl(\int_{0}^{Q_T r_{k_T}^{\frac{2+\delta}{4}}}\left(\log D\left(u,\mathcal S,\breve{\rho}_T\right)\right)^{1/2}du\Biggr)\\
			&\leq C\, E^\star_T\Biggl[\int_{0}^{Q_T r_{k_T}^{\frac{2+\delta}{4}}}\Biggl(\log \Biggl(\left(\frac{u}{Q_T}\right)^{-\frac{4}{2+\delta}}+1\Biggr)^d\Biggr)^{1/2}du\Biggr]\\
		&=C \,E^\star_T  [Q_T]\,\int_{0}^{r_{k_T}^{\frac{2+\delta}{4}}}\left(\log \left({u^{-\frac{4}{2+\delta}}}+1\right)^d\right)^{1/2}du\\
		&\leq\, C\,  L_T^{\frac{1+\delta}{2}}\,L_T^{-\frac{1+\delta}{4\delta}}\\
		&=\,C\, L_T^{-\frac{1-\delta-\delta^2}{4\delta}}\Tinfty 0,
	\end{align*}
 where the last inequality follows from Lemma~\ref{le: KT} and
 \begin{equation*}
 		\int_{0}^{r_{k_T}^{\frac{2+\delta}{4}}}\left(\log \left({u^{\frac{4}{2+\delta}}}+1\right)^d\right)^{1/2}du
 	\leq C_2\,L_T^{-\frac{1+\delta}{4\delta}}.
 \end{equation*}
Hence, applying Markov's inequality we have proven asymptotic negligibility of III in \eqref{eq: BS-T1aB}.
\end{proof}

\end{document}